\documentclass[11pt,reqno]{amsart}

\usepackage{amsmath,amssymb,amsthm,mathtools}
\numberwithin{equation}{section}
\usepackage{enumitem}
\usepackage{graphicx}
\usepackage[margin=1in]{geometry}
\usepackage[
    hidelinks,
    pdftitle={Long-Time Behavior of Zero-Noise Mean-Field Adam},
    pdfauthor={Ethan Huffman and Asaf Cohen}
]{hyperref}
\usepackage[capitalize,noabbrev]{cleveref}

\newcommand{\RR}{\mathbb{R}}
\newcommand{\NN}{\mathbb{N}}
\newcommand{\SCP}{\mathcal{P}}
\newcommand{\WW}{\mathcal{W}}
\newcommand{\EE}{\mathbb{E}}
\newcommand{\FF}{\mathcal{F}}
\newcommand{\Law}{\operatorname{Law}}
\newcommand{\Lip}{\operatorname{Lip}}
\newcommand{\defeq}{\vcentcolon=}
\newcommand{\from}{\colon}
\newcommand{\dd}{\,\mathrm{d}}

\newtheorem{thm}{Theorem}[section]
\newtheorem{lem}[thm]{Lemma}
\newtheorem{prop}[thm]{Proposition}
\newtheorem{cor}[thm]{Corollary}

\theoremstyle{definition}
\newtheorem{asmp}[thm]{Assumption}
\newtheorem{defn}[thm]{Definition}
\newtheorem{ex}[thm]{Example}

\crefname{thm}{Theorem}{Theorems}
\crefname{lem}{Lemma}{Lemmas}
\crefname{prop}{Proposition}{Propositions}
\crefname{cor}{Corollary}{Corollaries}
\crefname{asmp}{Assumption}{Assumptions}
\crefname{defn}{Definition}{Definitions}
\crefname{ex}{Example}{Examples}

\title{Long-Time Behavior of Zero-Noise Mean-Field Adam
}

\author{Ethan Huffman}
\address{Department of Mathematics, University of Michigan,
530 Church Street, Ann Arbor, MI 48109-1043, USA}
\email{ehuf@umich.edu}

\author{Asaf Cohen}
\address{Department of Mathematics, University of Michigan,
530 Church Street, Ann Arbor, MI 48109-1043, USA}
\email{shloshim@gmail.com}

\date{}

\begin{document}

\begin{abstract}
Adaptive methods such as Adam are the most commonly used methods for training neural networks, but their long-time behavior is not as well understood as that of non-adaptive methods.
This is particularly true in the mean-field limit, where the dynamics of infinitely many neurons are described as a distribution over parameters which satisfies a self-consistency condition.
In this paper, we study a deterministic continuous-time model of Adam in the mean-field setting.
Under a condition relating the averaging rates of the momentum and second-moment variables, the size of the gradients, and a regularization parameter, we prove that the distribution of parameters converges to the set of equilibria.
This condition extends the range known for finite-dimensional models, and we show that it is sharp.
Under an additional \L{}ojasiewicz-type inequality for distributions, we prove that individual parameter trajectories converge and obtain quantitative rates for this convergence.
We verify the qualitative convergence assumptions for a regularized two-layer neural network.
We also show that adding noise to the momentum prevents the system from approaching these equilibria and rules out many free-energy type Lyapunov functionals.
Finally, numerical experiments illustrate the energy estimates, give evidence of persistent oscillations outside the sufficient parameter range, and compare the observed decay with the predicted rates.

\end{abstract}

\maketitle

\section{Introduction}
\label{sec:introduction}

Adam \cite{KingmaBa2015} is a widely used method for training neural networks which maintains two exponential moving averages of past gradients.
The first acts as momentum, while the second rescales the step in each coordinate, in the spirit of AdaGrad \cite{DuchiEtAl2011}.
This enables Adam to adapt to the local geometry of the objective.
Despite the widespread use of Adam, its long-time behavior is understood far less completely than that of gradient descent.
Analyses of the discrete algorithm show that Adam can fail to converge \cite{ReddiEtAl2018} and, under suitable assumptions, bound the size of the gradients it produces \cite{ChenEtAl2019,DefossezEtAl2022}.
Continuous-time models have used energy methods to describe how a single parameter vector evolves over long time horizons \cite{BelottoGazeau2020,BarakatBianchi2021}.
At any fixed network width, the parameter vector can contain all neurons and their memory variables.
However, convergence at each fixed width does not by itself give rates that remain bounded away from zero as the width grows.
At infinite width, the natural object is the distribution of parameters across neurons, and this perspective underlies the mean-field theory of training by gradient descent \cite{MeiMontanariNguyen2018,ChizatBach2018,RotskoffVandenEijnden2022,SirignanoSpiliopoulos2020}.
For Adam, a long-time analysis in this setting must control the joint law of the parameters and both memory variables, while the force depends on the evolving parameter distribution.

In this paper, we carry out such an analysis for a deterministic continuous-time model of Adam acting on a probability distribution of parameters.
We show that the joint law of the parameters and their momentum and second-moment variables approaches the set of equilibria and, under a \L{}ojasiewicz-type inequality, that individual trajectories converge at quantitative rates.

\subsection{From Adam to a mean-field differential equation}

Let $x_k$ be a parameter vector and let $g_k$ be the estimate of the gradient at iteration $k$.
With learning rate $h>0$ and averaging coefficients $\beta_1,\beta_2\in[0,1)$, Adam updates according to
\begin{align*}
	m_{k+1}&=\beta_1m_k+(1-\beta_1)g_k,\\
	v_{k+1}&=\beta_2v_k+(1-\beta_2)g_k^2,\\
	x_{k+1}&=x_k-h\frac{\widehat m_{k+1}}{\sqrt{\widehat v_{k+1}}+\varepsilon_{\mathrm A}},
\end{align*}
where all vector operations are coordinatewise and $\varepsilon_{\mathrm A}>0$ prevents division by zero.
For the standard initialization $m_0=v_0=0$, the quantities $\widehat m_{k+1}=m_{k+1}/(1-\beta_1^{k+1})$ and $\widehat v_{k+1}=v_{k+1}/(1-\beta_2^{k+1})$ correct the initial bias of the averages toward zero.
Due to the presence of momentum, the objective can increase, even without gradient noise.

For a two-layer neural network with unit parameters $x_1,\ldots,x_N\in\RR^d$, the output at an input $\xi$ is
\begin{equation*}
	f_\nu(\xi)=\int_{\RR^d}\phi(x;\xi)\,\nu(\dd x),
	\qquad
	\nu=\frac1N\sum_{j=1}^N\delta_{x_j},
\end{equation*}
where $\phi(x;\xi)$ is the contribution of a neuron with parameters $x$.
Allowing $\nu$ to be a general probability measure, we are then able to describe a continuum of neurons and study their collective behavior.
The training objective becomes a functional $Q(\nu)$, and the gradient driving a unit at position $x$ becomes a \textit{force} $G(\nu,x)$ that depends on both $x$ and the distribution $\nu$ of all units.

For Adam, we attach a momentum and a second-moment variable to each parameter and track the joint law of the three variables.
Writing $X_t$, $M_t$, and $V_t$ for the position, momentum, and second-moment variable, $\mu_t=\Law(X_t,M_t,V_t)$ for their joint law, and $\mu_t^x=\Law(X_t)$ for its position marginal, the system studied in this paper is
\begin{equation*}
	\begin{aligned}
		\dot X_t&=-\frac{M_t}{\sqrt{V_t+\varepsilon}},\\
		\dot M_t&=\rho_1\bigl(G(\mu_t^x,X_t)-M_t\bigr),\\
		\dot V_t&=\rho_2\bigl(G(\mu_t^x,X_t)^2-V_t\bigr),
	\end{aligned}
\end{equation*}
where $\rho_1,\rho_2>0$ are averaging rates and $\varepsilon>0$ regularizes the denominator.
We formulate it precisely as \eqref{eq:adam-system} in \cref{sec:main}.
The system has \textit{zero noise} in the sense that the exact force drives the evolution, and randomness enters only through the initial law $\mu_0$.

This system is connected to the discrete-time Adam algorithm by the scaling $\beta_i=1-\rho_ih$ with $h\to0$.
Under this scaling, the momentum update is $(m_{k+1}-m_k)/h=\rho_1(g_k-m_k)$, and the second-moment update has the same form.
Scaling in this manner preserves a fixed time scale for the memory of both averaged variables, allowing them to remain dynamic in the limit.
This is in contrast to the regime with fixed $\beta_i$, in which the averages adjust quickly relative to the parameters \cite{MaWuE2022}.
We place the regularizer inside the square root, using $\sqrt{V_t+\varepsilon}$, as in Belotto da Silva and Gazeau \cite[Section~4.2]{BelottoGazeau2020}.
This makes the preconditioning factor smooth even when $V_t=0$.
This regularization is also available in the Optax library for JAX through its \texttt{eps\_root} parameter \cite{deepmind2020jax}.
For ease of modeling and because our focus is on long-time behavior, we omit the debiasing step, whose multiplicative correction factors tend to one as time increases.
In \cref{sec:numerics-setup}, we give rescalings under which the discretization \eqref{eq:discrete-adam} is a first-order consistent scheme for the system above.

\subsection{Main contributions}

Our analysis is built on the energy
\begin{equation*}
	\FF(\mu_t)=Q(\mu_t^x)+\frac{1}{2\rho_1}\EE\sum_{i=1}^d\frac{M_{t,i}^2}{\sqrt{V_{t,i}+\varepsilon}},
\end{equation*}
which adds a weighted momentum term to the objective.
This term accounts for the exchange of energy between the objective and the momentum.
In the finite-dimensional case, an energy of this form was introduced by Belotto da Silva and Gazeau \cite{BelottoGazeau2020} and shown to be nonincreasing when $\rho_2\leq4\rho_1$.
Our contributions are as follows.
\begin{enumerate}[label=\textup{(\roman*)}]
	\item \textit{Convergence of the joint law to equilibrium.}
	\Cref{thm:main} shows that every weak limit of $\mu_t$ as $t\to\infty$ has zero momentum, zero second-moment variable, and a position law $\nu$ satisfying $G(\nu,x)=0$ for $\nu$-almost every $x$.
	We call such a law a \textit{critical law}.
	Additionally, we demonstrate that the weighted momentum and the mean-square force tend to zero, and that the objective converges to a finite value.
	If the set of critical laws is totally disconnected, for example finite, then $\mu_t$ converges weakly to a single equilibrium.
	Notably, this result requires no finite moment of the initial second-moment variable.
	\Cref{cor:sublevel-convergence} verifies the necessary trajectory bounds using an objective sublevel set.
	\item \textit{An enlarged and sharp parameter range.}
	Denoting the bound on the force by $K$, with $K>0$, the theorem holds under the condition
	\begin{equation*}
		\rho_2\leq4\rho_1\left(1+\frac{\varepsilon}{K^2}\right),
	\end{equation*}
	which extends the usual range $\rho_2\leq4\rho_1$.
	The bounded force drives the second-moment variable toward the region $v_i\leq K^2$, where regularization permits the larger range of averaging rates.
	Any initial excess above this region results in an integrable error, which we control through a corrected energy.
	This condition is sharp as a uniform convergence guarantee over smooth confining potentials with coordinatewise gradient bound $K$.
    For every choice of parameters outside this range, \cref{prop:sharpness} constructs, already in dimension three, a law-independent force and a nonconstant periodic solution from suitable initial memory.
    Thus the converse is new even in the finite-dimensional deterministic setting, complementing the convergence analyses of \cite{BelottoGazeau2020,BarakatBianchi2021} by identifying a matching obstruction for the autonomous model studied here.
	\item \textit{Trajectory convergence and rates.}
	Under a weighted mean-field \L{}ojasiewicz inequality and suitable dissipation and regularity assumptions, \Cref{thm:kl-rates} shows that the position trajectories have finite length and that $(X_t,M_t,V_t)$ converges almost surely to $(X_\infty,0,0)$.
	This gives convergence to a single critical law even when the critical laws form a connected set.
	Under a power-type inequality we show that the objective gap, the weighted momentum and force, and the position error all obey exponential or algebraic bounds.
\end{enumerate}

An interesting and important feature of the mean-field setting is that critical laws must be self-consistent.
This allows several collective equilibria even when, for each fixed distribution, the force has a unique zero.
\Cref{ex:three-critical-laws} gives such a model with three critical laws.

We apply the qualitative convergence theorem to a regularized two-layer neural network in \cref{cor:neural-convergence}, with the standard initialization $M_0=V_0=0$.
The examples in \cref{sec:width-dependence} explain why a convergence rate at each finite particle count need not give a rate uniform in that count.
For a limiting interaction model, the slow local exponential rate decreases like $1/N$ as the number of particles $N$ grows.
In contrast, we give an example with a quartic potential which admits a mean-field \L{}ojasiewicz inequality with a constant independent of $N$.

In \cref{sec:noise}, we show that Brownian forcing of the momentum injects energy and keeps the time-averaged mean-square momentum bounded away from zero, so the noisy system does not approach the equilibria above.
For a class of one-dimensional potentials, a symmetry argument rules out convex, reflection-invariant strict Lyapunov functionals with a finite minimum, as stated in \cref{cor:no-free-energy}.
Experiments in \cref{sec:numerics} test the energy estimates and rate predictions, and give numerical evidence of persistent oscillations outside the sufficient averaging-rate range.

Our proofs proceed in three stages.
First, we use variation-of-constants formulas for the momentum and second-moment equations to obtain a priori bounds, and a fixed-point argument gives well-posedness.
Second, we compute the time derivative of the energy exactly.
When $\rho_2>4\rho_1$, large values of the second-moment variable can make this derivative positive.
Because the force bounds the input to the second-moment equation by $K^2$, any excess above this level is inherited from the initial state and decays exponentially.
\Cref{prop:corrected-energy} isolates this excess in an integrable error and shows that a corrected dissipation has finite integral.
Our weighted estimates retain the initial adaptive state in the denominator and give an integrable error without a moment bound on $V_0$.
Third, we translate the trajectory in time and use tightness to extract limiting trajectories defined for all real times, along which the corrected dissipation vanishes.
A rigidity argument shows that every such trajectory is an equilibrium.
At the boundary of the averaging-rate condition, zero dissipation does not immediately force the momentum to vanish.
Instead, the equations along the complete limiting trajectory give the additional restriction needed to identify an equilibrium.
This step is an analogue of LaSalle's invariance principle \cite{LaSalle1960} for a law-dependent second-order system.
For the convergence rates, a further correction of the energy controls both the momentum and the force, and the \L{}ojasiewicz inequality then yields finite length and a differential inequality for the corrected energy.
This correction connects the inequality on probability laws to dissipation in both the momentum and the adaptively weighted force, including times when the objective lies below its limiting value.

\subsection{Related work}

\textit{Analyses of discrete Adam.}
Reddi et~al.\ \cite{ReddiEtAl2018} give examples in online convex optimization for which Adam fails to converge and introduce AMSGrad to address this failure.
Chen et~al.\ \cite{ChenEtAl2019} and D\'efossez et~al.\ \cite{DefossezEtAl2022} establish nonconvex convergence bounds for Adam-type methods under assumptions such as bounded gradients and smoothness.
These bounds control stationarity measures along the iterates and depend explicitly on the algorithmic parameters.
Our results instead concern the small-step regime $\beta_i=1-\rho_ih$, in which we identify the limiting equilibria and prove convergence of the joint law and of individual trajectories.

\textit{Continuous-time models of Adam.}
The closest deterministic model to ours is the general system of Belotto da Silva and Gazeau \cite{BelottoGazeau2020}, which includes the finite-dimensional version of our model.
Barakat and Bianchi \cite{BarakatBianchi2021} connect stochastic Adam to a bias-corrected differential equation that depends explicitly on time.
They prove convergence toward the set of critical points and, under a \L{}ojasiewicz assumption, convergence to a single point with rates.
In their model, the adaptive variable averages the second moment of the stochastic gradient and can remain positive where the expected gradient vanishes, whereas our adaptive variable vanishes at equilibrium.
Relative to these works, the new difficulties here are the dependence of the force on the distribution, the need for compactness arguments in a space of probability measures, convergence of a joint law whose initial adaptive state may have heavy tails,  the extension of the rate condition beyond $\rho_2\leq4\rho_1$, and the identification of equilibria that must be self-consistent.
For our rate theorem, we also must formulate the \L{}ojasiewicz inequality on probability distributions and weight the force by the adaptive state.

\textit{Mean-field training dynamics.}
The distributional description of two-layer networks trained by gradient descent was developed by Mei et~al.\ \cite{MeiMontanariNguyen2018}, Chizat and Bach \cite{ChizatBach2018}, Rotskoff and Vanden-Eijnden \cite{RotskoffVandenEijnden2022}, and Sirignano and Spiliopoulos \cite{SirignanoSpiliopoulos2020}.
Chizat and Bach also prove global convergence for homogeneous models under conditions on the initial support.
We add momentum and coordinatewise adaptation to this setting, as is the case for almost all modern neural networks.
The limiting objects in our theorems are critical laws, and their global optimality requires additional structure of the objective, as in the work of Chizat and Bach.

\textit{Momentum and adaptive methods at large width.}
Krichene et~al.\ \cite{KricheneCaluyaHalder2020} study noisy second-order dynamics for two-layer networks through a free energy and its stationary distribution, in which the noise plays an essential role.
Wu et~al.\ \cite{WuKungurtsevMondelli2023} study heavy-ball methods for two- and three-layer networks, including mean-field approximation, dropout stability, and connectivity, and their global-convergence result for three layers assumes convergence of the dynamics.
Yang and Littwin \cite{YangLittwin2023} derive infinite-width limits of adaptive optimizers, including Adam, in both kernel and feature-learning regimes, describing the training dynamics as the width grows.
We instead give conditions under which the noiseless mean-field dynamics converge as time grows, while retaining Adam's adaptive state, which changes both the energy dissipation and the compactness estimates.

\textit{Invariance principles and \L{}ojasiewicz inequalities.}
Carrillo et~al.\ \cite{CarrilloGvalaniWu2020} develop an invariance principle for gradient flows on probability measures, using compactness and energy dissipation to identify stationary limits.
Our qualitative proof follows this approach, but the momentum and adaptive variables require their own analysis of trajectories with zero dissipation.
For damped second-order systems, \L{}ojasiewicz inequalities give convergence of bounded trajectories \cite{HarauxJendoubi1998} and rates through modified energies \cite{BegoutBolteJendoubi2015}.
Hauer and Maz\'on \cite{HauerMazon2019} relate such inequalities to finite length and rates for gradient flows in metric spaces.
Our quantitative argument adapts these ideas to a force that depends on a probability law and is weighted by the adaptive state.

\subsection{Organization}

\Cref{sec:main} introduces the model and the main results, together with examples.
\Cref{sec:energy} derives the exact energy identity and the corrected dissipation estimate.
\Cref{sec:omega,sec:lojasiewicz-proof} prove the qualitative convergence theorem and the convergence-rate results, respectively.
Auxiliary results and the proofs of well-posedness and sharpness appear in \cref{sec:auxiliary-proofs,sec:well-posedness,sec:sharpness}.
\Cref{sec:noise,sec:numerics} contain the analysis of the noisy system and our numerical experiments.

\section{Model and main results}\label{sec:main}
This section introduces the zero-noise mean-field Adam flow, its solution concepts, and examples, and then states our two main results, on convergence of the joint law to equilibrium and on convergence rates under a weighted mean-field \L{}ojasiewicz inequality.

\subsection{Notation}\label{sec:notation}

Fix $d\in\NN$ and let
\begin{equation*}
	E\defeq\RR^d\times\RR^d\times[0,\infty)^d.
\end{equation*}
The three vector components of $z=(x,m,v)\in E$ represent position, momentum, and the second-moment variable, respectively.
Powers, square roots, and quotients of vectors are always understood coordinatewise.

For a closed subset $S$ of a Euclidean space, $\SCP(S)$ denotes its Borel probability measures, always equipped with weak convergence unless another topology is specified.
We write $\nu_n\Longrightarrow\nu$ for weak convergence, $\delta_x$ for the Dirac mass at $x$, and $T_\#\nu$ for the pushforward of $\nu$ by a Borel map $T$.
For $r\geq1$, $\SCP_r(S)$ is the set of measures with finite $r$th moment, and $\WW_r$ is the $r$-Wasserstein distance.
If $\mu\in\SCP(E)$, then $\mu^x$ denotes its position marginal.

Recall that a family $\mathcal H\subset\SCP(S)$ is tight if, for every $\eta>0$, some compact set $C_\eta\subset S$ satisfies
\begin{equation*}
	\sup_{\nu\in\mathcal H}\nu(S\setminus C_\eta)<\eta.
\end{equation*}
By Prokhorov's theorem, tightness is equivalent to relative compactness for weak convergence \cite[Theorem~5.1.3]{AmbrosioGigliSavare2008}.

It will be convenient for us to use the bounded-Lipschitz metric
\begin{equation*}
	d_{\mathrm{BL}}(\mu,\nu)
	\defeq\sup_{\substack{\lVert\varphi\rVert_\infty\leq1\\ \Lip(\varphi)\leq1}}
	\left\lvert\int_S\varphi\,\dd\mu-\int_S\varphi\,\dd\nu\right\rvert,
\end{equation*}
which metrizes weak convergence \cite[Section~11.3]{Dudley2002}.
For $A\subset\SCP(S)$, we use the convention
\begin{equation*}
	d_{\mathrm{BL}}(\mu,A)\defeq\inf_{\nu\in A}d_{\mathrm{BL}}(\mu,\nu).
\end{equation*}
We write $\lvert\cdot\rvert$ for the Euclidean norm, $\lVert\cdot\rVert$ for the induced operator norm, and $\NN_0=\{0,1,2,\ldots\}$.
For $r\in\RR$, we define $r_+\defeq\max\{r,0\}$.
For a curve $t\mapsto\mu_t$, we will call the set $\{\mu_t:t\geq0\}$ its \textit{orbit} and the set of position marginals $\{\mu_t^x:t\geq0\}$ its \textit{position orbit}.

\subsection{The dynamics and its solutions}

To formulate the dynamics, we begin by specifying the relationship between the objective function and the force, which will later be used to understand the dissipation of the energy functional.

\begin{asmp}\label[asmp]{asmp:objective}
	The functional $Q\from\SCP_1(\RR^d)\to\RR$ has an associated force
	\begin{equation*}
		G\from\SCP(\RR^d)\times\RR^d\to\RR^d
	\end{equation*}
	with the following properties.
	\begin{enumerate}[label=\textup{(\roman*)},ref=\thethm(\roman*)]
		\item\label[asmp]{asmp:force-regularity} \textit{Boundedness and continuity.} The map $G$ is jointly continuous with respect to the weak topology on $\SCP(\RR^d)$ and
		\begin{equation*}
		K\defeq\sup_{\substack{\mu\in\SCP(\RR^d),\,x\in\RR^d\\1\leq i\leq d}}
		\lvert G_i(\mu,x)\rvert<\infty.
		\end{equation*}
		\item\label[asmp]{asmp:chain-rule} \textit{Chain rule.} On any probability space $(\Omega,\mathcal F,\mathbb P)$ and for every $T>0$, suppose that
		$X\in W^{1,1}(0,T;L^1(\Omega;\RR^d))$; explicitly, there is
		$\dot X\in L^1((0,T)\times\Omega;\RR^d)$ such that
		\begin{equation*}
			X_t=X_0+\int_0^t\dot X_s\,\dd s
			\qquad\text{in }L^1(\Omega;\RR^d).
		\end{equation*}
		Then $t\mapsto Q(\Law(X_t))$ is absolutely continuous and
		\begin{equation}\label{eq:chain-rule}
			\frac{\dd}{\dd t}Q(\Law(X_t))
			=\EE\bigl[G(\Law(X_t),X_t)\cdot\dot X_t\bigr]
		\end{equation}
		for almost every $t\in[0,T]$.
	\end{enumerate}
\end{asmp}

The larger domain of $G$ in \Cref{asmp:force-regularity} is intentional since a weak limit of measures in $\SCP_1(\RR^d)$ may not have a finite first moment, but we need the dynamics to still have a well-defined force.

The chain rule in \Cref{asmp:chain-rule} is the only differentiability property of $Q$ we require.
We use this instead of the usual Lions differentiability condition because it permits us to work with position laws that only have a first moment.
To relate it to the usual calculus on probability measures, let
$\widetilde Q(X)=Q(\Law(X))$ be the lift of $Q$ to square-integrable random variables.
If this lift is Fr\'echet differentiable and
\begin{equation*}
	D\widetilde Q(X)=G(\Law(X),X),
\end{equation*}
then the Hilbert-space chain rule gives \Cref{asmp:chain-rule} for absolutely continuous $L^2$ curves.

To define the dynamics, we fix constant averaging rates $\rho_1,\rho_2>0$ and a regularization parameter $\varepsilon>0$.
For an $E$-valued process $Z=(X,M,V)$, set
\(\mu_t=\Law(Z_t)\).
The zero-noise mean-field Adam system is then
\begin{equation}\label{eq:adam-system}
	\begin{aligned}
		\dot X_t
		&=-\frac{M_t}{\sqrt{V_t+\varepsilon}},\\
		\dot M_t
		&=\rho_1\bigl(G(\mu_t^x,X_t)-M_t\bigr),\\
		\dot V_t
		&=\rho_2\bigl(G(\mu_t^x,X_t)^2-V_t\bigr).
	\end{aligned}
\end{equation}
For $z=(x,m,v)\in E$ and $\nu\in\SCP(\RR^d)$, we denote the mean-field Adam drift by
\begin{equation*}
	b(z,\nu)
	\defeq\left(
	-\frac{m}{\sqrt{v+\varepsilon}},
	\rho_1\bigl(G(\nu,x)-m\bigr),
	\rho_2\bigl(G(\nu,x)^2-v\bigr)
	\right).
\end{equation*}

A \textit{global characteristic solution} is a jointly measurable $E$-valued process $Z=(Z_t)_{t\geq0}$ on a probability space whose sample paths are locally absolutely continuous on $[0,\infty)$ and satisfy \eqref{eq:adam-system} for almost every $t\geq0$, almost surely.
We state our main theorem below conditional on the existence of such a solution.
\Cref{prop:well-posedness} gives a convenient sufficient condition for existence that covers the standard initialization $M_0=V_0=0$.
Along a solution, we write
\begin{equation*}
	g_t\defeq G(\mu_t^x,X_t),
	\qquad
	a_{t,i}\defeq\sqrt{V_{t,i}+\varepsilon},
	\qquad 1\leq i\leq d,
\end{equation*}
for the force acting along the characteristic and the denominators in the position equation.

For later compactness arguments, we will also need a weak formulation of the dynamics.
Let $C_c^1(E)$ consist of the restrictions to $E$ of compactly supported $C^1$ functions defined on a neighborhood of $E$.
For a characteristic solution $Z_s$, $\varphi\in C_c^1(E)$, and $0\leq s\leq t$, the chain rule along each characteristic gives
\begin{equation*}
	\varphi(Z_t)-\varphi(Z_s)
	=\int_s^t\nabla\varphi(Z_r)\cdot b(Z_r,\mu_r^x)\,\dd r.
\end{equation*}
The integrand is bounded because $\nabla\varphi$ has compact support and $G$ is bounded.
Taking expectations and applying Fubini's theorem therefore yields
\begin{equation}\label{eq:weak-equation}
	\int_E\varphi\,\dd\mu_t-\int_E\varphi\,\dd\mu_s
	=\int_s^t\int_E
	\nabla\varphi(z)\cdot b(z,\mu_r^x)\,\mu_r(\dd z)\dd r.
\end{equation}

Thus every global characteristic solution induces a weakly continuous forward solution of the corresponding continuity equation on $[0,\infty)$.
For the argument in \Cref{sec:omega}, we will also define \textit{complete weak solutions} to be weakly continuous curves $\lambda\from\RR\to\SCP(E)$ satisfying \eqref{eq:weak-equation}, with $\mu$ replaced by $\lambda$, for every $s\leq t$ in $\RR$.
The complete solutions used in \Cref{sec:omega} arise, under the compactness hypotheses imposed below, as limits of late-time translates of the original forward solution.
Their extension to negative times is essential for ruling out bounded trajectories with nonzero second-moment components that decay exponentially forward in time and therefore grow exponentially backward in time.

\begin{asmp}\label[asmp]{asmp:well-posedness}
	For the force $G$ and an initial law $\mu_0\in\SCP(E)$, the following hold.
	\begin{enumerate}[label=\textup{(\roman*)},ref=\thethm(\roman*)]
		\item\label[asmp]{asmp:force-lipschitz} \textit{Lipschitz continuity.} There is $L>0$ such that
		\begin{equation}\label{eq:force-wasserstein-lipschitz}
			\lvert G(\mu,x)-G(\nu,y)\rvert
			\leq L\bigl(\WW_1(\mu,\nu)+\lvert x-y\rvert\bigr)
		\end{equation}
		for all $\mu,\nu\in\SCP_1(\RR^d)$ and $x,y\in\RR^d$.
		\item\label[asmp]{asmp:bounded-initial-momentum} \textit{Initial data.} The position marginal satisfies $\mu_0^x\in\SCP_1(\RR^d)$, and there is $B_0<\infty$ such that
		\begin{equation*}
			\mu_0\bigl(\{(x,m,v):\lvert m\rvert\leq B_0\}\bigr)=1.
		\end{equation*}
	\end{enumerate}
\end{asmp}

\begin{prop}\label[prop]{prop:well-posedness}
	Suppose that \Cref{asmp:force-regularity,asmp:well-posedness} hold.
	For every $E$-valued random variable $Z_0$ with law $\mu_0$, there is a unique global characteristic solution of \eqref{eq:adam-system} starting from $Z_0$, with law depending only on $\mu_0$.
	Moreover, with $B\defeq\max\{B_0,\sqrt dK\}$,
	\begin{equation}\label{eq:moment-explicit-bounds}
		\lvert M_t\rvert\leq B,
		\qquad
		0\leq V_{t,i}\leq e^{-\rho_2t}V_{0,i}+K^2,
		\qquad 1\leq i\leq d,
	\end{equation}
	almost surely for every $t\geq0$.
\end{prop}

The proof is a standard fixed-point argument for a law-dependent ordinary differential equation and is thus deferred to \Cref{sec:well-posedness}.
For the proof we use boundedness and continuity of the force, together with \Cref{asmp:well-posedness}, but do not require the chain rule or any parameter restrictions.
The initial law of $V_0$ may be any probability measure on $[0,\infty)^d$, with no moment assumption.

\subsection{Examples}

The following examples illustrate the relationship between the objective and the force.
In the first, an external potential $U$ acts on each position, while a symmetric pair interaction $W$ makes the force depend on the position law.
The second example is a confining potential with bounded force and no interaction.
The final example shows how the law dependence can produce several critical laws even when each frozen force has a unique zero.

\begin{ex}\label[ex]{ex:interaction}
	Suppose that $U\in C^1(\RR^d)$ is bounded from below, $\nabla U$ is bounded and globally Lipschitz, and
	\begin{equation*}
		U(x)\geq c\lvert x\rvert-C
	\end{equation*}
	for some $c>0$ and $C\geq0$.
	Let $W\in C^1(\RR^d\times\RR^d)$ be bounded and symmetric, and suppose that $\nabla_1W$ is bounded and globally Lipschitz in both variables.
	Then
	\begin{equation*}
		Q(\nu)
		=\int_{\RR^d}U(x)\,\nu(\dd x)
		+\frac12\int_{\RR^d\times\RR^d}W(x,y)\,\nu(\dd x)\nu(\dd y)
	\end{equation*}
	satisfies \Cref{asmp:objective} and has tight sublevel sets, with
	\begin{equation*}
		G(\nu,x)=\nabla U(x)+\int_{\RR^d}\nabla_1W(x,y)\,\nu(\dd y).
	\end{equation*}
	Here $G$ is the gradient in $x$ of the linear functional derivative $U(x)+\int_{\RR^d}W(x,y)\,\nu(\dd y)$.
	For a curve as in \Cref{asmp:chain-rule}, let $(Y,\dot Y)$ be an independent copy of $(X,\dot X)$.
	The bounded derivatives, the ordinary chain rule,  Fubini's theorem, and the symmetry of $W$ give
	\begin{equation*}
		\frac{\dd}{\dd t}\frac12\EE W(X_t,Y_t)
		=\EE\left[
		\left(\int_{\RR^d}\nabla_1W(X_t,y)\,\Law(X_t)(\dd y)\right)
		\cdot\dot X_t
		\right].
	\end{equation*}
	Together with the corresponding identity for $\EE U(X_t)$, this proves \eqref{eq:chain-rule}.
	The Lipschitz assumptions also give \Cref{asmp:force-lipschitz}.
	Moreover, $Q(\nu)\leq a$ implies
	\begin{equation*}
		c\int_{\RR^d}\lvert x\rvert\,\nu(\dd x)
		\leq a+C+\frac12\lVert W\rVert_\infty,
	\end{equation*}
	so every nonempty sublevel set is tight.
\end{ex}

\begin{ex}\label[ex]{ex:single-well}
	Consider
	\begin{equation*}
		Q(\nu)=\int_{\RR^d}\bigl(\sqrt{1+\lvert x\rvert^2}-1\bigr)\,\nu(\dd x),
		\qquad
		G(\nu,x)=\frac{x}{\sqrt{1+\lvert x\rvert^2}}.
	\end{equation*}
	The hypotheses of \Cref{ex:interaction} hold with $W=0$, and the coordinatewise force bound is $K=1$.
	The force is independent of $\nu$ and vanishes exactly at the origin, while the potential grows linearly as $\lvert x\rvert\to\infty$.
\end{ex}

\begin{ex}
\label[ex]{ex:three-critical-laws}
	Let $d=1$, $c>1$, and $\vartheta(\nu)\defeq\int_{\RR}\tanh y\,\nu(\dd y)$.
	Consider
	\begin{align*}
		Q(\nu)&\defeq\int_{\RR}\log\cosh x\,\nu(\dd x)-\frac c2\vartheta(\nu)^2,
		\qquad \nu\in\SCP_1(\RR), \\
		G(\nu,x)&=\tanh x-c\vartheta(\nu)\operatorname{sech}^2x.
	\end{align*}
	This is an instance of \cref{ex:interaction} with $U(x)=\log\cosh x$ and $W(x,y)=-c\tanh x\tanh y$.
	For every fixed law $\nu$, the force has exactly one zero, since
	\begin{equation*}
		G(\nu,x)=0\quad\Longleftrightarrow\quad\frac12\sinh(2x)=c\vartheta(\nu).
	\end{equation*}
	A critical law, satisfying $G(\nu,x)=0$ for $\nu$-almost every $x$, must therefore be a Dirac mass $\delta_x$.
	The self-consistency equation $G(\delta_x,x)=\tanh x(1-c\operatorname{sech}^2x)=0$ gives exactly three such laws:
	\begin{equation*}
		\delta_{-x_*},\quad\delta_0,\quad\delta_{x_*},
		\qquad x_*\defeq\operatorname{arctanh}\sqrt{1-\frac1c}.
	\end{equation*}

	This illustrates the additional self-consistency constraint in the mean-field system.
	For the finite-dimensional Adam dynamics of Barakat and Bianchi \cite[Theorem~3.2]{BarakatBianchi2021}, the position approaches the zeros of a fixed gradient $\nabla F$ under their hypotheses, and therefore a unique zero determines the limit.
	Here the unique zero of $G(\nu,\cdot)$ moves with $\nu$, and three equilibrium laws exist despite uniqueness of the zero for each frozen law.
\end{ex}

\subsection{First main result: qualitative convergence}
\label{sec:convergence-framework}
For $r>0$ and $\mu\in\SCP(E)$, define the \textit{weighted second moment of the momentum} by
\begin{equation*}
	\mathcal A_r(\mu)
	\defeq\int_E\sum_{i=1}^d\frac{m_i^2}{(v_i+\varepsilon)^r}\,\mu(\dd x,\dd m,\dd v).
\end{equation*}
For the law $\mu_t$ of a characteristic solution, the position equation gives
\begin{equation*}
	\mathcal A_1(\mu_t)
	=\EE\sum_{i=1}^d\frac{M_{t,i}^2}{V_{t,i}+\varepsilon}
	=\EE\lvert\dot X_t\rvert^2.
\end{equation*}
Thus $\mathcal A_1$ can be directly interpreted as the mean-square speed.
In contrast, $\mathcal A_{1/2}$ is the kinetic correction whose derivative cancels the force--momentum term arising from the time derivative of $Q(\mu_t^x)$.

Whenever $\mu^x\in\SCP_1(\RR^d)$ and $\mathcal A_{1/2}(\mu)<\infty$, we define the energy functional
\begin{equation*}
	\FF(\mu)
	\defeq Q(\mu^x)+\frac{1}{2\rho_1}\mathcal A_{1/2}(\mu).
\end{equation*}

For our convergence analysis, we further define the parameters
\begin{equation*}
	\alpha\defeq\frac{\rho_2}{4\rho_1},
	\qquad
	\beta\defeq(\alpha-1)_+,
	\qquad
	p\defeq\min\left\{\frac12,\frac{2\rho_1}{\rho_2}\right\},
	\qquad
	\delta\defeq\varepsilon-\beta K^2.
\end{equation*}
The ratio $\alpha$ compares the averaging rates, and $\beta$ records its excess above one.
The exponent $p$ determines the weighted initial moment propagated by the equations.
Its definition ensures both that $p\leq1/2$, which allows bounds on $\mathcal A_p$ to control $\mathcal A_{1/2}$ and $\mathcal A_1$, and that $p\leq2\rho_1/\rho_2$, so that
\begin{equation*}
	e^{-(2\rho_1-p\rho_2)t}\leq1, \qquad t\geq0.
\end{equation*}
This prevents the contribution of the initial momentum from growing when it is divided by the possibly decaying denominator $(V_{t,i}+\varepsilon)^p$.
Finally, $\delta$ measures the margin in the dissipation condition below. We require the following conditions on the initial law, the averaging parameters, and the long-time behavior of the position laws.

\begin{asmp}\label[asmp]{asmp:convergence}
	For a global characteristic solution $Z = (X,M,V)$ with law $\mu_t$ and initial law $\mu_0\in\SCP(E)$, the following hold.
	\begin{enumerate}[label=\textup{(\roman*)},ref=\thethm(\roman*)]
		\item\label[asmp]{asmp:initial-integrability} \textit{Initial integrability.} The initial law satisfies
		\begin{equation*}
			\mu_0^x\in\SCP_1(\RR^d),
			\qquad
			\mathcal A_p(\mu_0)<\infty.
		\end{equation*}
		\item\label[asmp]{asmp:rate-condition} \textit{Averaging rates.} The parameters satisfy
		\begin{equation*}
			\left(\frac{\rho_2}{4\rho_1}-1\right)_+K^2\leq\varepsilon,
		\end{equation*}
		or, equivalently, $\delta\geq0$.
		\item\label[asmp]{asmp:trajectory-bounds} \textit{Bounds along the trajectory.} The family $(\mu_t^x)_{t\geq0}$ is tight, and
		\begin{equation}\label{eq:orbit-lower-bound}
			\inf_{t\geq0}Q(\mu_t^x)>-\infty.
		\end{equation}
	\end{enumerate}
\end{asmp}

Initial integrability justifies the energy estimates, and the rate condition makes the dissipation nonnegative apart from an integrable error due to the initial second-moment variable.
The trajectory bounds provide compactness of the position laws and a lower bound for the energy.
\Cref{cor:sublevel-convergence} gives a criterion for these bounds in terms of the objective and the initial law.

The initial-data condition in \Cref{asmp:bounded-initial-momentum} implies \Cref{asmp:initial-integrability}, since $\mathcal A_p(\mu_0)\leq\varepsilon^{-p}B_0^2$.
When $K>0$, \Cref{asmp:rate-condition} is equivalent to
\begin{equation}\label{eq:extended-rate-condition}
	\rho_2\leq4\rho_1\left(1+\frac{\varepsilon}{K^2}\right).
\end{equation}
It holds, in particular, whenever $\rho_2\leq4\rho_1$, which is the standard range for convergence results on continuous-time Adam \cite{BelottoGazeau2020,BarakatBianchi2021}.
The regularization permits larger values of $\rho_2$ relative to $\rho_1$, with the admissible excess depending on the force bound $K$.
Apart from its finiteness, the proofs below use $K$ only through the bound $\lvert G_i(\mu_t^x,X_t)\rvert\leq K$ along the solution.
The results that assume \Cref{asmp:rate-condition} therefore remain valid if $K$ is replaced, in that condition and in the definitions of $\beta$ and $\delta$, by any constant $K'$ such that $\lvert G_i(\mu_t^x,X_t)\rvert\leq K'$ for all $t\geq0$ and $1\leq i\leq d$, almost surely.

We next identify the stationary laws that can occur as long-time limits and state our main convergence theorem.

\begin{defn}
	The set of \textit{critical laws} is
	\begin{equation*}
		\mathcal C
		\defeq\bigl\{\nu\in\SCP(\RR^d):G(\nu,x)=0\text{ for }\nu\text{-almost every }x\bigr\}.
	\end{equation*}
	The corresponding set of \textit{equilibrium laws} on $E$ is
	\begin{equation*}
		\mathcal S
		\defeq\bigl\{\nu\otimes\delta_0\otimes\delta_0:\nu\in\mathcal C\bigr\}.
	\end{equation*}
\end{defn}

A critical law is supported on the set of positions where its associated force vanishes.
With this position law and $m=v=0$, every component of the drift is zero, which explains the form of the equilibrium laws in $\mathcal S$.
The continuity assumption on $G$ makes $\mathcal C$ and $\mathcal S$ weakly closed, which is proved in \Cref{lem:closed-critical}.
Recall for use in the following theorem that a set is \textit{totally disconnected} if each of its connected subsets contains at most one point.

\begin{thm}\label{thm:main}
	Suppose that \Cref{asmp:objective} holds, and let $Z$ be a global characteristic solution satisfying \Cref{asmp:convergence}.
	Then:
	\begin{enumerate}[label=\textup{(\roman*)}]
		\item The family $(\mu_t)_{t\geq0}$ is tight in $\SCP(E)$ and every weak limit point of $\mu_t$ as $t\to\infty$ belongs to $\mathcal S$.
		In particular,
		\begin{equation*}
			\lim_{t\to\infty}d_{\mathrm{BL}}(\mu_t,\mathcal S)=0.
		\end{equation*}
		\item The weighted second moment of the momentum and the mean-square force tend to zero:
		\begin{equation}\label{eq:momentum-force-limits}
			\lim_{t\to\infty}\mathcal A_p(\mu_t)=0,
			\qquad
			\lim_{t\to\infty}\int_{\RR^d}\lvert G(\mu_t^x,x)\rvert^2\,\mu_t^x(\dd x)=0.
		\end{equation}
		There is $q_\infty\in\RR$ such that
		\begin{equation}\label{eq:objective-energy-limits}
			\lim_{t\to\infty}Q(\mu_t^x)
			=\lim_{t\to\infty}\FF(\mu_t)=q_\infty.
		\end{equation}
		\item If
		\begin{equation}\label{eq:intrinsic-critical-set}
			\mathcal C\cap
			\overline{\{\mu_t^x:t\geq0\}}
		\end{equation}
		where the closure is taken in the weak topology, is totally disconnected, then there is
		\begin{equation*}
			\nu_\infty\in
			\mathcal C\cap\overline{\{\mu_t^x:t\geq0\}}
		\end{equation*}
		such that
		\begin{equation}\label{eq:weak-convergence}
			\mu_t\Longrightarrow\nu_\infty\otimes\delta_0\otimes\delta_0
			\qquad\text{as }t\to\infty.
		\end{equation}
	\end{enumerate}
\end{thm}

The total-disconnectedness hypothesis in \textup{(iii)} holds, for instance, when the set \eqref{eq:intrinsic-critical-set} is finite.
We note that this concerns critical laws rather than zeros of the force, as even if $x\mapsto G(\nu,x)$ has a unique zero for every fixed $\nu$, there may be several critical laws, as in \Cref{ex:three-critical-laws}.
When the force is independent of the law and vanishes at two distinct positions, every mixture of the corresponding Dirac masses is critical, producing a connected family of critical laws.
We also remark that the limit $q_\infty$ need not equal $Q(\nu_\infty)$ since a weak limit may have infinite first moment, and even when $\nu_\infty\in\SCP_1(\RR^d)$, weak convergence alone need not preserve $Q$.
Later, in \Cref{thm:kl-rates}, we show convergence in $\WW_1$, where $Q$ is continuous by \eqref{eq:q-wasserstein-lipschitz}.

The trajectory bounds in \Cref{asmp:trajectory-bounds} can be verified using properties of a sublevel set of $Q$.
For $a\in\RR$, define
\begin{equation*}
	\mathcal L_a\defeq\{\nu\in\SCP_1(\RR^d):Q(\nu)\leq a\},
	\qquad
	\mathcal C_a\defeq\mathcal C\cap\overline{\mathcal L_a}.
\end{equation*}
Thus $\mathcal L_a$ is the sublevel set of $Q$ at level $a$.

The criterion below uses $\mathcal L_a$ at a particular finite level $a_0$ which is determined by the initial law and the model parameters.
Specifically, $a_0$ equals the initial energy $\FF(\mu_0)$ plus a nonnegative correction accounting for any initial excess of the second-moment variable above $K^2$.
The estimate in \Cref{prop:corrected-energy} gives $Q(\mu_t^x)\leq a_0$ for every $t\geq0$, so the entire position orbit lies in $\mathcal L_{a_0}$.
We give the explicit formula in \eqref{eq:a0}, where the correction is derived.
Note that for the standard initialization $M_0=V_0=0$, the correction and initial momentum term vanish, giving $a_0=Q(\mu_0^x)$.

\begin{cor}\label[cor]{cor:sublevel-convergence}
	Suppose that \Cref{asmp:objective} holds and that $Z$ is a global characteristic solution satisfying \Cref{asmp:initial-integrability,asmp:rate-condition}.
	If the sublevel set $\mathcal L_{a_0}$ is tight, then \Cref{asmp:trajectory-bounds} holds, and all conclusions of \Cref{thm:main} follow.
	If $\mathcal C_{a_0}$ is totally disconnected, then  $\mu_t\Longrightarrow\nu_\infty\otimes\delta_0\otimes\delta_0$ for some $\nu_\infty\in\mathcal C_{a_0}$.
\end{cor}

\subsection{Sharpness of the averaging-rate condition}

The next result shows that \cref{asmp:rate-condition} cannot be replaced by a weaker condition on $\rho_1$, $\rho_2$, $\varepsilon$, and $K$.
When it fails, the energy $\FF$ can increase even when every second-moment variable is at most $K^2$, and the next result shows that these increases can sustain a periodic solution.

\begin{prop}\label[prop]{prop:sharpness}
	Let $d\geq3$, and let $\rho_1,\rho_2,\varepsilon,K>0$ be such that \eqref{eq:extended-rate-condition} does not hold.
	Then there exists $U\in C^\infty(\RR^d)$ with bounded Hessian and $C\geq0$ such that
	\begin{equation*}
		\sup_{\substack{y\in\RR^d\\1\leq i\leq d}}\lvert\partial_iU(y)\rvert=K,
		\qquad
		U(y)\geq K\lvert y\rvert-C
		\quad\text{for }y\in\RR^d,
	\end{equation*}
	and such that, for $Q(\nu)=\int_{\RR^d}U\,\dd\nu$, the system \eqref{eq:adam-system} has a deterministic, periodic global characteristic solution $Z_t=(x(t),m(t),v(t))$.
	This objective and solution satisfy all hypotheses of \cref{thm:main,cor:sublevel-convergence} except \cref{asmp:rate-condition}, but the law $\mu_t=\delta_{Z_t}$ satisfies
	\begin{equation*}
		\inf_{t\geq0}\mathcal A_p(\mu_t)>0,
		\qquad
		\inf_{t\geq0}d_{\mathrm{BL}}(\mu_t,\mathcal S)>0.
	\end{equation*}
\end{prop}

In the construction, given in \cref{sec:sharpness}, we use a rapidly oscillating force to raise one second-moment variable while producing little momentum, and a weak force then produces momentum while this variable is still large, so that the energy increases.
Two further coordinates rotate with positive dissipation, which balances this increase over each period.

Since the solution is a single particle and the force $\nabla U$ does not depend on the law, the example also shows that \eqref{eq:extended-rate-condition} is sharp for the finite-dimensional version of \eqref{eq:adam-system}.
The same behavior can be realized in dimensions one and two, using three particles in disjoint regions and an objective that depends on three localized first moments of their position law.
It can also be realized from zero initial memory, by adding a slowly decaying coordinate that gradually switches on a translated and time-reparametrized copy of the construction, so that the solution approaches the periodic orbit.
We omit the details of these two modifications.

\subsection{A regularized two-layer network}\label{sec:neural-convergence}

The qualitative theorem applies to the teacher--student model analyzed numerically in \cref{sec:numerics-teacher}.
Let $\xi\sim N(0,I_8)$, let $y(\xi)$ be a bounded target, and let $\gamma>0$.
For neuron parameters $x=(u,w)\in\RR\times\RR^8$, set
\begin{equation*}
	\phi(x;\xi)\defeq2\tanh(u)\operatorname{erf}(w\cdot\xi),
\end{equation*}
and define
\begin{equation}\label{eq:neural-objective}
	\begin{aligned}
		f_\nu(\xi)&\defeq\int_{\RR^9}\phi(x;\xi)\,\nu(\dd x),\\
		Q(\nu)&\defeq\frac12\EE\lvert f_\nu(\xi)-y(\xi)\rvert^2
		+\gamma\int_{\RR^9}\bigl(\sqrt{1+\lvert x\rvert^2}-1\bigr)\,\nu(\dd x).
	\end{aligned}
\end{equation}
The force is
\begin{equation*}
	G(\nu,x)=\EE\bigl[(f_\nu(\xi)-y(\xi))\nabla_x\phi(x;\xi)\bigr]
	+\gamma\frac{x}{\sqrt{1+\lvert x\rvert^2}}.
\end{equation*}
All expectations in these formulas are over the input $\xi$.

\begin{cor}\label[cor]{cor:neural-convergence}
	Let $Q$ be given by \eqref{eq:neural-objective}, let $K>0$ bound each coordinate of its force, and suppose that \cref{asmp:rate-condition} holds.
	For every initial position law $\nu_0\in\SCP_1(\RR^9)$ and the standard initialization $M_0=V_0=0$, there is a unique global characteristic solution of \eqref{eq:adam-system}.
	Its joint law satisfies
	\begin{equation*}
		\lim_{t\to\infty}d_{\mathrm{BL}}(\mu_t,\mathcal S)=0,
		\qquad
		\lim_{t\to\infty}\int_{\RR^9}\lvert G(\mu_t^x,x)\rvert^2\,\mu_t^x(\dd x)=0,
	\end{equation*}
	and $Q(\mu_t^x)$ converges to a finite value.
\end{cor}

\begin{proof}
	Expanding the square in \eqref{eq:neural-objective} gives the interaction form of \cref{ex:interaction}, up to the constant $\frac12\EE\lvert y(\xi)\rvert^2$, with
	\begin{equation*}
		\begin{aligned}
			U(x)&=-\EE\bigl[y(\xi)\phi(x;\xi)\bigr]+\gamma\bigl(\sqrt{1+\lvert x\rvert^2}-1\bigr),\\
			W(x,z)&=\EE\bigl[\phi(x;\xi)\phi(z;\xi)\bigr].
		\end{aligned}
	\end{equation*}
	The function $\phi$ is bounded, and its first and second derivatives in $x$ are bounded by constants times $1+\lvert\xi\rvert$ and $1+\lvert\xi\rvert^2$, respectively.
	The bounded target and Gaussian moments therefore give bounded, globally Lipschitz gradients for $U$ and $W$, and $U$ is confining.
	Thus \cref{asmp:objective,asmp:force-lipschitz} hold, and \cref{prop:well-posedness} gives a unique global solution.
	Moreover, the nonnegative square loss gives
	\begin{equation*}
		\int_{\RR^9}\lvert x\rvert\,\nu(\dd x)\leq1+\frac{Q(\nu)}{\gamma},
	\end{equation*}
	so the objective sublevel sets are tight.
	For zero initial memory, $a_0=Q(\nu_0)$, and \cref{cor:sublevel-convergence} gives the conclusions.
\end{proof}

This applies both to empirical initial laws with any finite number of neurons and to general initial laws in $\SCP_1(\RR^9)$.
Note, however, that the limiting critical laws do not necessarily minimize the objective, and the solution may not approach a single law.
This latter conclusion requires further assumptions, such as those of the next theorem.

\subsection{Second main result: convergence rates under a mean-field
\texorpdfstring{\L{}ojasiewicz}{Lojasiewicz} inequality}

A \L{}ojasiewicz inequality relates the objective gap to the size of the force \cite{Kurdyka1998,HauerMazon2019}.
Together with the dissipation and regularity conditions below, it gives convergence to a single critical law and a rate of approach, even when the critical set is connected.
To allow unbounded second-moment variables, we use the weighted force
\begin{equation*}
	\mathcal B(\mu)
	\defeq\int_E\sum_{i=1}^d
	\frac{G_i(\mu^x,x)^2}{\sqrt{v_i+\varepsilon}}\,\mu(\dd x,\dd m,\dd v),
	\qquad \mu\in\SCP(E).
\end{equation*}
This is finite by the force bound in \Cref{asmp:force-regularity}.

\begin{asmp}\label[asmp]{asmp:mean-field-kl}
	Let $q\in\RR$ and $\mathcal U\subset\{\mu\in\SCP(E):\mu^x\in\SCP_1(\RR^d)\}$.
	There are $r>0$ and a continuous, concave function $\varphi\from[0,r)\to[0,\infty)$, with $\varphi(0)=0$, such that $\varphi$ is continuously differentiable on $(0,r)$, $\varphi'>0$, and
	\begin{equation}\label{eq:mean-field-kl}
		\varphi'\bigl(Q(\mu^x)-q\bigr)\mathcal B(\mu)^{1/2}\geq1
	\end{equation}
	for every $\mu\in\mathcal U$ with $0<Q(\mu^x)-q<r$.
\end{asmp}

We call \eqref{eq:mean-field-kl} a \textit{weighted mean-field \L{}ojasiewicz inequality} and $\varphi$ its \textit{desingularizing function}.
The family $\mathcal U$ need only contain the eventual orbit, and $q$ does not have to be a minimum of $Q$.
The inequality is required only above the threshold $q$ as the proof controls times with $Q(\mu_t^x)\leq q$ through the momentum terms.
The power choice
\begin{equation}\label{eq:kl-power}
	\varphi(s)=\frac{C_{\mathrm L}}{1-\theta}s^{1-\theta},
	\qquad C_{\mathrm L}>0,\quad 0\leq\theta<1,
\end{equation}
gives the more familiar inequality
$\bigl(Q(\mu^x)-q\bigr)^\theta\leq C_{\mathrm L}\mathcal B(\mu)^{1/2}$ for positive gaps smaller than $r$.
When $v_i\leq\bar v$ almost surely under every law in $\mathcal U$, the weighted inequality \eqref{eq:mean-field-kl} is equivalent, up to rescaling $\varphi$, to the usual mean-field inequality
\begin{equation}\label{eq:unweighted-mean-field-kl}
	\varphi'\bigl(Q(\nu)-q\bigr)
	\left(\int_{\RR^d}\lvert G(\nu,x)\rvert^2\,\nu(\dd x)\right)^{1/2}\geq1
\end{equation}
on the corresponding position laws with $0<Q(\nu)-q<r$.

To connect this force estimate with the energy dissipation, we impose two conditions along the given solution.
Here $\mathcal D$ is the exact dissipation defined in \eqref{eq:dissipation}, with $\frac{\dd}{\dd t}\FF(\mu_t)=-\mathcal D(\mu_t)$.

\begin{asmp}\label[asmp]{asmp:kl-dynamics}
	There are $T_0\geq0$, $\kappa>0$, and $L_g<\infty$ such that:
	\begin{enumerate}[label=\textup{(\roman*)}]
		\item \textit{Coercive dissipation.} For almost every $t\geq T_0$,
		\begin{equation}\label{eq:kl-coercivity}
			\mathcal D(\mu_t)\geq\kappa\mathcal A_{1/2}(\mu_t).
		\end{equation}
		\item \textit{Variation of the force.} The curve $g$ is locally absolutely continuous from $[T_0,\infty)$ into $L^2(\Omega;\RR^d)$, and
		\begin{equation}\label{eq:kl-force-variation}
			\lVert\dot g_t\rVert_{L^2(\Omega)}
			\leq L_g\mathcal A_{1/2}(\mu_t)^{1/2}
		\end{equation}
		for almost every $t\geq T_0$.
	\end{enumerate}
\end{asmp}

The theorem below shows that these conditions follow from the Lipschitz condition on the force in \Cref{asmp:force-lipschitz} under restrictions on the averaging parameters.
\begin{thm}\label{thm:kl-rates}
	Suppose that the objective and force satisfy \Cref{asmp:objective}, and let $Z$ be a global characteristic solution satisfying \Cref{asmp:convergence}.
	Suppose also that \Cref{asmp:kl-dynamics} holds, or that the force satisfies the Lipschitz condition in \Cref{asmp:force-lipschitz} and one of the following holds:
	\begin{enumerate}[label=\textup{(\alph*)}]
		\item $\rho_2<4\rho_1$;
		\item $\delta>0$, and $V_{0,i}\leq\bar v_0$ almost surely for every $i$, for some $\bar v_0<\infty$.
	\end{enumerate}
	Let $q_\infty$ be the common limit of the objective and energy in \eqref{eq:objective-energy-limits}.
	Suppose that the mean-field \L{}ojasiewicz inequality in \Cref{asmp:mean-field-kl} holds with $q=q_\infty$ and that $\mu_t\in\mathcal U$ for all sufficiently large $t$.
	Then:
	\begin{enumerate}[label=\textup{(\roman*)}]
		\item The position has finite length in the sense that
		\begin{equation}\label{eq:kl-finite-length}
			\int_0^\infty\lVert\dot X_t\rVert_{L^2(\Omega)}\,\dd t<\infty.
		\end{equation}
		There is $X_\infty\in L^1(\Omega;\RR^d)$ such that
		\begin{equation*}
			(X_t,M_t,V_t)\longrightarrow(X_\infty,0,0)
			\quad\text{almost surely},
		\end{equation*}
		with $\nu_\infty\defeq\Law(X_\infty)\in\mathcal C$ and $Q(\nu_\infty)=q_\infty$.
		Moreover, $\mu_t^x\to\nu_\infty$ in $\WW_1$, and also in $\WW_2$ if $\mu_0^x\in\SCP_2(\RR^d)$.
		\item If the desingularizing function has the power form \eqref{eq:kl-power}, there are $T\geq0$ and $c,C>0$ such that, for $t\geq T$,
		\begin{align}
			\lvert Q(\mu_t^x)-q_\infty\rvert
			+\mathcal A_{1/2}(\mu_t)+\mathcal B(\mu_t)
			&\leq C E_\theta(t-T),\label{eq:kl-objective-rate}\\
			\lVert X_t-X_\infty\rVert_{L^2(\Omega)}
			\leq\int_t^\infty\lVert\dot X_s\rVert_{L^2(\Omega)}\,\dd s
			&\leq C R_\theta(t-T),\label{eq:kl-position-rate}
		\end{align}
		where
		\begin{equation}\label{eq:kl-power-rates}
			\begin{array}{c|c|c}
			& E_\theta(\tau)&R_\theta(\tau)\\ \hline
			0\leq\theta\leq\tfrac12
			&e^{-c\tau}&e^{-c\tau/2}\\[2pt]
			\tfrac12<\theta<1
			&(1+\tau)^{-1/(2\theta-1)}
			&(1+\tau)^{-(1-\theta)/(2\theta-1)}.
			\end{array}
		\end{equation}
		\item In case \textup{(b)}, if the desingularizing function has the power form \eqref{eq:kl-power}, then there are $T\geq0$ and $c,C>0$ such that, for all $t\geq T$,
		\begin{equation}\label{eq:kl-state-moment-rates}
			\EE\lvert M_t\rvert^2+\EE\lvert g_t\rvert^2
			+\EE\lvert V_t\rvert+\EE\lvert V_t\rvert^2
			\leq C E_\theta(t-T).
		\end{equation}
		Writing $\mu_\infty=\nu_\infty\otimes\delta_0\otimes\delta_0$, we also have
		\begin{equation}\label{eq:kl-full-law-rate}
			\WW_1(\mu_t,\mu_\infty)\leq C R_\theta(t-T).
		\end{equation}
		If $\mu_0^x\in\SCP_2(\RR^d)$, the same estimate holds for $\WW_2$.
	\end{enumerate}
\end{thm}

The displacement estimate \eqref{eq:kl-position-rate} requires only $X_0\in L^1$, since $X_t-X_\infty$ is square integrable even when $X_t$ is not.

Case \textup{(b)} uses the strict margin $\delta>0$ to verify coercivity.
At $\delta=0$, the convergence-rate result in \Cref{thm:kl-rates} still applies if the dissipation and force-variation bounds in \Cref{asmp:kl-dynamics} hold along the trajectory.
Case \textup{(a)} requires no additional moment of $V_0$, but a rate for the full state law needs further tail control, as even with $G=M_0=0$, one has $V_t=e^{-\rho_2t}V_0$, whose convergence in $d_{\mathrm{BL}}$ need not be exponential.
We provide proofs in \Cref{sec:lojasiewicz-proof}.

\subsection{Dependence on the number of particles}\label{sec:width-dependence}

The rate theorem uses an inequality on probability laws.
An inequality at each finite particle count need not have a constant uniform in the number of particles, even when every finite system has a positive local exponential rate.
The interacting model of \cref{sec:numerics-meanfield} gives an explicit example.
Its confinement is $a\log\cosh x$, and its repulsion approaches the kernel $-\lvert x-y\rvert$ as a smoothing parameter $\varsigma$ tends to zero.
For the limiting interaction, the objective on $N$ particles is
\begin{equation*}
	Q_N(x_1,\ldots,x_N)
	=\frac aN\sum_{j=1}^N\log\cosh x_j
	-\frac1{2N^2}\sum_{j,k=1}^N\lvert x_j-x_k\rvert,
	\qquad a\geq1.
\end{equation*}
The ordered critical configuration is
\begin{equation*}
	x_j^*=\operatorname{arctanh}\frac{2j-1-N}{aN},\qquad 1\leq j\leq N.
\end{equation*}
In a neighborhood where the particles remain distinct and ordered, the force Jacobian $\mathsf H=N\nabla^2Q_N$ is diagonal.
At this critical configuration, its smallest eigenvalue is
\begin{equation*}
	\lambda_N=a\left(1-\left(\frac{N-1}{aN}\right)^2\right).
\end{equation*}
The slowest linearized position mode has exponential envelope rate
\begin{equation*}
	r_N=\frac12\left(\rho_1-\operatorname{Re}\sqrt{\rho_1^2-4\rho_1\lambda_N/\sqrt\varepsilon}\right),
\end{equation*}
as derived in \cref{sec:numerics-meanfield}.
For $a>1$, $\lambda_N\to a-a^{-1}>0$ and these local rates remain bounded away from zero.
For $a=1$, however,
\begin{equation*}
	\lambda_N\sim\frac2N,
	\qquad
	r_N\sim\frac{2}{N\sqrt\varepsilon}.
\end{equation*}
In the local quadratic approximation, the constant in the \L{}ojasiewicz inequality with exponent $1/2$ grows like $\lambda_N^{-1/2}$.
Thus exponential convergence at every fixed particle count does not supply an exponential rate uniform in $N$.

\begin{figure}[t]
	\centering
	\includegraphics[width=\linewidth]{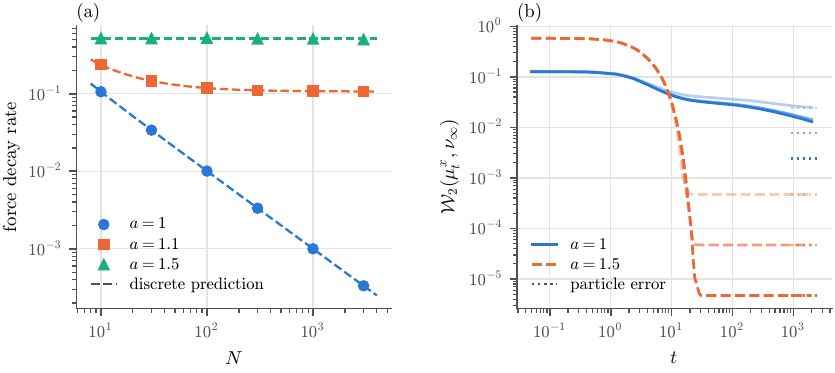}
	\caption{Dependence on the number of particles in the interacting model.
	(a) Measured force decay rates, with predictions from the discrete linearization (dashed).
	For $a=1$ the rate decreases like $1/N$, while for $a>1$ it stays bounded away from zero over the tested range.
	(b) Distance to the continuous law \eqref{eq:repulsive-critical-law}, for $N=10^3,10^4,10^5$ from light to dark.
	The short dotted segments show $\WW_2(\nu_{\infty,N},\nu_\infty)$, the finite-particle approximation errors.}
	\label{fig:numerics-meanfield}
\end{figure}

\Cref{fig:numerics-meanfield} compares these local predictions with computations using positive smoothing.
For this smoothed model, the force Jacobian differs from that of the limiting interaction by terms bounded by a constant times $\varsigma^{-1}e^{-\Delta/\varsigma}$, where $\Delta$ is the smallest particle spacing.
The calculation therefore describes the computations when $\Delta\gg\varsigma$ and its limit as $N\to\infty$ concerns the limiting interaction, not fixed positive smoothing.

The quartic example of \cref{sec:numerics-rates} gives a complementary case.
For $Q(\nu)=\int_\RR((1+x^4)^{1/4}-1)\,\nu(\dd x)$ and laws supported in $[-1,1]$, a pointwise estimate and Jensen's inequality give a mean-field \L{}ojasiewicz inequality with exponent $3/4$ and a constant independent of the number of particles.
With zero initial memory, the adaptive state is bounded, so the weighted inequality holds as well.

\section{Energy estimates}\label{sec:energy}
The position evolves according to the averaged force $M_t$, which is not necessarily aligned with the current force $G(\mu_t^x,X_t)$, so the objective alone may not decrease.
Thus $Q(\mu_t^x)$ is not, by itself, a Lyapunov functional for the Adam flow.
In this section, we derive an exact energy identity, in which the derivative of the weighted momentum term in $\FF$ cancels the force--momentum product, and we show that a corrected dissipation has finite time integral and that the energy converges.
This allows for our $\omega$-limit analysis in \Cref{sec:omega}.

\subsection{The exact identity}

After cancellation, momentum relaxation contributes to energy decay, while the changing weight $(V_{t,i}+\varepsilon)^{-1/2}$ contributes a term whose sign can vary.
When $\alpha\leq1$, the relaxation controls this additional term for every value of $V$.
For larger $\alpha$, the sign depends on the second-moment variable and the force.
These contributions are collected in the following functional.

For $\mu\in\SCP(E)$ with $\mathcal A_{1/2}(\mu)<\infty$, define
\begin{equation}\label{eq:dissipation}
	\mathcal D(\mu)
	\defeq\int_E\sum_{i=1}^d
	\frac{m_i^2\bigl(\varepsilon+(1-\alpha)v_i+\alpha G_i(\mu^x,x)^2\bigr)}
	{(v_i+\varepsilon)^{3/2}}\,\mu(\dd x,\dd m,\dd v).
\end{equation}
The force bound makes this integral absolutely convergent: for a constant $C$ depending only on $\alpha$, $K$, and $\varepsilon$,
\begin{equation*}
	\sum_{i=1}^d
	\frac{m_i^2\left\lvert\varepsilon+(1-\alpha)v_i+\alpha G_i(\mu^x,x)^2\right\rvert}
	{(v_i+\varepsilon)^{3/2}}
	\leq C\sum_{i=1}^d\frac{m_i^2}{\sqrt{v_i+\varepsilon}}.
\end{equation*}
The following proposition identifies $\mathcal D$ as the rate of energy dissipation.
Its first conclusion propagates the weighted momentum moment, which bounds the momentum term in the energy and the mean-square position velocity, providing the integrability needed for the subsequent differentiation.

\begin{prop}\label[prop]{prop:dissipation}
	Suppose that \Cref{asmp:objective} holds and that $Z$ is a global characteristic solution satisfying \Cref{asmp:initial-integrability}.
	Then:
	\begin{enumerate}[label=\textup{(\roman*)}]
		\item The weighted momentum moment satisfies
		\begin{equation}\label{eq:weighted-moment-propagation}
			\sup_{t\geq0}\mathcal A_p(\mu_t)
			\leq2\mathcal A_p(\mu_0)+2dK^2\varepsilon^{-p}.
		\end{equation}
		\item The energy $t\mapsto\FF(\mu_t)$ is locally absolutely continuous, and
		\begin{equation*}
			\FF(\mu_t)+\int_s^t\mathcal D(\mu_r)\,\dd r=\FF(\mu_s),
			\qquad 0\leq s\leq t.
		\end{equation*}
		\item If $\rho_2\leq4\rho_1$, then $\mathcal D\geq0$ and, for every $\mu$ in its domain,
		\begin{equation*}
			\mathcal D(\mu)=0
			\quad\Longleftrightarrow\quad
			\mu\bigl(\{m=0\}\bigr)=1.
		\end{equation*}
	\end{enumerate}
\end{prop}

\begin{proof}
	Summing the pathwise bound \eqref{eq:a-priori-weighted} over $i$ and taking expectations proves (i).

	We next establish the integrability needed to differentiate the two terms of $\FF$.
	For $r\geq p$, the inequality
	$(v_i+\varepsilon)^{-r}\leq\varepsilon^{p-r}(v_i+\varepsilon)^{-p}$ implies
	\begin{equation}\label{eq:weighted-moment-comparison}
		\mathcal A_r(\mu_t)
		\leq\varepsilon^{p-r}\mathcal A_p(\mu_t).
	\end{equation}
	Since $p\leq1/2$, part (i) therefore bounds both $\mathcal A_{1/2}(\mu_t)$ and $\mathcal A_1(\mu_t)$ uniformly in time.
	By the identity $\EE\lvert\dot X_t\rvert^2=\mathcal A_1(\mu_t)$ from \Cref{sec:convergence-framework}, $\dot X$ belongs to $L^2([0,T]\times\Omega)$ and hence to $L^1([0,T]\times\Omega)$ for every $T<\infty$.
	Together with $X_0\in L^1$, this shows that
	\begin{equation*}
		X_t=X_0+\int_0^t\dot X_s\,\dd s
		\qquad\text{in }L^1(\Omega;\RR^d).
	\end{equation*}
	Thus $X$ satisfies \Cref{asmp:chain-rule}.
	Substituting $\dot X_{t,i}=-M_{t,i}/a_{t,i}$ into the chain rule yields
	\begin{equation*}
		\frac{\dd}{\dd t}Q(\mu_t^x)
		=-\EE\sum_{i=1}^d\frac{g_{t,i}M_{t,i}}{a_{t,i}}
	\end{equation*}
	for almost every $t\geq0$.

	For the weighted momentum term, we differentiate along almost every sample path.
	Using $\dot M_{t,i}=\rho_1(g_{t,i}-M_{t,i})$ and
	$\dot a_{t,i}=\rho_2(g_{t,i}^2-V_{t,i})/(2a_{t,i})$, the product rule gives
	\begin{equation*}
		\frac{\dd}{\dd t}\left(\frac{M_{t,i}^2}{2\rho_1a_{t,i}}\right)
		=\frac{M_{t,i}(g_{t,i}-M_{t,i})}{a_{t,i}}
		-\alpha\frac{M_{t,i}^2(g_{t,i}^2-V_{t,i})}{a_{t,i}^3}.
	\end{equation*}
	To take expectations in the integrated form of this equality, we bound its right-hand side.
	For $\lvert g\rvert\leq K$ and $v\geq0$, Young's inequality and $v/(v+\varepsilon)\leq1$ give
	\begin{equation*}
		\left\lvert
		\frac{m(g-m)}{\sqrt{v+\varepsilon}}
		-\alpha\frac{m^2(g^2-v)}{(v+\varepsilon)^{3/2}}
		\right\rvert
		\leq C\left(1+\frac{m^2}{\sqrt{v+\varepsilon}}\right),
	\end{equation*}
	where $C$ depends only on $K,\alpha$, and $\varepsilon$.
	The uniform bound on $\mathcal A_{1/2}(\mu_t)$ makes this bound integrable over $[0,T]\times\Omega$ after summing over $i$.
	Fubini's theorem therefore shows that $t\mapsto\mathcal A_{1/2}(\mu_t)/(2\rho_1)$ is locally absolutely continuous, with derivative equal to the expectation of the sum of the pathwise derivatives above.

	Adding this derivative to that of $Q(\mu_t^x)$ cancels the terms containing $g_{t,i}M_{t,i}/a_{t,i}$.
	The remaining terms give, for almost every $t\geq0$,
	\begin{align*}
		\frac{\dd}{\dd t}\FF(\mu_t)
		&=-\EE\sum_{i=1}^d\left(
		\frac{M_{t,i}^2}{a_{t,i}}
		+\alpha\frac{M_{t,i}^2(g_{t,i}^2-V_{t,i})}{a_{t,i}^3}\right)\\
		&=-\EE\sum_{i=1}^d
		\frac{M_{t,i}^2\bigl(\varepsilon+(1-\alpha)V_{t,i}+\alpha g_{t,i}^2\bigr)}
		{(V_{t,i}+\varepsilon)^{3/2}}
		=-\mathcal D(\mu_t).
	\end{align*}
	Here we used $a_{t,i}^2=V_{t,i}+\varepsilon$ to combine the fractions.
	Integrating from $s$ to $t$ proves (ii).

	Finally, $\rho_2\leq4\rho_1$ is equivalent to $\alpha\leq1$.
	In this case the factor $\varepsilon+(1-\alpha)v_i+\alpha G_i(\mu^x,x)^2$ is at least $\varepsilon>0$.
	Each summand defining $\mathcal D$ is therefore nonnegative and vanishes exactly when $m_i=0$.
	Its integral is zero precisely when every coordinate of $m$ vanishes $\mu$-almost surely, proving (iii).
\end{proof}

\subsection{The corrected dissipation}

When $\alpha>1$, the coefficient of $v_i$ in \eqref{eq:dissipation} is negative, but \Cref{asmp:rate-condition} ensures that each summand is nonnegative whenever $v_i\leq K^2$.
The input $G_i^2$ to the second-moment equation is bounded by $K^2$, so the dynamics cannot generate an excess above this level.
Any initial excess is bounded by $(V_{0,i}-K^2)_+e^{-\rho_2t}$, though it need not disappear in finite time.
We use this bound to separate a nonnegative dissipation from an error whose time integral is finite.

Under \Cref{asmp:rate-condition}, we have $\delta\geq0$ and define the nonnegative, possibly infinite functional
\begin{equation}\label{eq:corrected-dissipation}
	\mathcal D_*(\mu)
	\defeq\int_E\sum_{i=1}^d
	\frac{m_i^2\bigl(\delta+\beta(K^2-v_i)_+
	+\alpha G_i(\mu^x,x)^2\bigr)}
	{(v_i+\varepsilon)^{3/2}}\,\mu(\dd x,\dd m,\dd v),
	\qquad \mu\in\SCP(E).
\end{equation}
For a global characteristic solution $Z$, we bound the contribution of the initial excess by
\begin{equation*}
	\mathcal R(t)
	\defeq\EE\sum_{i=1}^d
	\frac{\beta M_{t,i}^2e^{-\rho_2t}(V_{0,i}-K^2)_+}
	{(V_{t,i}+\varepsilon)^{3/2}}.
\end{equation*}
Exponential decay of the initial excess alone does not establish integrability of this error, since $V_0$ may not have a finite first moment.
Therefore, we must retain the dependence on $V_0$ in the denominator when estimating its time integral.
Separating the decaying initial momentum from the bounded force contribution leads to the two integrals, for $r\geq0$,
\begin{align*}
	\mathsf I(r)
	&\defeq\int_0^\infty
	\frac{e^{-(2\rho_1+\rho_2)t}(r-K^2)_+}
	{(e^{-\rho_2t}r+\varepsilon)^{3/2}}\,\dd t,\\
	\mathsf J(r)
	&\defeq\int_0^\infty
	\frac{e^{-\rho_2t}(r-K^2)_+}
	{(e^{-\rho_2t}r+\varepsilon)^{3/2}}\,\dd t.
\end{align*}
When $\beta>0$, we will bound $\mathsf I(r)$ by a multiple of $(r+\varepsilon)^{-p}$ and $\mathsf J(r)$ by a constant independent of $r$.
The first bound matches the assumed weighted momentum moment.
Under \Cref{asmp:initial-integrability}, set
\begin{align}
	\mathcal R_0
	&\defeq 2\int_E\sum_{i=1}^d
	\beta\bigl(m_i^2\mathsf I(v_i)+K^2\mathsf J(v_i)\bigr)
	\,\mu_0(\dd x,\dd m,\dd v),
	\label{eq:R0}\\
	a_0
	&\defeq\FF(\mu_0)+\mathcal R_0.
	\label{eq:a0}
\end{align}
The next proposition shows that $\mathcal R_0$ bounds the accumulated error, so $a_0$ bounds both the energy and the objective along the trajectory.
If the objective is also bounded from below, adding the remaining error $\int_t^\infty\mathcal R(s)\,\dd s$ to the energy gives a nonincreasing quantity bounded from below.
This yields convergence of the energy and finite total corrected dissipation, which will allow us to identify the limiting trajectories in \Cref{sec:omega}.

\begin{prop}\label[prop]{prop:corrected-energy}
	Suppose that \Cref{asmp:objective} holds and that $Z$ is a global characteristic solution satisfying \Cref{asmp:initial-integrability,asmp:rate-condition}.
	Then:
	\begin{enumerate}[label=\textup{(\roman*)}]
		\item For almost every $t\geq0$,
		\begin{equation*}
			\frac{\dd}{\dd t}\FF(\mu_t)
			\leq-\mathcal D_*(\mu_t)+\mathcal R(t).
		\end{equation*}
		\item The total error satisfies
		\begin{equation*}
			\int_0^\infty\mathcal R(t)\,\dd t\leq\mathcal R_0<\infty,
		\end{equation*}
		and
		\begin{equation*}
			Q(\mu_t^x)\leq\FF(\mu_t)\leq a_0
			\qquad\text{for every }t\geq0.
		\end{equation*}
		\item If the objective is bounded from below along the trajectory as in \eqref{eq:orbit-lower-bound}, then
		\begin{equation}\label{eq:finite-corrected-dissipation}
			\int_0^\infty\mathcal D_*(\mu_t)\,\dd t<\infty,
		\end{equation}
		and $\FF(\mu_t)$ has a finite limit as $t\to\infty$.
	\end{enumerate}
\end{prop}

\begin{proof}
	We first compare $\mathcal D(\mu_t)$ with $\mathcal D_*(\mu_t)$.
	Since $1-\alpha\geq-\beta$ and $\delta=\varepsilon-\beta K^2$, the factor multiplying $M_{t,i}^2/(V_{t,i}+\varepsilon)^{3/2}$ in $\mathcal D(\mu_t)$ satisfies, by \eqref{eq:a-priori-second-moment},
	\begin{align*}
		\varepsilon+(1-\alpha)V_{t,i}+\alpha g_{t,i}^2
		&\geq\delta+\beta(K^2-V_{t,i})+\alpha g_{t,i}^2\\
		&\geq\delta+\beta(K^2-V_{t,i})_++\alpha g_{t,i}^2
		-\beta e^{-\rho_2t}(V_{0,i}-K^2)_+.
	\end{align*}
	Multiplying by $M_{t,i}^2/(V_{t,i}+\varepsilon)^{3/2}$, summing over $i$, and taking expectations yields
	\begin{equation*}
		\mathcal D(\mu_t)\geq\mathcal D_*(\mu_t)-\mathcal R(t).
	\end{equation*}
	To use this comparison in the energy identity, we check that both terms on its right-hand side are finite.
	The bounds $(K^2-V_{t,i})_+\leq K^2$, $g_{t,i}^2\leq K^2$, and $e^{-\rho_2t}V_{0,i}\leq V_{t,i}$ give
	\begin{equation*}
		\mathcal D_*(\mu_t)
		\leq\bigl(\delta+(\beta+\alpha)K^2\bigr)\mathcal A_{3/2}(\mu_t),
		\qquad
		\mathcal R(t)\leq\beta\mathcal A_{1/2}(\mu_t),
	\end{equation*}
	and the right-hand sides are finite by \eqref{eq:weighted-moment-propagation} and \eqref{eq:weighted-moment-comparison}.
	Combining this comparison with the identity $\frac{\dd}{\dd t}\FF(\mu_t)=-\mathcal D(\mu_t)$ from \Cref{prop:dissipation} proves (i).

	We next bound the total error.
	Inserting the bounds $M_{t,i}^2\leq2e^{-2\rho_1t}M_{0,i}^2+2K^2$ and $V_{t,i}\geq e^{-\rho_2t}V_{0,i}$ from \Cref{lem:a-priori} into the definition of $\mathcal R(t)$ gives
	\begin{equation*}
		\mathcal R(t)\leq2\beta\EE\sum_{i=1}^d
		\frac{\bigl(e^{-2\rho_1t}M_{0,i}^2+K^2\bigr)e^{-\rho_2t}(V_{0,i}-K^2)_+}
		{(e^{-\rho_2t}V_{0,i}+\varepsilon)^{3/2}}.
	\end{equation*}
	Integrating the resulting two nonnegative terms in time and applying Tonelli's theorem gives
	\begin{equation*}
		\int_0^\infty\mathcal R(t)\,\dd t
		\leq2\EE\sum_{i=1}^d\beta\bigl(
		M_{0,i}^2\mathsf I(V_{0,i})+K^2\mathsf J(V_{0,i})\bigr)
		=\mathcal R_0.
	\end{equation*}

	We now show that $\mathcal R_0$ is finite.
	If $\beta=0$, then $\mathcal R_0=0$.
	Suppose therefore that $\beta>0$, so $p=2\rho_1/\rho_2<1/2$.
	For $r>0$, replacing $(r-K^2)_+$ by $r$  and substituting $u=e^{-\rho_2t}r$, so that $\dd t=-\dd u/(\rho_2u)$ and, since $p=2\rho_1/\rho_2$, $e^{-2\rho_1t}=(u/r)^p$, gives
	\begin{equation*}
		\mathsf I(r)
		\leq\frac{1}{\rho_2r^p}\int_0^r\frac{u^p}{(u+\varepsilon)^{3/2}}\,\dd u.
	\end{equation*}
	Bounding $(u+\varepsilon)^{-3/2}$ by $\varepsilon^{-3/2}$, or extending the integral to $(0,\infty)$, where it converges because $p<1/2$, shows that $\mathsf I(r)$ is at most a constant multiple of $\min\{r,r^{-p}\}$.
	Using $r+\varepsilon\leq2\varepsilon$  for $r\leq\varepsilon$ and $r+\varepsilon\leq2r$  for $r>\varepsilon$, and $\mathsf I(0)=0$, we conclude that
	\begin{equation*}
		\mathsf I(r)\leq\frac{C}{(r+\varepsilon)^p},
		\qquad r\geq0,
	\end{equation*}
	where $C$ depends only on $\rho_1,\rho_2$, and $\varepsilon$.

	For the contribution involving $\mathsf J$, a uniform bound suffices.
	Again replace $(r-K^2)_+$ by $r$, and now set $u=e^{-\rho_2t}r$.
	Since $e^{-\rho_2t}r\,\dd t=-\dd u/\rho_2$, this gives, for $r>0$,
	\begin{align*}
		\mathsf J(r)
		&\leq\int_0^\infty\frac{e^{-\rho_2t}r}{(e^{-\rho_2t}r+\varepsilon)^{3/2}}\,\dd t\\
		&=\frac1{\rho_2}\int_0^r\frac{\dd u}{(u+\varepsilon)^{3/2}}
		=\frac2{\rho_2}\left(\frac1{\sqrt{\varepsilon}}-\frac1{\sqrt{r+\varepsilon}}\right)
		\leq\frac2{\rho_2\sqrt{\varepsilon}}.
	\end{align*}
	This bound also holds at $r=0$, where $\mathsf J(0)=0$.
	Thus
	\begin{equation*}
		\mathcal R_0
		\leq2\beta C\mathcal A_p(\mu_0)
		+\frac{4\beta dK^2}{\rho_2\sqrt{\varepsilon}}
		<\infty.
	\end{equation*}
	The initial energy $\FF(\mu_0)$ is finite by \Cref{asmp:initial-integrability} and the weighted moment comparison \eqref{eq:weighted-moment-comparison}, so $a_0<\infty$ as well.
	Since $\mathcal D_*\geq0$, integrating (i) now gives
	\begin{equation*}
		Q(\mu_t^x)\leq\FF(\mu_t)
		\leq\FF(\mu_0)+\int_0^t\mathcal R(s)\,\dd s
		\leq a_0,
	\end{equation*}
	which completes (ii).

	Finally, suppose that $Q(\mu_t^x)\geq Q_*>-\infty$ for every $t\geq0$.
	Adding the remaining error to the energy defines a locally absolutely continuous function
	\begin{equation*}
		H(t)\defeq\FF(\mu_t)+\int_t^\infty\mathcal R(s)\,\dd s.
	\end{equation*}
	By (i), $H'(t)\leq-\mathcal D_*(\mu_t)\leq0$ almost everywhere, while $H(t)\geq Q_*$.
	Thus $H$ has a finite limit, and $\FF(\mu_t)$ has the same limit because the tail integral tends to zero.
	Moreover, integrating this derivative bound gives
	\begin{equation*}
		\int_0^T\mathcal D_*(\mu_t)\,\dd t
		\leq H(0)-H(T)\leq a_0-Q_*.
	\end{equation*}
	Letting $T\to\infty$ proves (iii).
\end{proof}

\section{Proof of the qualitative convergence theorem:
\texorpdfstring{$\omega$}{Omega}-limit analysis}
\label{sec:omega}

We prove \cref{thm:main,cor:sublevel-convergence} by identifying complete limiting trajectories.

\Cref{prop:corrected-energy} shows that the corrected dissipation has
finite integral over $[0,\infty)$. Consequently, the dissipation accumulated
on any fixed time window tends to zero when that window is translated toward
infinity. A limit of the laws at one isolated time does not retain enough
information to exploit this observation, so we translate the entire
trajectory and use tightness and equicontinuity to extract a complete limiting
trajectory defined on $\RR$.
Lower semicontinuity gives zero corrected dissipation almost everywhere along the limiting trajectory.
We then show that every such trajectory is an equilibrium.

The main difficulty is that, at the boundary value $\delta=0$, zero
corrected dissipation does not immediately force the momentum to vanish.
\Cref{lem:complete-zero-momentum} resolves this degeneracy by combining the
zero-dissipation identities with the equations for the momentum and
second-moment variables. Once the momentum vanishes, the weak equation forces
the position law to be constant and critical, and the availability of the
trajectory at arbitrarily early times forces the second-moment variable to
vanish as well. It follows that every $\omega$-limit point belongs to the
equilibrium set $\mathcal S$. Finally, connectedness of the $\omega$-limit
set upgrades convergence to $\mathcal S$ into convergence to a single
equilibrium whenever the relevant critical set is totally disconnected.

\subsection{Compactness of time translates}

The next lemma supplies the compactness needed to take these limits.
Tightness allows us to extract limiting laws, while the common modulus of continuity allows us to extract an entire limiting trajectory.
The final assertion describes the support of the momentum and second-moment variables for such a limiting trajectory, showing that the contributions of the initial values decay, leaving bounds determined by the force alone.

\begin{lem}\label[lem]{lem:compactness}
	Let $Z$ be a global characteristic solution satisfying \Cref{asmp:objective,asmp:convergence}.
	Then:
	\begin{enumerate}[label=\textup{(\roman*)}]
		\item The family $(\mu_t)_{t\geq0}$ is tight in $\SCP(E)$.
		\item There is a nondecreasing function $\varpi\from[0,\infty)\to[0,2]$, with $\varpi(h)\to0$ as $h\downarrow0$, such that
		\begin{equation*}
			d_{\mathrm{BL}}(\mu_t,\mu_s)\leq\varpi(\lvert t-s\rvert),
			\qquad s,t\geq0.
		\end{equation*}
		\item Suppose that $t_n\to\infty$ and $\bar\mu\in C(\RR;\SCP(E))$ satisfy
		\begin{equation*}
			\sup_{s\in I}d_{\mathrm{BL}}(\mu_{t_n+s},\bar\mu_s)\longrightarrow0
		\end{equation*}
		for every compact interval $I\subset\RR$.
		The translated curves are defined on $I$ for all sufficiently large $n$.
		Then
		\begin{equation}\label{eq:limit-bounded-support}
			\bar\mu_s\left(
			\RR^d\times[-K,K]^d\times[0,K^2]^d
			\right)=1,
			\qquad s\in\RR.
		\end{equation}
	\end{enumerate}
\end{lem}

\begin{proof}
	We first prove (i) by bounding the three marginals separately.
	\Cref{lem:a-priori} and the triangle inequality give, almost surely for every $t\geq0$,
	\begin{equation*}
		\lvert M_t\rvert\leq e^{-\rho_1t}\lvert M_0\rvert+\sqrt dK,
		\qquad
		\lvert V_t\rvert\leq e^{-\rho_2t}\lvert V_0\rvert+\sqrt dK^2.
	\end{equation*}
	Thus $\lvert M_t\rvert$ and $\lvert V_t\rvert$ are bounded, uniformly in $t$, by almost surely finite random variables, so the momentum and second-moment marginals of $(\mu_t)_{t\geq0}$ are tight.
	The position marginals are tight by \Cref{asmp:trajectory-bounds}, and a family of laws on a product space with tight marginals is tight.
	This proves (i).

	For (ii), we estimate increments over $[s,s+h]$, with $s,h\geq0$.
	The weighted moment estimates \eqref{eq:weighted-moment-propagation} and \eqref{eq:weighted-moment-comparison} bound $\mathcal A_1(\mu_r)$ uniformly in $r$.
	Since the position equation gives $\EE\lvert\dot X_r\rvert^2=\mathcal A_1(\mu_r)$, integration and the Cauchy--Schwarz inequality yield
	\begin{equation*}
		\EE\lvert X_{s+h}-X_s\rvert
		\leq\int_s^{s+h}\EE\lvert\dot X_r\rvert\,\dd r
		\leq h\sup_{r\geq0}\mathcal A_1(\mu_r)^{1/2}
		\leq Ch.
	\end{equation*}
	Here and below, $C$ denotes a finite constant independent of $s$ and $h$.
	For the other two variables, \Cref{lem:variation-of-constants} gives
	\begin{align*}
		M_{s+h}-M_s
		&=(e^{-\rho_1h}-1)M_s
		+\rho_1\int_s^{s+h}e^{-\rho_1(s+h-r)}g_r\,\dd r,\\
		V_{s+h}-V_s
		&=(e^{-\rho_2h}-1)V_s
		+\rho_2\int_s^{s+h}e^{-\rho_2(s+h-r)}g_r^2\,\dd r.
	\end{align*}
	The bounds $\lvert g_r\rvert\leq\sqrt dK$ and $\lvert g_r^2\rvert\leq\sqrt dK^2$ therefore imply
	\begin{align*}
		\lvert M_{s+h}-M_s\rvert
		&\leq(1-e^{-\rho_1h})(\lvert M_s\rvert+\sqrt dK)
		\leq C(1-e^{-\rho_1h})(1+\lvert M_0\rvert),\\
		\lvert V_{s+h}-V_s\rvert
		&\leq(1-e^{-\rho_2h})(\lvert V_s\rvert+\sqrt dK^2)
		\leq C(1-e^{-\rho_2h})(1+\lvert V_0\rvert),
	\end{align*}
	where the last inequalities use the bounds from (i).

	The right-hand sides need not be integrable, so we use boundedness of the test functions in $d_{\mathrm{BL}}$ before taking expectations.
	For every $\varphi$ with $\lVert\varphi\rVert_\infty\leq1$ and $\Lip(\varphi)\leq1$,
	\begin{align*}
		\lvert\varphi(Z_{s+h})-\varphi(Z_s)\rvert
		&\leq\min\{2,\lvert Z_{s+h}-Z_s\rvert\}\\
		&\leq\lvert X_{s+h}-X_s\rvert
		+\min\{2,\lvert M_{s+h}-M_s\rvert\}
		+\min\{2,\lvert V_{s+h}-V_s\rvert\}.
	\end{align*}
	Taking expectations, inserting the increment bounds, and then taking the supremum over $\varphi$ yields
	\begin{equation}\label{eq:time-equicontinuity-bound}
		\begin{aligned}
			d_{\mathrm{BL}}(\mu_{s+h},\mu_s)
			&\leq Ch
			+\EE\min\{2,C(1-e^{-\rho_1h})(1+\lvert M_0\rvert)\}\\
			&\quad
			+\EE\min\{2,C(1-e^{-\rho_2h})(1+\lvert V_0\rvert)\}.
		\end{aligned}
	\end{equation}
	Let $\varpi(h)$ be the minimum of $2$ and the right-hand side of \eqref{eq:time-equicontinuity-bound}.
	Since $d_{\mathrm{BL}}\leq2$, this is still an upper bound for the distance.
	It is nondecreasing in $h$ and independent of $s$.
	Each truncated random variable is bounded by $2$ and tends to zero almost surely as $h\downarrow0$.
	Dominated convergence therefore gives $\varpi(h)\to0$.
	Taking $h=\lvert t-s\rvert$ and using symmetry of the distance proves (ii).

	To prove (iii), fix $s\in\RR$ and a sequence as in the statement.
	For all sufficiently large $n$, $t_n+s\geq0$.
	Define $\pi\from E\to E$ by
	\begin{equation*}
		\pi(x,m,v)\defeq\bigl(x,\max\{-K,\min\{m,K\}\},\min\{v,K^2\}\bigr),
	\end{equation*}
	with the maximum and minimum taken coordinatewise, so that $\pi$ takes values in the closed set $F\defeq\RR^d\times[-K,K]^d\times[0,K^2]^d$.
	By \Cref{lem:a-priori}, $(\lvert M_{t,i}\rvert-K)_+\leq e^{-\rho_1t}\lvert M_{0,i}\rvert$ and $(V_{t,i}-K^2)_+\leq e^{-\rho_2t}V_{0,i}$, so
	\begin{equation*}
		\lvert Z_t-\pi(Z_t)\rvert
		\leq e^{-\rho_1t}\lvert M_0\rvert+e^{-\rho_2t}\lvert V_0\rvert,
	\end{equation*}
	which tends to zero almost surely as $t\to\infty$.
	For every $\varphi$ with $\lVert\varphi\rVert_\infty\leq1$ and $\Lip(\varphi)\leq1$, the difference $\lvert\varphi(Z_t)-\varphi(\pi(Z_t))\rvert$ is at most $\min\{2,\lvert Z_t-\pi(Z_t)\rvert\}$.
	Dominated convergence therefore gives $d_{\mathrm{BL}}(\mu_t,\pi_\#\mu_t)\to0$, and hence $\pi_\#\mu_{t_n+s}\Longrightarrow\bar\mu_s$.
	Since each $\pi_\#\mu_{t_n+s}$ is supported in $F$, the Portmanteau theorem gives $\bar\mu_s(F)\geq\limsup_{n\to\infty}\pi_\#\mu_{t_n+s}(F)=1$, which is \eqref{eq:limit-bounded-support}.
\end{proof}

\subsection{Stability of the weak formulation}

We next establish the properties that pass to a limit of time translates.
The limiting curve must still solve the weak equation.
We also justify the moment tests used to identify its dynamics.

\begin{lem}\label[lem]{lem:weak-stability}
	Suppose that \Cref{asmp:force-regularity} holds.
	Let $I\subset\RR$ be a compact interval, let
	$\lambda^n,\lambda\in C(I;\SCP(E))$, and suppose that
	\begin{equation*}
		\sup_{s\in I}d_{\mathrm{BL}}(\lambda_s^n,\lambda_s)\longrightarrow0.
	\end{equation*}
	If each $\lambda^n$ satisfies the integrated weak equation \eqref{eq:weak-equation} on $I$, then $\lambda$ satisfies the same equation.
\end{lem}

\begin{proof}
	Fix $\varphi\in C_c^1(E)$ and set
	\begin{equation*}
		f_s^n(z)\defeq\nabla\varphi(z)\cdot b(z,(\lambda_s^n)^x),
		\qquad
		f_s(z)\defeq\nabla\varphi(z)\cdot b(z,\lambda_s^x).
	\end{equation*}
	These functions vanish outside the fixed compact support of $\nabla\varphi$.
	Since $G$ is bounded, there is a constant $C<\infty$ such that $\lvert f_s^n(z)\rvert,\lvert f_s(z)\rvert\leq C$ for every $n$, $s\in I$, and $z\in E$.

	We first show convergence of the spatial integrals at a fixed time $s\in I$.
	The convergence $\lambda_s^n\Longrightarrow\lambda_s$ also gives weak convergence of the position marginals.
	Joint continuity of $G$ therefore implies that $f_s^n$ converges continuously to $f_s$, in the sense of \Cref{lem:continuous-convergence}.
	Since these functions are bounded by $C$, \Cref{lem:continuous-convergence-bounded} gives
	\begin{equation*}
		\int_E f_s^n\,\dd\lambda_s^n
		\longrightarrow\int_E f_s\,\dd\lambda_s.
	\end{equation*}

	It remains to pass to the limit in the integrated weak equation on any interval $[a,b]\subset I$.
	The spatial integrals above are bounded in absolute value by $C$, so dominated convergence gives convergence of their integrals over $[a,b]$.
	The endpoint integrals of $\varphi$ converge by weak convergence of $\lambda_a^n$ and $\lambda_b^n$, since $\varphi$ is bounded and continuous.
	Passing to the limit on both sides of \eqref{eq:weak-equation} proves the claim.
\end{proof}

The weak formulation uses compactly supported test functions, whereas the moment identities below require  test functions without compact support.
The next lemma justifies these tests when the $(m,v)$ variables remain in a fixed compact set.

\begin{lem}\label[lem]{lem:extended-tests}
	Suppose that \Cref{asmp:force-regularity} holds.
	Let $I\subset\RR$ be a compact interval, and let $\lambda\in C(I;\SCP(E))$ satisfy \eqref{eq:weak-equation}.
	Suppose that there is a compact set $B\subset\RR^d\times[0,\infty)^d$ such that
	\begin{equation*}
		\lambda_s(\RR^d\times B)=1,
		\qquad s\in I.
	\end{equation*}
	Then \eqref{eq:weak-equation} also holds for every function $\varphi(x,m,v)=\psi_1(x)\psi_2(m,v)$ such that $\psi_1\in C_c^1(\RR^d)$ or $\psi_1\equiv1$, and $\psi_2$ is $C^1$ on a neighborhood of $B$.
	For each such $\varphi$, the map $s\mapsto\int_E\varphi\,\dd\lambda_s$ is Lipschitz on $I$.
\end{lem}

\begin{proof}
	Choose $\chi\in C_c^\infty(\RR^d\times\RR^d)$ equal to one on a neighborhood of $B$ and supported in a larger neighborhood on which $\psi_2$ is $C^1$.
	If $\psi_1\in C_c^1(\RR^d)$, then $\psi_1\chi\psi_2$ is an admissible test function in \eqref{eq:weak-equation}.
	It agrees with $\varphi$, together with its gradient, on $\RR^d\times B$, which contains the support of $\lambda_s$, so \eqref{eq:weak-equation} holds for $\varphi$.
	If $\psi_1\equiv1$, we apply this case to $\zeta_R(x)\psi_2(m,v)$, where $\zeta_R$ is a smooth cutoff equal to one on $\{\lvert x\rvert\leq R\}$, supported in $\{\lvert x\rvert\leq2R\}$, and satisfying $\lvert\nabla\zeta_R\rvert\leq C/R$.
	The drift is bounded on $\RR^d\times B$, so every term without a derivative of $\zeta_R$ converges by dominated convergence as $R\to\infty$.
	The term containing $\nabla\zeta_R$ is bounded by $C/R$ after integration against $\lambda_s$, uniformly for $s\in I$.
	It therefore vanishes in the limit.
	Finally, $\nabla\varphi$ and the drift are bounded on $\RR^d\times B$, so the integrand on the right-hand side of \eqref{eq:weak-equation} is bounded and $s\mapsto\int_E\varphi\,\dd\lambda_s$ is Lipschitz.
\end{proof}

\subsection{Complete trajectories with zero dissipation}

If $\delta>0$, zero corrected dissipation immediately forces the momentum to vanish.
At the boundary $\delta=0$, the dissipation can also vanish at states with nonzero momentum, so we must use the dynamics to rule them out.
Under the support condition below, such states have $v_i=K^2$ and $G_i=0$ in any coordinate with $m_i\neq0$.
The equation for $v_i$ then makes this variable decrease, giving a positive contribution to the derivative of the mean of $m_i^2(K^2-v_i)$.
The next proof shows that this contradicts the vanishing of that mean along a trajectory with zero dissipation.

\begin{lem}\label[lem]{lem:complete-zero-momentum}
	Suppose that \Cref{asmp:force-regularity} holds.
	Let $\lambda\from\RR\to\SCP(E)$ be a weakly continuous solution of
	\eqref{eq:weak-equation}, with $\mu$ there replaced by $\lambda$.
	Suppose also that \Cref{asmp:rate-condition} holds.
	Suppose that
	\begin{equation}\label{eq:complete-trajectory-box}
		\lambda_s\bigl(\RR^d\times[-K,K]^d\times[0,K^2]^d\bigr)=1,
		\qquad s\in\RR,
	\end{equation}
	and that $\mathcal D_*(\lambda_s)=0$ for almost every $s\in\RR$.
	Then
	\begin{equation}\label{eq:complete-zero-momentum}
		\lambda_s(\{m=0\})=1,
		\qquad s\in\RR.
	\end{equation}
\end{lem}

\begin{proof}
	The support condition \eqref{eq:complete-trajectory-box} and  \Cref{lem:extended-tests}, applied on compact intervals with $B=[-K,K]^d\times[0,K^2]^d$, imply that
	\begin{equation*}
		s\longmapsto\int_E m_i^2\,\lambda_s(\dd z)
	\end{equation*}
	is continuous for each coordinate $i$.
	It therefore suffices to show that these nonnegative moments vanish for almost every $s$.

	If $\delta>0$, every summand in \eqref{eq:corrected-dissipation} is nonnegative and vanishes only when the corresponding $m_i$ is zero.
	Thus $\mathcal D_*(\lambda_s)=0$ gives that $m=0$ holds $\lambda_s$-almost surely, for almost every $s$, as required.

	Suppose now that $\delta=0$, so $\beta K^2=\varepsilon>0$.
	In particular, $K>0$ and $\beta>0$.
	Fix a coordinate $i$.
	For almost every $s$, zero dissipation and the support bound $v_i\leq K^2$ give
	\begin{equation}\label{eq:boundary-zero-set}
		m_i^2\bigl(\beta(K^2-v_i)+\alpha G_i(\lambda_s^x,x)^2\bigr)=0
		\qquad\lambda_s\text{-almost surely}.
	\end{equation}
	Both terms in parentheses are nonnegative, with $\alpha,\beta>0$.
	Consequently, at such a time, $\lambda_s$-almost every point with $m_i\neq0$ satisfies $v_i=K^2$ and $G_i(\lambda_s^x,x)=0$.

	To rule out momentum on this set, consider the weighted moment
	\begin{equation*}
		J_i(s)\defeq\int_E m_i^2(K^2-v_i)\,\lambda_s(\dd z).
	\end{equation*}
	The support condition makes its integrand nonnegative, and \eqref{eq:boundary-zero-set} implies $J_i(s)=0$ for almost every $s$.
	By \Cref{lem:extended-tests} again, $J_i$ is continuous.
	Hence $J_i$ vanishes at every time.

	We next compute the derivative of this identically zero moment.
	By \Cref{lem:extended-tests}, we may test with $m_i^2(K^2-v_i)$ on every compact time interval.
	The derivatives of the test function with respect to $m_i$ and $v_i$ give, respectively, the two terms in
	\begin{equation}\label{eq:boundary-moment-derivative}
		\begin{aligned}
			J_i'(s)
			&=2\rho_1\int_E m_i\bigl(G_i(\lambda_s^x,x)-m_i\bigr)(K^2-v_i)\,\lambda_s(\dd z)\\
			&\quad+\rho_2\int_E m_i^2\bigl(v_i-G_i(\lambda_s^x,x)^2\bigr)\,\lambda_s(\dd z)
		\end{aligned}
	\end{equation}
	for almost every $s$.
	Choose a time at which both \eqref{eq:boundary-zero-set} and \eqref{eq:boundary-moment-derivative} hold.
	On $\{m_i=0\}$, both integrands in the derivative formula vanish.
	On $\{m_i\neq0\}$, the identities $v_i=K^2$ and $G_i(\lambda_s^x,x)=0$ show that the first integrand vanishes and the second equals $K^2m_i^2$, $\lambda_s$-almost surely.
	Since $J_i$ is identically zero, we obtain
	\begin{equation*}
		0=J_i'(s)=\rho_2K^2\int_E m_i^2\,\lambda_s(\dd z).
	\end{equation*}
	Because $\rho_2K^2>0$, the momentum moment vanishes for almost every $s$.
	The continuity established at the start of the proof extends this conclusion to every $s$.
	Applying the argument to each coordinate proves \eqref{eq:complete-zero-momentum}.
\end{proof}

\subsection{Identification of limit points}

We now combine the preceding results to show that every limiting law is an equilibrium.

\begin{lem}\label[lem]{lem:complete-limit}
	Let $Z$ be a global characteristic solution satisfying \Cref{asmp:objective,asmp:convergence}.
	Let $t_n\to\infty$ and suppose that $\mu_{t_n}$ converges weakly to $\bar\mu_0$.
	Then $\bar\mu_0\in\mathcal S$.
\end{lem}

\begin{proof}
	\emph{Extraction of a complete trajectory.}
	For each integer $T\geq1$, the translate
	\(\mu_s^n\defeq\mu_{t_n+s}\) is defined on $[-T,T]$ for all sufficiently large $n$.
	By \Cref{lem:compactness}, these curves take values in a common tight, hence relatively compact, subset of $\SCP(E)$ and have a common modulus of continuity in $d_{\mathrm{BL}}$.
	The metric Arzel\`a--Ascoli theorem \cite[Proposition~3.3.1]{AmbrosioGigliSavare2008} on each $[-T,T]$, whose limits inherit the common modulus of continuity, and a diagonal argument over $T=1,2,\ldots$ give a weakly continuous curve
	\begin{equation*}
		\bar\mu\from\RR\to\SCP(E)
	\end{equation*}
	such that $\mu_{t_n+s}\to\bar\mu_s$ weakly, locally uniformly in $s\in\RR$.
	In particular, $\bar\mu_0$ is the prescribed limit.
	By \eqref{eq:limit-bounded-support},
	$\bar\mu$ satisfies the support condition \eqref{eq:complete-trajectory-box}.

	Each shifted curve satisfies the integrated weak equation \eqref{eq:weak-equation}.
	By \Cref{lem:weak-stability}, $\bar\mu$ satisfies that equation on the whole real line.

	\emph{Vanishing dissipation and momentum.}
	For $a<b$ and $n$ large enough that $t_n+a\geq0$, finite dissipation in \eqref{eq:finite-corrected-dissipation} gives
	\begin{equation}\label{eq:vanishing-window-dissipation}
		0\leq\int_a^b\mathcal D_*(\mu_{t_n+s})\,\dd s
		=\int_{t_n+a}^{t_n+b}\mathcal D_*(\mu_r)\,\dd r
		\leq\int_{t_n+a}^\infty\mathcal D_*(\mu_r)\,\dd r
		\longrightarrow0.
	\end{equation}
	For each $s$, joint continuity of $G$ shows that the nonnegative integrand defining $\mathcal D_*(\mu_{t_n+s})$ in \eqref{eq:corrected-dissipation} converges continuously to the integrand defining $\mathcal D_*(\bar\mu_s)$.
	Hence \Cref{lem:continuous-convergence-lsc} gives
	\begin{equation*}
		\mathcal D_*(\bar\mu_s)
		\leq\liminf_{n\to\infty}\mathcal D_*(\mu_{t_n+s}).
	\end{equation*}
	Fatou's lemma in the time variable and \eqref{eq:vanishing-window-dissipation} give
	\begin{equation*}
		\int_a^b\mathcal D_*(\bar\mu_s)\,\dd s
		\leq\liminf_{n\to\infty}
		\int_a^b\mathcal D_*(\mu_{t_n+s})\,\dd s=0.
	\end{equation*}
	Hence $\mathcal D_*(\bar\mu_s)=0$ for almost every $s$.

	Applying \Cref{lem:complete-zero-momentum} to $\bar\mu$ gives
	\begin{equation}\label{eq:limit-zero-m}
		\bar\mu_s(\{m=0\})=1,
		\qquad s\in\RR.
	\end{equation}

	\emph{Constancy and criticality of the position law.}
	Let $\nu_s=\bar\mu_s^x$.
	We first show that integrals of position test functions are constant in time.
	For $\psi\in C_c^1(\RR^d)$, \Cref{lem:extended-tests} applies to the test function $\varphi(x,m,v)=\psi(x)$.
	The gradient of $\varphi$ has no momentum or second-moment components, so the corresponding components of the drift do not appear.
	The position component contributes
	\begin{equation*}
		-\sum_{i=1}^d\frac{m_i}{\sqrt{v_i+\varepsilon}}\,\partial_i\psi(x),
	\end{equation*}
	which vanishes $\bar\mu_s$-almost surely by \eqref{eq:limit-zero-m}.
	Thus the integrated weak equation gives
	\begin{equation*}
		\int_{\RR^d}\psi\,\dd\nu_t=\int_{\RR^d}\psi\,\dd\nu_s,
		\qquad \psi\in C_c^1(\RR^d),\quad s<t.
	\end{equation*}
	Thus, $\nu_s$ is independent of $s$, and we denote the common law by $\nu$.

	We next use the momentum equation to show that $\nu$ is critical.
	Fix a coordinate $i$ and test with
	\begin{equation*}
		\varphi(x,m,v)=\psi(x)m_i.
	\end{equation*}
	Both endpoint integrals vanish because $m=0$ holds $\bar\mu_s$-almost surely at every time.
	In the drift integral, the position component again vanishes on $\{m=0\}$.
	Differentiating the factor $m_i$ therefore gives the only remaining term,
	\begin{equation*}
		\rho_1\bigl(G_i(\nu,x)-m_i\bigr)\psi(x)
		=\rho_1G_i(\nu,x)\psi(x)
		\qquad\bar\mu_s\text{-almost surely}.
	\end{equation*}
	Its integral depends only on the position law $\nu$, which is constant in time.
	Consequently,
	\begin{equation*}
		0=\rho_1(t-s)\int_{\RR^d}G_i(\nu,x)\psi(x)\,\nu(\dd x),
		\qquad s<t.
	\end{equation*}
	Since $\rho_1(t-s)>0$, the spatial integral vanishes for every $\psi\in C_c^1(\RR^d)$.
	The boundedness of $G_i$ makes $G_i(\nu,\cdot)\nu$ a finite signed measure, and the vanishing of all these test integrals implies that this measure is zero.
	Hence $G_i(\nu,x)=0$ for $\nu$-almost every $x$.
	Applying this to each coordinate gives $G(\nu,x)=0$ for $\nu$-almost every $x$, so $\nu\in\mathcal C$.

	\emph{Vanishing of the second-moment variable.}
	Since the position law is critical, its forcing term $G_i(\nu,x)^2$ is zero almost surely.
	We will combine the resulting exponential decay with the uniform support bound at arbitrarily early times.
	For each coordinate, set
	\begin{equation*}
		u_i(s)\defeq\int_E v_i\,\bar\mu_s(\dd z).
	\end{equation*}
	Then $0\leq u_i(s)\leq K^2$ for every $s\in\RR$.
	By \Cref{lem:extended-tests}, testing with $v_i$ shows that $u_i$ is locally absolutely continuous.
	The vanishing forcing term reduces its equation to
	\begin{equation*}
		u_i'(s)=-\rho_2u_i(s)
	\end{equation*}
	for almost every $s$.
	Solving between $s-T$ and $s$ gives, for every $T>0$,
	\begin{equation*}
		0\leq u_i(s)=e^{-\rho_2T}u_i(s-T)\leq K^2e^{-\rho_2T}.
	\end{equation*}
	Letting $T\to\infty$ yields $u_i(s)=0$.
	Nonnegativity of $v_i$ now implies that $v=0$ almost surely under $\bar\mu_s$.
	Consequently,
	\begin{equation*}
		\bar\mu_s=\nu\otimes\delta_0\otimes\delta_0,
		\qquad s\in\RR,
	\end{equation*}
	and in particular $\bar\mu_0\in\mathcal S$.
\end{proof}

\subsection{Proof of the convergence theorem}

To obtain convergence to a single law, we use the connectedness of the $\omega$-limit set.

\begin{proof}[Proof of \Cref{thm:main}]
	By \Cref{lem:compactness} and Prokhorov's theorem \cite[Theorem~5.1.3]{AmbrosioGigliSavare2008}, the continuous curve $t\mapsto\mu_t$ has relatively compact range in $(\SCP(E),d_{\mathrm{BL}})$.
	Define its $\omega$-limit set by
	\begin{equation*}
		\omega(\mu_\cdot)
		\defeq\bigcap_{T\geq0}\overline{\{\mu_t:t\geq T\}}.
	\end{equation*}
	The topological argument in \cite[Theorem~A.3(a)--(b)]{CarrilloGvalaniWu2020} applies to any continuous curve with relatively compact range.
	It gives that $\omega(\mu_\cdot)$ is nonempty, compact, and connected, and that
	\begin{equation*}
		d_{\mathrm{BL}}(\mu_t,\omega(\mu_\cdot))\longrightarrow0.
	\end{equation*}
	Every point of $\omega(\mu_\cdot)$ is a subsequential limit as $t\to\infty$, so \Cref{lem:complete-limit} gives $\omega(\mu_\cdot)\subset\mathcal S$.
	The distance to $\mathcal S$ is therefore no larger than the distance to $\omega(\mu_\cdot)$, proving \textup{(i)}.

	To prove \textup{(ii)}, define
	\begin{equation*}
		\Lambda(\mu)\defeq\int_E\min\{1,\lvert m\rvert^2\}\,\mu(\dd x,\dd m,\dd v),
		\qquad
		\Gamma(\nu)\defeq\int_{\RR^d}\lvert G(\nu,x)\rvert^2\,\nu(\dd x).
	\end{equation*}
	Both functionals are continuous for weak convergence.
	This is immediate for $\Lambda$, whose integrand is bounded and continuous.
	For $\Gamma$, suppose that $\nu_n\Longrightarrow\nu$.
	Joint continuity of $G$ shows that $\lvert G(\nu_n,\cdot)\rvert^2$ converges continuously to $\lvert G(\nu,\cdot)\rvert^2$, and these functions are bounded by $dK^2$, so \Cref{lem:continuous-convergence-bounded} gives $\Gamma(\nu_n)\to\Gamma(\nu)$.
	Both functionals also vanish at equilibria since $\Lambda(\nu\otimes\delta_0\otimes\delta_0)=\Gamma(\nu)=0$ for every $\nu\in\mathcal C$.
	By tightness and \textup{(i)}, every sequence of times tending to infinity has a subsequence along which $\mu_t$ converges weakly to a law in $\mathcal S$.
	Hence
	\begin{equation*}
		\lim_{t\to\infty}\Lambda(\mu_t)=0,
		\qquad
		\lim_{t\to\infty}\Gamma(\mu_t^x)=0,
	\end{equation*}
	where the second limit is the limit of the force in \eqref{eq:momentum-force-limits}.
	By Markov's inequality, the first limit shows that $M_t\to0$ in probability.
	Since $(V_{t,i}+\varepsilon)^p\geq\varepsilon^p$, the same holds for $\sum_iM_{t,i}^2/(V_{t,i}+\varepsilon)^p$.
	Summing \eqref{eq:a-priori-weighted} over $i$ bounds this sum, uniformly in $t$, by a random variable that is integrable by \Cref{asmp:initial-integrability}.
	Dominated convergence therefore gives the momentum limit in \eqref{eq:momentum-force-limits}.
	By \eqref{eq:weighted-moment-comparison}, the same holds for $\mathcal A_r(\mu_t)$ for every $r\geq p$, in particular for the kinetic energy with $r=1/2$.

	Finally, \Cref{prop:corrected-energy} supplies a finite limit for $\FF(\mu_t)$.
	The vanishing kinetic energy proves \eqref{eq:objective-energy-limits}.

	For \textup{(iii)}, recall that $\omega(\mu_\cdot)$ is connected and contained in $\mathcal S$.
	The continuous position-marginal map sends $\omega(\mu_\cdot)$ into
	the set in \eqref{eq:intrinsic-critical-set}.
	If that set is totally disconnected, the connected image of
	$\omega(\mu_\cdot)$ is a singleton, say $\{\nu_\infty\}$.
	The position-marginal map is one-to-one on $\mathcal S$.
	Hence
	$\omega(\mu_\cdot)
	=\{\nu_\infty\otimes\delta_0\otimes\delta_0\}$.
	Convergence of the distance to $\omega(\mu_\cdot)$ now proves \eqref{eq:weak-convergence}.
\end{proof}

\begin{proof}[Proof of \Cref{cor:sublevel-convergence}]
	By \Cref{prop:corrected-energy}, $a_0<\infty$ and $Q(\mu_t^x)\leq a_0$ for every $t\geq0$.
	Thus the position laws lie in the tight sublevel set $\mathcal L_{a_0}$, which is nonempty since it contains $\mu_0^x$.
	Applying \Cref{lem:objective-lower-bound} to this sublevel set shows that $Q$ is bounded from below on $\SCP_1(\RR^d)$.
	This proves \Cref{asmp:trajectory-bounds}, so \Cref{thm:main} applies.

	Moreover, containment of the position laws in $\mathcal L_{a_0}$ gives
	\begin{equation*}
		\mathcal C\cap\overline{\{\mu_t^x:t\geq0\}}
		\subset\mathcal C\cap\overline{\mathcal L_{a_0}}
		=\mathcal C_{a_0}.
	\end{equation*}
	If $\mathcal C_{a_0}$ is totally disconnected, so is the set on the left.
	Assertion \textup{(iii)} of \Cref{thm:main} therefore gives convergence to $\nu_\infty\otimes\delta_0\otimes\delta_0$ with $\nu_\infty\in\mathcal C_{a_0}$.
\end{proof}

\section{Proofs of the mean-field
\texorpdfstring{\L{}ojasiewicz}{Lojasiewicz} convergence-rate results}
\label{sec:lojasiewicz-proof}

We prove \Cref{thm:kl-rates} in three steps.
First, a momentum--force correction to the energy supplies force dissipation in addition to momentum dissipation.
Second, the mean-field \L{}ojasiewicz inequality yields finite trajectory length and a scalar differential inequality for the modified energy.
Finally, we verify the dynamical hypotheses under  conditions \textup{(a)} and \textup{(b)} of \Cref{thm:kl-rates} and use boundedness of $V$ in case \textup{(b)} to obtain full-state Wasserstein rates.
Throughout this section, write
\begin{equation*}
	A(t)\defeq\mathcal A_{1/2}(\mu_t),\qquad
	B(t)\defeq\mathcal B(\mu_t),\qquad S(t)\defeq A(t)+B(t).
\end{equation*}
The norm $\lVert\cdot\rVert_2$ is taken in $L^2(\Omega;\RR^d)$.
Constants $c,C>0$ may change from line to line but are independent of time.

\subsection{A modified energy controlling momentum and force}

\begin{lem}\label[lem]{lem:kl-energy}
	Suppose that the objective and force satisfy \Cref{asmp:objective}, and let $Z$ be a global characteristic solution satisfying \Cref{asmp:convergence} and \Cref{asmp:kl-dynamics}.
	For sufficiently small $\eta>0$, the modified energy
	\begin{equation}\label{eq:kl-modified-energy}
		h(t)\defeq\FF(\mu_t)-q_\infty
		-\eta\EE\sum_{i=1}^d
		\frac{M_{t,i}g_{t,i}}{\sqrt{V_{t,i}+\varepsilon}},
		\qquad t\geq T_0,
	\end{equation}
	is locally absolutely continuous, nonnegative, and decreases to zero.
	It satisfies
	\begin{align}
		cS(t)\leq-h'(t)&\leq CS(t),\label{eq:kl-strong-dissipation}\\
		S(t)&\leq Ch(t),\label{eq:kl-S-by-energy}\\
		\left\lvert h(t)-\bigl(Q(\mu_t^x)-q_\infty\bigr)\right\rvert
		&\leq CS(t),\label{eq:kl-gap-comparison}
	\end{align}
	where the derivative bounds hold almost everywhere and the other bounds hold for every $t\geq T_0$.
	Moreover, either $h$ is identically zero or $h(t)>0$ for every $t\geq T_0$.
\end{lem}

\begin{proof}
	\textit{Step 1: estimate the energy correction.}
	Set
	\begin{equation*}
		J(t)\defeq\EE\sum_{i=1}^d
		\frac{M_{t,i}g_{t,i}}{\sqrt{V_{t,i}+\varepsilon}},
		\qquad \lvert J(t)\rvert\leq\sqrt{A(t)B(t)}.
	\end{equation*}
	Subtracting $\eta J$ adds force dissipation because the momentum equation contributes $\rho_1B$ to its derivative:
	\begin{align*}
		J'
		&=\rho_1B-\rho_1J
		+\EE\sum_{i=1}^d\frac{M_{t,i}\dot g_{t,i}}{\sqrt{V_{t,i}+\varepsilon}}\\
		&\quad-\frac{\rho_2}{2}\EE\sum_{i=1}^d
		\frac{M_{t,i}g_{t,i}(g_{t,i}^2-V_{t,i})}{(V_{t,i}+\varepsilon)^{3/2}}.
	\end{align*}
	To bound the last two terms, use the force-variation bound \eqref{eq:kl-force-variation} and boundedness of $g$:
	\begin{equation*}
		\EE\sum_i\frac{\lvert\dot g_{t,i}\rvert^2}{\sqrt{V_{t,i}+\varepsilon}}\leq CA,
		\qquad
		\left\lvert\frac{g_{t,i}^2-V_{t,i}}{V_{t,i}+\varepsilon}\right\rvert
		\leq1+\frac{K^2}{\varepsilon}.
	\end{equation*}
	Weighted Cauchy--Schwarz therefore bounds these terms by $CA$ and $C\sqrt{AB}$, respectively.
	Young's inequality gives
	\begin{equation*}
		J'\geq\rho_1B-C\sqrt{AB}-CA\geq\frac{\rho_1}{2}B-CA.
	\end{equation*}

	To control $S$, we differentiate $A$ and $B$.
	The product rule gives, for almost every $t\geq T_0$,
	\begin{align*}
		A'(t)&=2\rho_1\EE\sum_i\frac{M_{t,i}(g_{t,i}-M_{t,i})}{a_{t,i}}
		-\frac{\rho_2}{2}\EE\sum_i\frac{M_{t,i}^2(g_{t,i}^2-V_{t,i})}{a_{t,i}^3},\\
		B'(t)&=2\EE\sum_i\frac{g_{t,i}\dot g_{t,i}}{a_{t,i}}
		-\frac{\rho_2}{2}\EE\sum_i\frac{g_{t,i}^2(g_{t,i}^2-V_{t,i})}{a_{t,i}^3}.
	\end{align*}
	The preceding bounds on $\dot g$ and $(g_i^2-V_i)/(V_i+\varepsilon)$, together with weighted Cauchy--Schwarz, imply
	\begin{equation*}
		\lvert A'\rvert\leq C\bigl(A+\sqrt{AB}\bigr),
		\qquad \lvert B'\rvert\leq C\bigl(B+\sqrt{AB}\bigr).
	\end{equation*}
	Combining these with the formula for $J'$ and $\sqrt{AB}\leq(A+B)/2$ yields
	\begin{equation}\label{eq:kl-S-derivative}
		\lvert J'(t)\rvert+\lvert S'(t)\rvert\leq CS(t).
	\end{equation}

	To justify these differentiations, fix a compact interval $I=[a,b]\subset[T_0,\infty)$.
	Absolute continuity of $g$ in $L^2$ gives
	\begin{equation*}
		g_t=g_a+\int_a^t\dot g_s\,\dd s
		\qquad\text{in }L^2(\Omega;\RR^d).
	\end{equation*}
	Since $\dot g\in L^1(I\times\Omega)$, the right-hand side defines a version with almost surely absolutely continuous paths.
	Joint continuity of $G$ and continuity of $X$ and $\mu^x$ make the original process $g_t=G(\mu_t^x,X_t)$ continuous.
	Equality at rational times outside a common null set therefore identifies the two versions on $I$.
	The bounds above make the product-rule terms integrable: $A$ is bounded by \eqref{eq:weighted-moment-propagation} and \eqref{eq:weighted-moment-comparison}, and $B$ by the force bound.
	Fubini's theorem then gives local absolute continuity of $A,B,J$ and the stated derivative formulas.

	\textit{Step 2: combine momentum and force dissipation.}
	The exact dissipation formula \eqref{eq:dissipation} and coercivity assumption \eqref{eq:kl-coercivity} give $\kappa A\leq\mathcal D(\mu_t)\leq CA$.
	Together with the energy identity from \Cref{prop:dissipation}, Step 1 gives
	\begin{equation*}
		h'=-\mathcal D(\mu_t)-\eta J'
		\leq-(\kappa-\eta C)A-\frac{\eta\rho_1}{2}B.
	\end{equation*}
	Choose $\eta>0$ so that $\eta C\leq\kappa/2$.
	Then $-h'\geq cS$, while $\lvert\mathcal D\rvert+\eta\lvert J'\rvert\leq CS$ gives the upper dissipation bound in \eqref{eq:kl-strong-dissipation}.
	The momentum and force limits \eqref{eq:momentum-force-limits}, together with the weighted-moment comparison \eqref{eq:weighted-moment-comparison}, imply $A,B\to0$ and hence $J\to0$.
	The energy limit \eqref{eq:objective-energy-limits} now gives $h\to0$.
	Since $h$ is nonincreasing, it is nonnegative.

	\textit{Step 3: control the objective gap and dissipation by $h$.}
	Integrating $-h'\geq cS$ controls the time integral of $S$.
	The derivative bound \eqref{eq:kl-S-derivative} prevents $S$ from decreasing too quickly, so for $s\geq t$ we have $S(s)\geq e^{-C(s-t)}S(t)$.
	Consequently,
	\begin{equation*}
		h(t)\geq c\int_t^\infty S(s)\,\dd s
		\geq cS(t)\int_0^1e^{-Cu}\,\dd u.
	\end{equation*}
	This is the pointwise estimate $S\leq Ch$ in \eqref{eq:kl-S-by-energy}.
	The definition of $h$ also gives the objective-gap comparison \eqref{eq:kl-gap-comparison}:
	\begin{equation*}
		\left\lvert h-\bigl(Q(\mu_t^x)-q_\infty\bigr)\right\rvert
		\leq\frac{A}{2\rho_1}+\eta\sqrt{AB}\leq CS.
	\end{equation*}
	Finally, $-h'\leq CS\leq Ch$ gives $h(t)\geq h(T_0)e^{-C(t-T_0)}$.
	Thus $h$ stays positive if $h(T_0)>0$; otherwise nonnegativity and monotonicity give $h\equiv0$.
\end{proof}

\subsection{Finite length and quantitative rates}

\begin{proof}[Proof of \Cref{thm:kl-rates}  \textup{(i)} and \textup{(ii)} under \Cref{asmp:kl-dynamics}]
	\textit{Step 1: transfer the mean-field inequality to the modified energy.}
	Let $h$ be the modified energy \eqref{eq:kl-modified-energy} from \Cref{lem:kl-energy}.
	Its bounds on $S$ and the objective gap imply
	\begin{equation}\label{eq:kl-energy-control}
		\lvert Q(\mu_t^x)-q_\infty\rvert+A(t)+B(t)\leq Ch(t).
	\end{equation}
	Choose $T\geq T_0$ so that $\mu_t\in\mathcal U$, $\lvert Q(\mu_t^x)-q_\infty\rvert<r$, and $h(t)<\min\{1,2r\}$ for $t\geq T$.
	We derive the scalar estimates when $h>0$.
	If $h\equiv0$, then $S=0$ and the objective and length bounds below hold directly.

	The mean-field inequality controls the objective gap, whereas dissipation is expressed through $h$.
	To connect them, write $e(t)\defeq Q(\mu_t^x)-q_\infty$ and use the comparison $h\leq e+CS$ from \eqref{eq:kl-gap-comparison}.
	If $e\geq h/2$, then $0<e<r$, and the mean-field inequality \eqref{eq:mean-field-kl} and concavity of $\varphi$ give
	\begin{equation*}
		\varphi'(h/2)\sqrt S\geq\varphi'(e)\sqrt B\geq1.
	\end{equation*}
	If $e<h/2$, the correction must account for at least half of $h$, so $S\geq h/(2C)$ and $\sqrt S/\sqrt h\geq c$.
	The two cases therefore suggest combining $\varphi$ with a square root: for $0\leq s<2r$, set
	\begin{equation*}
		\Psi(s)\defeq2\varphi(s/2)+\sqrt s,
		\qquad \Psi'(h)\sqrt S\geq c.
	\end{equation*}
	The square-root term also covers $e\leq0$, so this argument requires neither a nonnegative objective gap nor a local minimum at level $q_\infty$.
	Using the dissipation bound $-h'\geq cS$ from \eqref{eq:kl-strong-dissipation}, we obtain
	\begin{equation*}
		-\frac{\dd}{\dd t}\Psi(h)=\Psi'(h)(-h')
		\geq c\Psi'(h)S\geq c\sqrt S.
	\end{equation*}
	The chain rule applies on finite intervals because $h$ stays positive.
	Integration and $\Psi(h(t))\to0$ give the tail estimate
	\begin{equation}\label{eq:kl-root-S-tail}
		\int_t^\infty\sqrt{S(s)}\,\dd s\leq C\Psi(h(t)),\qquad t\geq T.
	\end{equation}

	\textit{Step 2: deduce convergence from finite length.}
	The position equation and weighted-moment comparison \eqref{eq:weighted-moment-comparison} give $\lVert\dot X_t\rVert_2\leq\varepsilon^{-1/4}\sqrt{A(t)}$.
	The tail estimate \eqref{eq:kl-root-S-tail}, together with the moment bound \eqref{eq:weighted-moment-propagation} on $[0,T]$, therefore proves the finite-length assertion \eqref{eq:kl-finite-length}.
	By Tonelli's theorem, $\EE\int_0^\infty\lvert\dot X_t\rvert\,\dd t<\infty$, so almost every position path has a limit $X_\infty\in L^1$.
	Minkowski's inequality on finite intervals and then Fatou's lemma yield
	\begin{equation}\label{eq:kl-length-rate}
		\lVert X_t-X_\infty\rVert_2
		\leq\int_t^\infty\lVert\dot X_s\rVert_2\,\dd s
		\leq C\Psi(h(t)).
	\end{equation}
	Writing $\nu_\infty\defeq\Law(X_\infty)$, this coupling gives convergence in $\WW_1$.
	If $X_0\in L^2$, finite length gives $X_\infty\in L^2$ and convergence in $\WW_2$.
	The Lipschitz bound \eqref{eq:q-wasserstein-lipschitz} for $Q$ then gives $Q(\nu_\infty)=q_\infty$.

	To identify the full limiting state, joint continuity of $G$ gives $g_t\to G(\nu_\infty,X_\infty)$ almost surely.
	Boundedness of $g$ and the mean-square force limit \eqref{eq:momentum-force-limits} show that this limit is zero, so $\nu_\infty\in\mathcal C$.
	The variation-of-constants formulas in \Cref{lem:variation-of-constants} now give $M_t,V_t\to0$ almost surely: in each forcing integral, the contribution from a fixed initial interval decays exponentially, while the remainder is bounded by a tail supremum of $\lvert g\rvert$ or $\lvert g\rvert^2$.
	Letting the splitting time tend to infinity proves the claim; the homogeneous terms vanish because $M_0,V_0$ are finite almost surely.

	\textit{Step 3: extract the rates from a scalar differential inequality.}
	Step 1 gives $S\geq c/\Psi'(h)^2$, hence $-h'\geq c/\Psi'(h)^2$ whenever $h>0$.
	Integrating yields the general rate bound
	\begin{equation}\label{eq:kl-implicit-rate}
		\int_{h(t)}^{h(T)}\Psi'(s)^2\,\dd s\geq c(t-T).
	\end{equation}
	The left-hand side, viewed as a function of its lower endpoint, decreases continuously from infinity to zero, since $\Psi'(s)\geq1/(2\sqrt s)$.
	Inverting the scalar bound \eqref{eq:kl-implicit-rate} therefore gives an upper bound for $h(t)$.

	For the power desingularizing function \eqref{eq:kl-power}, put $\vartheta=\max\{\theta,1/2\}$.
	On $0<s\leq h(T)<1$,
	\begin{equation*}
		\Psi'(s)\leq Cs^{-\vartheta},\qquad
		\Psi(s)\leq Cs^{1-\vartheta},\qquad
		-h'\geq ch^{2\vartheta}.
	\end{equation*}
	When $\vartheta=1/2$, the last inequality gives exponential decay of $h$.
	When $\vartheta>1/2$, it gives $(h^{1-2\vartheta})'\geq c(2\vartheta-1)$ and hence $h(t)\leq C(1+t-T)^{-1/(2\vartheta-1)}$.
	Thus $h(t)\leq CE_\theta(t-T)$ and $\Psi(h(t))\leq CR_\theta(t-T)$, with the rates in \eqref{eq:kl-power-rates}.
	These bounds hold trivially if $h\equiv0$.
	The energy comparison \eqref{eq:kl-energy-control} and length estimate \eqref{eq:kl-length-rate} now prove the objective and position rates \eqref{eq:kl-objective-rate} and \eqref{eq:kl-position-rate}.
\end{proof}

\subsection{Checkable dynamical conditions and full-state rates}

\begin{proof}[Proof of \Cref{thm:kl-rates} under conditions \textup{(a)} and \textup{(b)}]
	\textit{Step 1: verify the dynamical assumptions.}
	Lipschitz continuity gives the force-variation bound; the parameter conditions give coercivity.
	For the former, apply the force Lipschitz estimate \eqref{eq:force-wasserstein-lipschitz} and couple $X_s,X_t$ on the solution space:
	\begin{align*}
		\lVert g_t-g_s\rVert_2
		&\leq L\left(\WW_1(\mu_t^x,\mu_s^x)+\lVert X_t-X_s\rVert_2\right)\\
		&\leq2L\int_s^t\lVert\dot X_a\rVert_2\,\dd a
		\leq2L\varepsilon^{-1/4}\int_s^t\sqrt{A(a)}\,\dd a.
	\end{align*}
	The position increments are square integrable by the weighted-moment bound \eqref{eq:weighted-moment-propagation}, even if $X_s$ is only integrable.
	Since $A$ is locally bounded, the displayed estimate makes $g$ locally absolutely continuous in $L^2$; differentiation at Lebesgue points gives the force-variation bound \eqref{eq:kl-force-variation} with $L_g=2L\varepsilon^{-1/4}$.

	In case \textup{(a)}, $\alpha<1$, and the coefficient in the exact dissipation \eqref{eq:dissipation} satisfies
	\begin{equation*}
		\frac{\varepsilon+(1-\alpha)v_i+\alpha G_i(\mu^x,x)^2}{v_i+\varepsilon}
		\geq1-\alpha.
	\end{equation*}
	This proves the coercivity bound \eqref{eq:kl-coercivity} with $\kappa=1-\alpha$ and $T_0=0$.
	In case \textup{(b)},  \Cref{lem:a-priori} gives
	\begin{equation*}
		(V_{t,i}-K^2)_+\leq e^{-\rho_2t}(\bar v_0-K^2)_+,\qquad
		V_{t,i}\leq\bar v\defeq\max\{K^2,\bar v_0\}.
	\end{equation*}
	Using $1-\alpha\geq-\beta$ in the dissipation coefficient, we obtain
	\begin{equation*}
		\varepsilon+(1-\alpha)V_{t,i}+\alpha g_{t,i}^2
		\geq\delta-\beta e^{-\rho_2t}(\bar v_0-K^2)_+.
	\end{equation*}
	Since $\delta>0$, this is at least $\delta/2$ after some time $T_0$.
	Dividing by $V_{t,i}+\varepsilon\leq\bar v+\varepsilon$ proves coercivity with $\kappa=\delta/(2(\bar v+\varepsilon))$.
	Thus \Cref{asmp:kl-dynamics} holds in both cases, and parts \textup{(i)} and \textup{(ii)} follow from their proof above.

	\textit{Step 2: pass from weighted estimates to full-state rates.}
	For the rest of the proof, assume case \textup{(b)}.
	The uniform bound on $V$ makes the weights comparable to constants:
	\begin{equation*}
		(\bar v+\varepsilon)^{-1/2}\EE\lvert g_t\rvert^2
		\leq B(t)\leq\varepsilon^{-1/2}\EE\lvert g_t\rvert^2.
	\end{equation*}
	The same comparison holds for every law supported on $\{v_i\leq\bar v\}$.
	Thus the unweighted mean-field inequality \eqref{eq:unweighted-mean-field-kl} implies the weighted inequality \eqref{eq:mean-field-kl} after multiplying $\varphi$ by $(\bar v+\varepsilon)^{1/4}$; the upper bound gives the converse after rescaling.
	Hence parts \textup{(i)} and \textup{(ii)} also hold under the unweighted inequality alone.

	Suppose now that $\varphi$ has the power form \eqref{eq:kl-power}.
	The weighted-moment rate \eqref{eq:kl-objective-rate} gives
	\begin{equation*}
		\EE\lvert M_t\rvert^2+\EE\lvert g_t\rvert^2
		\leq\sqrt{\bar v+\varepsilon}\,S(t)\leq CE_\theta(t-T).
	\end{equation*}
	The variable $V$ inherits this decay through its linear relaxation equation:
	\begin{equation*}
		\EE\sum_i V_{t,i}
		=e^{-\rho_2(t-T)}\EE\sum_i V_{T,i}
		+\rho_2\int_T^t e^{-\rho_2(t-s)}\EE\lvert g_s\rvert^2\,\dd s.
	\end{equation*}
	For the exponential rate, choose $c<\rho_2$ and integrate to obtain $\EE\sum_iV_{t,i}\leq CE_\theta(t-T)$.
	For the algebraic rate, split the convolution at $(t+T)/2$: on the first half the kernel contributes exponential decay, and on the second half $E_\theta(s-T)\leq CE_\theta(t-T)$.
	This gives the same bound.
	Since $\lvert V_t\rvert\leq\sum_iV_{t,i}$ and $\lvert V_t\rvert^2\leq\bar v\sum_iV_{t,i}$, the full-state moment estimate \eqref{eq:kl-state-moment-rates} follows.

	Finally, the coupling of $Z_t$ with $(X_\infty,0,0)$ gives
	\begin{equation*}
		\WW_1(\mu_t,\mu_\infty)
		\leq\lVert X_t-X_\infty\rVert_2+\lVert M_t\rVert_2+\lVert V_t\rVert_2
		\leq C\bigl(R_\theta(t-T)+E_\theta(t-T)^{1/2}\bigr).
	\end{equation*}
	We used the position rate \eqref{eq:kl-position-rate} and the moment estimate above.
	Because $E_\theta^{1/2}\leq R_\theta$, this is the Wasserstein rate \eqref{eq:kl-full-law-rate}.
	If $X_0\in L^2$, the same coupling bounds $\WW_2$: the initial weighted-moment assumption and boundedness of $V_0$ give $M_0\in L^2$, and finite position length ensures that both laws have second moments.
\end{proof}

\appendix

\section{Auxiliary results}\label{sec:auxiliary-proofs}
\begin{lem}[Variation of constants]\label[lem]{lem:variation-of-constants}
	Let $I\subset\RR$ be an interval, let $g\in L^2_{\mathrm{loc}}(I;\RR^d)$, and suppose that $M,V\from I\to\RR^d$ are locally absolutely continuous and satisfy
	\begin{equation*}
		\dot M_t=\rho_1(g_t-M_t),
		\qquad
		\dot V_t=\rho_2(g_t^2-V_t)
	\end{equation*}
	for almost every $t\in I$, with the square taken coordinatewise.
	Then, for every $s,t\in I$ with $s\leq t$,
	\begin{align*}
		M_t&=e^{-\rho_1(t-s)}M_s
		+\rho_1\int_s^te^{-\rho_1(t-r)}g_r\,\dd r,\\
		V_t&=e^{-\rho_2(t-s)}V_s
		+\rho_2\int_s^te^{-\rho_2(t-r)}g_r^2\,\dd r.
	\end{align*}
\end{lem}

\begin{proof}
	Multiplying the two equations by their integrating factors gives
	\begin{equation*}
		\frac{\dd}{\dd t}\bigl(e^{\rho_1t}M_t\bigr)=\rho_1e^{\rho_1t}g_t,
		\qquad
		\frac{\dd}{\dd t}\bigl(e^{\rho_2t}V_t\bigr)=\rho_2e^{\rho_2t}g_t^2
	\end{equation*}
	for almost every $t\in I$.
	Integrating from $s$ to $t$ and dividing by the respective exponential factors proves the formulas.
\end{proof}

\begin{lem}[A priori bounds]\label[lem]{lem:a-priori}
	Let $T\in(0,\infty]$, and let $g$, $M$, and $V$ be as in \Cref{lem:variation-of-constants} with $I=[0,T)$.
	Suppose that $V_0\in[0,\infty)^d$ and that $\lvert g_{t,i}\rvert\leq K$ for almost every $t\in I$ and every $1\leq i\leq d$.
	Then, for every $t\in I$ and $1\leq i\leq d$,
	\begin{equation}\label{eq:a-priori-momentum}
		\lvert M_{t,i}\rvert\leq e^{-\rho_1t}\lvert M_{0,i}\rvert+\bigl(1-e^{-\rho_1t}\bigr)K,
		\qquad
		\lvert M_t\rvert\leq e^{-\rho_1t}\lvert M_0\rvert+\bigl(1-e^{-\rho_1t}\bigr)\sqrt dK,
	\end{equation}
	and
	\begin{equation}\label{eq:a-priori-second-moment}
		e^{-\rho_2t}V_{0,i}\leq V_{t,i}\leq e^{-\rho_2t}V_{0,i}+\bigl(1-e^{-\rho_2t}\bigr)K^2,
		\qquad
		(V_{t,i}-K^2)_+\leq e^{-\rho_2t}(V_{0,i}-K^2)_+.
	\end{equation}
	Moreover, with $p$ as in \cref{sec:convergence-framework},
	\begin{equation}\label{eq:a-priori-weighted}
		\frac{M_{t,i}^2}{(V_{t,i}+\varepsilon)^p}
		\leq\frac{2M_{0,i}^2}{(V_{0,i}+\varepsilon)^p}
		+2K^2\varepsilon^{-p}.
	\end{equation}
\end{lem}

\begin{proof}
	By \Cref{lem:variation-of-constants} with $s=0$, the integral term in the formula for $M_t$ has coordinates of absolute value at most $(1-e^{-\rho_1t})K$ and norm at most $(1-e^{-\rho_1t})\sqrt dK$.
	The integral term in the formula for $V_{t,i}$ lies in $[0,(1-e^{-\rho_2t})K^2]$.
	This proves \eqref{eq:a-priori-momentum} and the first two inequalities in \eqref{eq:a-priori-second-moment}.
	Subtracting $K^2$ from the upper bound for $V_{t,i}$ gives $V_{t,i}-K^2\leq e^{-\rho_2t}(V_{0,i}-K^2)$, and hence the last inequality.

	For \eqref{eq:a-priori-weighted}, the momentum bound gives $M_{t,i}^2\leq2e^{-2\rho_1t}M_{0,i}^2+2K^2$.
	We divide the term containing $M_{0,i}$ by $(V_{t,i}+\varepsilon)^p\geq e^{-p\rho_2t}(V_{0,i}+\varepsilon)^p$ and the remaining term by $(V_{t,i}+\varepsilon)^p\geq\varepsilon^p$.
	This gives
	\begin{equation*}
		\frac{M_{t,i}^2}{(V_{t,i}+\varepsilon)^p}
		\leq2e^{-(2\rho_1-p\rho_2)t}
		\frac{M_{0,i}^2}{(V_{0,i}+\varepsilon)^p}
		+2K^2\varepsilon^{-p}.
	\end{equation*}
	Since $p\leq2\rho_1/\rho_2$, the exponential factor is at most one.
\end{proof}

Under \Cref{asmp:force-regularity}, \Cref{lem:a-priori} applies along almost every sample path of a global characteristic solution, with $g_t=G(\mu_t^x,X_t)$, and to the frozen equations in \cref{sec:well-posedness}, with $g_t=G(\eta_t,X_t)$.

The next lemma lets us pass to the limit in integrals when both the measures and the integrands vary.
We say that functions $f_n$ \textit{converge continuously} to $f$ on a set $S$ if $f_n(z_n)\to f(z)$ whenever $z_n\to z$ in $S$.

\begin{lem}\label[lem]{lem:continuous-convergence}
	Let $S$ be a closed subset of a Euclidean space, let $\lambda_n\Longrightarrow\lambda$ in $\SCP(S)$, and suppose that continuous functions $f_n\from S\to\RR$ converge continuously to a continuous function $f$ on $S$.
	\begin{enumerate}[label=\textup{(\roman*)},ref=\thethm(\roman*)]
		\item\label[lem]{lem:continuous-convergence-lsc} If every $f_n$ is nonnegative, then
		\begin{equation*}
			\int_Sf\,\dd\lambda
			\leq\liminf_{n\to\infty}\int_Sf_n\,\dd\lambda_n.
		\end{equation*}
		\item\label[lem]{lem:continuous-convergence-bounded} If $\sup_n\sup_{z\in S}\lvert f_n(z)\rvert<\infty$, then
		\begin{equation*}
			\int_Sf_n\,\dd\lambda_n\longrightarrow\int_Sf\,\dd\lambda.
		\end{equation*}
	\end{enumerate}
\end{lem}

\begin{proof}
	Continuous convergence gives $f(z)=\liminf_{n\to\infty,\,z'\to z}f_n(z')$ for every $z\in S$.
	Part (i) is therefore Fatou's lemma for weakly converging probabilities \cite[Theorem~1.1]{FeinbergKasyanovZadoianchuk2014}.
	For (ii), let $C$ bound every $\lvert f_n\rvert$, so that also $\lvert f\rvert\leq C$.
	Applying (i) to $C+f_n$ and to $C-f_n$, and using that all the measures are probability measures, gives
	\begin{equation*}
		\limsup_{n\to\infty}\int_Sf_n\,\dd\lambda_n
		\leq\int_Sf\,\dd\lambda
		\leq\liminf_{n\to\infty}\int_Sf_n\,\dd\lambda_n.
		\qedhere
	\end{equation*}
\end{proof}

\begin{lem}\label[lem]{lem:closed-critical}
	Under \Cref{asmp:force-regularity}, the sets $\mathcal C\subset\SCP(\RR^d)$ and $\mathcal S\subset\SCP(E)$ are closed for weak convergence.
\end{lem}

\begin{proof}
	Let $\nu_n\in\mathcal C$ and $\nu_n\Longrightarrow\nu$.
	Joint continuity of $G$ shows that $\min\{1,\lvert G(\nu_n,\cdot)\rvert\}$ converges continuously to $\min\{1,\lvert G(\nu,\cdot)\rvert\}$, so \Cref{lem:continuous-convergence-bounded} gives
	\begin{equation*}
		\int_{\RR^d}\min\{1,\lvert G(\nu_n,x)\rvert\}\,\nu_n(\dd x)
		\longrightarrow
		\int_{\RR^d}\min\{1,\lvert G(\nu,x)\rvert\}\,\nu(\dd x).
	\end{equation*}
	The integrals on the left vanish, so $\nu\in\mathcal C$.
	This proves that $\mathcal C$ is closed.
	If $\nu_n\otimes\delta_0\otimes\delta_0$ converges weakly in $\SCP(E)$, its limit has the same form with position marginal equal to the weak limit of $(\nu_n)$.
	The closedness of $\mathcal S$ follows.
\end{proof}

\begin{lem}\label[lem]{lem:objective-lower-bound}
	Under \Cref{asmp:objective},
	\begin{equation}\label{eq:q-wasserstein-lipschitz}
		\lvert Q(\nu)-Q(\mu)\rvert
		\leq\sqrt dK\WW_1(\mu,\nu),
		\qquad \mu,\nu\in\SCP_1(\RR^d).
	\end{equation}
	Moreover, if $\mathcal L_a$ is nonempty and tight for some $a\in\RR$, then $Q$ is bounded from below on $\SCP_1(\RR^d)$.
\end{lem}

\begin{proof}[Proof of \Cref{lem:objective-lower-bound}]
	Let $(X,Y)$ be any coupling of $\mu$ and $\nu$, and set
	$X_t=(1-t)X+tY$ for $0\leq t\leq1$.
	The chain rule and the coordinatewise force bound give
	\begin{align*}
		\lvert Q(\nu)-Q(\mu)\rvert
		&\leq\int_0^1
		\EE\bigl[\lvert G(\Law(X_t),X_t)\rvert\lvert Y-X\rvert\bigr]\dd t\\
		&\leq\sqrt dK\EE\lvert Y-X\rvert.
	\end{align*}
	Taking the infimum over all couplings proves \eqref{eq:q-wasserstein-lipschitz}.

	If $K=0$, the estimate shows that $Q$ is constant.
	Suppose that $K>0$ and, toward a contradiction, that $Q$ is not bounded from below.
	Choose $\nu_n\in\SCP_1(\RR^d)$ such that $Q(\nu_n)\to-\infty$.
	After omitting finitely many terms, $\nu_n\in\mathcal L_a$ and $Q(\nu_n)<0$.
	Hence $(\nu_n)$ is tight.
	Fix a unit vector $e$ and define
	\begin{equation*}
		r_n\defeq-\frac{Q(\nu_n)}{2\sqrt dK},
		\qquad
		\widetilde\nu_n\defeq(x\mapsto x+r_ne)_\#\nu_n.
	\end{equation*}
	For all sufficiently large $n$, \eqref{eq:q-wasserstein-lipschitz} gives
	\begin{equation*}
		Q(\widetilde\nu_n)
		\leq Q(\nu_n)+\sqrt dKr_n
		=\frac12Q(\nu_n)
		\leq a.
	\end{equation*}
	Thus $\widetilde\nu_n\in\mathcal L_a$.
	On the other hand, $r_n\to\infty$.
	To see that $(\widetilde\nu_n)$ is not tight, choose $R<\infty$ such that
	$\nu_n(B_R)\geq3/4$ for every $n$.
	For any compact set $C$, the translate $C-r_ne$ is disjoint from $B_R$ for all sufficiently large $n$, and hence
	$\widetilde\nu_n(C)=\nu_n(C-r_ne)\leq1/4$.
	This contradicts the tightness of $\mathcal L_a$.
\end{proof}

\section{Proof of well-posedness}\label{sec:well-posedness}

\begin{proof}[Proof of \Cref{prop:well-posedness}]
	We follow the fixed-point construction for McKean--Vlasov equations with Lipschitz coefficients \cite[Chapter~I]{Sznitman1991}, which freezes the curve of position laws, solves the resulting ordinary differential equation, and recovers the law as a fixed point.
	The two points that require care are that the drift is Lipschitz only for bounded momentum and that $V_0$ may not have any moments.

	Fix $T>0$ and $\eta\in C([0,T];\SCP_1(\RR^d))$, where continuity in $\SCP_1(\RR^d)$ is understood with respect to $\WW_1$.
	The \emph{frozen equation} is obtained from \eqref{eq:adam-system} by replacing $\mu_t^x$ with $\eta_t$.
	By \eqref{eq:force-wasserstein-lipschitz}, its vector field is continuous in time and locally Lipschitz in the state on $\RR^d\times\RR^d\times(-\varepsilon,\infty)^d$, so it has a unique maximal solution from each $z_0=(x_0,m_0,v_0)\in E$.
	On its interval of existence, \Cref{lem:a-priori} with $g_t=G(\eta_t,X_t)$ gives
	\begin{equation*}
		0\leq V_{t,i}\leq e^{-\rho_2t}v_{0,i}+K^2,
		\qquad
		\lvert M_t\rvert\leq B(z_0)\defeq\max\{\lvert m_0\rvert,\sqrt dK\},
	\end{equation*}
	and hence $\lvert X_t-x_0\rvert\leq B(z_0)t/\sqrt\varepsilon$.
	The solution therefore exists on all of $[0,T]$.
	If $z_0$ has law $\mu_0$, then $B(z_0)\leq B$ almost surely, which yields \eqref{eq:moment-explicit-bounds}.

	Since $G$ is bounded by $K$ and satisfies \eqref{eq:force-wasserstein-lipschitz}, and since $(m,v)\mapsto m/\sqrt{v+\varepsilon}$ is Lipschitz on $\{\lvert m\rvert\leq B\}\times[0,\infty)^d$, the drift satisfies
	\begin{equation}\label{eq:drift-lipschitz}
		\lvert b(z,\eta)-b(z',\xi)\rvert
		\leq C\bigl(\lvert z-z'\rvert+\WW_1(\eta,\xi)\bigr)
	\end{equation}
	whenever $\lvert m\rvert,\lvert m'\rvert\leq B$, where $C$ depends only on $B$, $K$, $L$, $\rho_1$, $\rho_2$, and $\varepsilon$.
	Let $Z^\eta$ and $Z^\xi$ be the frozen solutions driven by $\eta$ and $\xi$ from the same initial variable $Z_0$ with law $\mu_0$.
	Integrating \eqref{eq:drift-lipschitz} and applying Gronwall's inequality pathwise gives, for $0\leq t\leq T$,
	\begin{equation}\label{eq:frozen-stability}
		\sup_{0\leq s\leq t}
		\lvert Z_s^\eta-Z_s^\xi\rvert
		\leq Ce^{CT}\int_0^t\WW_1(\eta_s,\xi_s)\,\dd s.
	\end{equation}
	The right-hand side is deterministic, so this estimate requires no moment of $V_0$.

	Let $\mathcal X_T$ be the set of curves $\eta\in C([0,T];\SCP_1(\RR^d))$ with $\eta_0=\mu_0^x$, which is complete for the uniform $\WW_1$ metric, and define $\Phi_T(\eta)_t\defeq\Law(X_t^\eta)$.
	The bound $\lvert\dot X_t^\eta\rvert\leq B/\sqrt\varepsilon$ shows that $\Phi_T$ maps $\mathcal X_T$ into itself, and \eqref{eq:frozen-stability}, with the two outputs coupled through $Z_0$, gives
	\begin{equation*}
		\sup_{0\leq t\leq T}
		\WW_1\bigl(\Phi_T(\eta)_t,\Phi_T(\xi)_t\bigr)
		\leq Ce^{CT}T
		\sup_{0\leq t\leq T}\WW_1(\eta_t,\xi_t).
	\end{equation*}
	For $T=\tau$ sufficiently small, depending only on the constants in \Cref{prop:well-posedness}, Banach's fixed-point theorem gives a unique position-law curve on $[0,\tau]$.
	At time $\tau$, we again have $\lvert M_\tau\rvert\leq B$, $V_\tau\in[0,\infty)^d$, and $\EE\lvert X_\tau\rvert<\infty$, so the construction can be repeated from $Z_\tau$ with the same $\tau$.
	Since the solution on each interval is unique, the pieces combine into a global characteristic solution.

	For the fixed-point curve $\eta$, continuity of the frozen flow $(t,z_0)\mapsto\Phi^\eta_{t,0}(z_0)$ makes $(t,\omega)\mapsto\Phi^\eta_{t,0}(Z_0(\omega))$ jointly measurable, and the law $\mu_t=(\Phi^\eta_{t,0})_\#\mu_0$ depends only on $\mu_0$.
	Finally, if two solutions start from the same $Z_0$, then $\WW_1(\Law(X_t),\Law(X_t'))\leq\EE\lvert X_t-X_t'\rvert$, and \eqref{eq:drift-lipschitz} and Gronwall's inequality show that the two solutions coincide.
\end{proof}

\section{Proof of sharpness of the averaging-rate condition}\label{sec:sharpness}

In this appendix, we prove \cref{prop:sharpness}.
We construct the periodic trajectory before constructing the potential.
After normalizing the parameters and deriving the compatibility condition for a potential, we prescribe the momentum, recover the force and the periodic second-moment variables, and integrate the position equation to obtain a closed curve.
The first coordinate produces negative dissipation, which we balance against positive dissipation in the other two so that the total work over a period is zero.
These two coordinates also make the position curve embedded and provide a phase parameter for the construction of the potential.
We then extend the prescribed potential values and gradient from the curve to a confining potential on $\RR^3$, preserving the force bound, and verify the remaining conclusions of the proposition.

\begin{proof}[Proof of \cref{prop:sharpness}]
	\textit{Step 1: normalization and the compatibility condition.}\par\nopagebreak[4]

	We first construct $U$ and the periodic solution $Z$, and verify the remaining assertions at the end of the proof.
	It suffices to do this for $\rho_1=K=1$, as the parameters $(1,\rho_2/\rho_1,\varepsilon/K^2,1)$ in place of $(\rho_1,\rho_2,\varepsilon,K)$ also violate \eqref{eq:extended-rate-condition}.
	If $\widetilde U$ and $(\tilde x,\tilde m,\tilde v)$ provide the desired properties for these parameters, with constant $\widetilde C$, then
	\begin{equation*}
		U(y)\defeq\frac{K}{\rho_1}\widetilde U(\rho_1y),
		\qquad
		\bigl(x(t),m(t),v(t)\bigr)\defeq\left(\frac{\tilde x(\rho_1t)}{\rho_1},\,K\tilde m(\rho_1t),\,K^2\tilde v(\rho_1t)\right)
	\end{equation*}
	are as required for the original parameters.

	With $\rho_1=K=1$, the hypothesis becomes $\alpha>1+\varepsilon$.
	Set
	\begin{equation*}
		v_*\defeq\frac{\varepsilon}{\alpha-1}<1,
		\qquad\text{so that}\qquad
		\varepsilon+(1-\alpha)v+\alpha g^2=(\alpha-1)(v_*-v)+\alpha g^2.
	\end{equation*}
	Thus a coordinate in which the second-moment variable exceeds $v_*$ and the force is small contributes negatively to the dissipation \eqref{eq:dissipation}.

	Fix a period $P>0$ and a smooth $P$-periodic function $g\from\RR\to\RR^3$ that is \textit{antiperiodic}, meaning that $g(t+P/2)=-g(t)$.
	Let $m$ and $v$ be the unique $P$-periodic solutions of
	\begin{equation*}
		\dot m=g-m,
		\qquad
		\dot v=\rho_2(g^2-v),
	\end{equation*}
	namely
	\begin{equation*}
		m(t)=\int_{-\infty}^te^{-(t-s)}g(s)\,\dd s,
		\qquad
		v(t)=\rho_2\int_{-\infty}^te^{-\rho_2(t-s)}g(s)^2\,\dd s.
	\end{equation*}
	Then $m$ is antiperiodic, while $v$ has period $P/2$ because $g^2$ does.
	Hence the velocity $w\defeq-m/\sqrt{v+\varepsilon}$ is antiperiodic, and the curve
	\begin{equation*}
		x(t)\defeq\frac12\int_{t-P/2}^tw(s)\,\dd s
	\end{equation*}
	is $P$-periodic and antiperiodic, with $\dot x(t)=\frac12\bigl(w(t)-w(t-P/2)\bigr)=w(t)$.
	If a potential $U$ satisfies
	\begin{equation}\label{eq:prescribed-gradient}
		\nabla U(x(t))=g(t),
		\qquad t\in\RR,
	\end{equation}
	then $Z_t\defeq(x(t),m(t),v(t))$ solves \eqref{eq:adam-system} for $Q(\nu)=\int_{\RR^3}U\,\dd\nu$, whose force is $G(\nu,y)=\nabla U(y)$.
	Equivalently, we may prescribe a smooth, $P$-periodic, antiperiodic $m$ and set $g\defeq m+\dot m$, since $m$ is then the periodic solution for this force.

	If \eqref{eq:prescribed-gradient} holds, then $\frac{\dd}{\dd t}U(x(t))=g(t)\cdot\dot x(t)$, and periodicity of $x$ forces $\int_0^Pg\cdot\dot x\,\dd t=0$.
	To express this in terms of the dissipation, define, for $1\leq i\leq3$,
	\begin{equation}\label{eq:period-dissipation}
		D_i\defeq\int_0^P\frac{m_i^2\bigl((\alpha-1)(v_*-v_i)+\alpha g_i^2\bigr)}{(v_i+\varepsilon)^{3/2}}\,\dd t.
	\end{equation}
	This is the integral over one period of the $i$th summand in \eqref{eq:dissipation}, evaluated along the curve.
	The derivative of the weighted momentum term computed in the proof of \cref{prop:dissipation} gives
	\begin{equation*}
		\frac{\dd}{\dd t}\frac{m_i^2}{2\sqrt{v_i+\varepsilon}}
		=\frac{g_im_i}{\sqrt{v_i+\varepsilon}}
		-\frac{m_i^2\bigl((\alpha-1)(v_*-v_i)+\alpha g_i^2\bigr)}{(v_i+\varepsilon)^{3/2}}.
	\end{equation*}
	The left-hand side has zero integral over a period, and $\dot x_i=-m_i/\sqrt{v_i+\varepsilon}$, so
	\begin{equation}\label{eq:period-dissipation-work}
		D_i=\int_0^P\frac{g_im_i}{\sqrt{v_i+\varepsilon}}\,\dd t
		=-\int_0^Pg_i\dot x_i\,\dd t.
	\end{equation}
	The compatibility condition is therefore $D_1+D_2+D_3=0$.

	\textit{Step 2: negative dissipation in the first coordinate.}\par\nopagebreak[4]

	Throughout the remainder of the construction, $\eta\in(0,1)$ is a small parameter.
	We write $c$ and $C$ for positive constants that may change from line to line and may depend on $\alpha$, $\varepsilon$, and the functions and numbers fixed below, but not on $\eta$.
	A statement holds \textit{for small $\eta$} if it holds for every $\eta$ in some interval $(0,\eta_0)$, where $\eta_0>0$ depends only on the same quantities.

	In the first coordinate, the force consists of a rapid oscillation with amplitude close to one, followed by a weak pulse.
	The oscillation raises $v_1$ above $v_*$, but since it has mean zero, it produces only a small momentum.
	The pulse then produces momentum while $v_1$ is still above $v_*$ and the force is small, so the dissipation is negative.

	To construct the oscillation, choose a smooth $1$-periodic function $q$ with mean zero, $\max\lvert q\rvert<1$, and
	\begin{equation*}
		\ell\defeq\int_0^1q(s)^2\,\dd s>v_*.
	\end{equation*}
	For instance, $q(s)=a\tanh(N\sin(2\pi s))$ has these properties if $v_*<a^2<1$ and $N$ is sufficiently large, since $\int_0^1\tanh^2(N\sin(2\pi s))\,\dd s\to1$ as $N\to\infty$.
	Since $q$ has mean zero, it has a $1$-periodic primitive $\bar q$.

	To construct the time profiles of the oscillation and the pulse, set $\gamma\defeq(\ell-v_*)/8$ and choose $T_1>1$ such that $\ell e^{-\rho_2(T_1-1)}\leq4\gamma$.
	Let $\chi_0,\chi_1\from\RR\to[0,1]$ be smooth functions such that $\chi_0=1$ on $[1,T_1]$, $\chi_0$ has compact support in $(0,T_1+\gamma/\rho_2)$, and $\chi_1\not\equiv0$ has compact support in $(T_1+\gamma/\rho_2,T)$, where $T\defeq T_1+2\gamma/\rho_2$.
	The function
	\begin{equation*}
		\bar v(t)\defeq\rho_2\ell\int_0^te^{-\rho_2(t-s)}\chi_0(s)^2\,\dd s
	\end{equation*}
	is the second-moment variable started from zero with input $\ell\chi_0^2$, which is the local average of the squared oscillation constructed below.
	It satisfies $\bar v(T_1)\geq\ell\bigl(1-e^{-\rho_2(T_1-1)}\bigr)\geq v_*+4\gamma$.
	Moreover, $0\leq\bar v\leq\ell<1$, so $\lvert\bar v'\rvert=\rho_2\lvert\ell\chi_0^2-\bar v\rvert\leq\rho_2$.
	Hence
	\begin{equation}\label{eq:sharpness-charged}
		\bar v(t)\geq v_*+2\gamma,
		\qquad t\in[T_1,T],
	\end{equation}
	and this interval contains the support of $\chi_1$.

	Fix $P>2T$ and set $\omega\defeq2\pi/P$.
	On $[0,P/2]$, define
	\begin{equation*}
		m_1(t)\defeq\eta^2\chi_0(t)\bar q(t/\eta^2)+\eta\chi_1(t),
	\end{equation*}
	and extend $m_1$ to $\RR$ by $m_1(t+P/2)=-m_1(t)$.
	Since $m_1$ vanishes near $0$ and $P/2$, this extension is smooth, $P$-periodic, and antiperiodic.
	On $[0,P/2]$, the corresponding force $g_1=m_1+\dot m_1$ is
	\begin{equation*}
		g_1(t)=\chi_0(t)q(t/\eta^2)
		+\eta^2\bigl(\chi_0+\chi_0'\bigr)(t)\,\bar q(t/\eta^2)
		+\eta\bigl(\chi_1+\chi_1'\bigr)(t).
	\end{equation*}
	The first term is the oscillation and the last is the pulse.
	The supports of $\chi_0$ and $\chi_1$ are disjoint, $\lvert m_1\rvert\leq C\eta^2$ on the support of $\chi_0$, and $\lvert m_1\rvert+\lvert g_1\rvert\leq C\eta$ on the support of $\chi_1$.
	Hence, for small $\eta$,
	\begin{equation}\label{eq:pulse-bounds}
		\lvert g_1\rvert\leq k\defeq\frac{1+\max\lvert q\rvert}{2}<1,
		\qquad
		\lvert\dot g_1\rvert\leq C\eta^{-2},
		\qquad
		\lvert m_1\rvert\leq C\eta,
		\qquad
		\lvert g_1m_1\rvert\leq C\eta^2.
	\end{equation}
	Since $v_1$ is a weighted average of $g_1^2$, it follows that $0\leq v_1\leq k^2<1$.

	We next show that $v_1>v_*$ during the pulse.
	For $t\in[0,P/2]$, \cref{lem:variation-of-constants} gives
	\begin{equation*}
		v_1(t)=e^{-\rho_2t}v_1(0)+\rho_2\int_0^te^{-\rho_2(t-s)}g_1(s)^2\,\dd s.
	\end{equation*}
	On $[0,P/2]$, we have $g_1^2\geq\chi_0^2q(\cdot/\eta^2)^2-C\eta^2$, because the two sides differ only by the small middle term on the support of $\chi_0$.
	Dropping the nonnegative first term therefore gives
	\begin{equation*}
		v_1(t)\geq\rho_2\int_0^te^{-\rho_2(t-s)}\chi_0(s)^2q(s/\eta^2)^2\,\dd s-C\eta^2.
	\end{equation*}
	To compare the integral with $\bar v(t)$, let $\Psi$ be a $1$-periodic primitive of $q^2-\ell$, which exists because $q^2-\ell$ has mean zero.
	Since $q(s/\eta^2)^2-\ell=\eta^2\frac{\dd}{\dd s}\Psi(s/\eta^2)$, an integration by parts shows that the integral differs from $\bar v(t)$ by at most $C\eta^2$.
	Together with \eqref{eq:sharpness-charged}, this gives $v_1\geq v_*+\gamma$ on the support of $\chi_1$, for small $\eta$.

	On the support of $\chi_1$, we have $m_1=\eta\chi_1$ and $g_1=\eta(\chi_1+\chi_1')$, so, for small $\eta$,
	\begin{equation*}
		(\alpha-1)(v_*-v_1)+\alpha g_1^2
		\leq-(\alpha-1)\gamma+C\eta^2
		\leq-\frac{(\alpha-1)\gamma}{2}.
	\end{equation*}
	Since also $(v_1+\varepsilon)^{3/2}\leq(1+\varepsilon)^{3/2}$, the part of the integral \eqref{eq:period-dissipation} for $i=1$ over the support of $\chi_1$ is at most
	\begin{equation*}
		-\frac{(\alpha-1)\gamma\eta^2}{2(1+\varepsilon)^{3/2}}\int_\RR\chi_1(t)^2\,\dd t
		=-c\eta^2.
	\end{equation*}
	On the support of $\chi_0$, we have $m_1^2\leq C\eta^4$, so this part of the integral is at most $C\eta^4$ in absolute value, and $m_1=0$ on the rest of $[0,P/2]$.
	Since $m_1^2$, $g_1^2$, and $v_1$ have period $P/2$, the integral over $[P/2,P]$ is the same.
	Finally, the integrand in \eqref{eq:period-dissipation} is bounded by $Cm_1^2\leq C\eta^2$.
	Hence, for small $\eta$,
	\begin{equation}\label{eq:negative-dissipation}
		-C\eta^2\leq D_1\leq-c\eta^2.
	\end{equation}

	\textit{Step 3: balancing the dissipation.}\par\nopagebreak[4]

	In the remaining two coordinates, the momentum rotates with a small constant amplitude.
	The force is then small, the second-moment variables remain far below $v_*$, and the dissipation is positive.
	We identify $\RR^2$ with $\mathbb C$ and write $\hat y\defeq y_2+iy_3$ for $y\in\RR^3$, and similarly $\hat x$, $\hat m$, and $\hat g$.
	The dot denotes the Euclidean inner product in $\RR^2$, so that $\hat g\cdot\hat m=\operatorname{Re}(\hat g\overline{\hat m})$.
	For $B\in(0,1]$, prescribe
	\begin{equation*}
		\hat m(t)\defeq Be^{i\omega t},
		\qquad
		\hat g\defeq\hat m+\dot{\hat m}=(1+i\omega)\hat m.
	\end{equation*}
	Both functions are $P$-periodic and antiperiodic, since $e^{i\omega P/2}=-1$, and $\hat g\cdot\hat m=B^2$.
	The corresponding periodic second-moment variables satisfy $0\leq v_2,v_3\leq\lvert\hat g\rvert^2=(1+\omega^2)B^2$.
	Hence $\bigl\lvert(v_i+\varepsilon)^{-1/2}-\varepsilon^{-1/2}\bigr\rvert\leq CB^2$ for $i=2,3$, and the first identity in \eqref{eq:period-dissipation-work} gives
	\begin{equation*}
		\left\lvert D_2+D_3-\frac{PB^2}{\sqrt\varepsilon}\right\rvert\leq CB^4.
	\end{equation*}
	Since $\hat g$ is proportional to $B$, the variables $v_2$ and $v_3$ are proportional to $B^2$, so $D_2+D_3$ depends continuously on $B$.
	The first coordinate does not depend on $B$.
	By \eqref{eq:negative-dissipation} and the intermediate value theorem, for small $\eta$ there is $B$ with
	\begin{equation*}
		D_2+D_3=-D_1,
		\qquad
		c\eta\leq B\leq C\eta.
	\end{equation*}
	We fix this amplitude, so that $D_1+D_2+D_3=0$.

	\textit{Step 4: geometry of the closed curve.}\par\nopagebreak[4]

	The force $g=(g_1,\hat g)$ is smooth, $P$-periodic, and antiperiodic, so the construction above gives a closed curve $x=(x_1,\hat x)$.
	We next show that $\hat x$ winds once around a circle of radius comparable to $\eta$.
	Since $v_2,v_3\leq CB^2$ and $\lvert\hat m\rvert=B$,
	\begin{equation*}
		\bigl\lvert\dot{\hat x}(t)+\omega Re^{i\omega t}\bigr\rvert
		=\left\lvert\dot{\hat x}(t)+\frac{\hat m(t)}{\sqrt\varepsilon}\right\rvert
		\leq CB^3,
		\qquad
		R\defeq\frac{B}{\omega\sqrt\varepsilon}.
	\end{equation*}
	The difference $\hat x(t)-iRe^{i\omega t}$ is antiperiodic, so each of its real components vanishes somewhere in $[0,P]$.
	Its derivative is bounded by $CB^3$, so $\lvert\hat x(t)-iRe^{i\omega t}\rvert\leq CB^3$.
	Since $R$ is comparable to $\eta$, we have $B^3\leq C\eta^2R$.
	Hence, for small $\eta$, $\hat x$ does not vanish, and we may write $\hat x=\lvert\hat x\rvert e^{i\vartheta}$ with a smooth function $\vartheta$ satisfying
	\begin{equation}\label{eq:sharpness-winding}
		\bigl\lvert\lvert\hat x(t)\rvert-R\bigr\rvert\leq\frac R8,
		\qquad
		\vartheta'(t)\geq\frac\omega2,
		\qquad
		\vartheta(t+P)=\vartheta(t)+2\pi.
	\end{equation}
	For the second bound, $\vartheta'=\operatorname{Im}(\dot{\hat x}/\hat x)$, and the two estimates above give $\lvert\dot{\hat x}/\hat x-i\omega\rvert\leq C\eta^2$.
	For the third, $\vartheta(t+P)-\vartheta(t)$ is a multiple of $2\pi$, and it differs from $\omega P=2\pi$ by at most $C\eta^2$.

	\textit{Step 5: construction of the potential.}\par\nopagebreak[4]

	By \eqref{eq:pulse-bounds} and the bounds $\lvert\hat g\rvert,\lvert\dot{\hat x}\rvert\leq CB\leq C\eta$, the function $g\cdot\dot x$ is bounded by $C\eta^2$.
	Its integral over a period is $-(D_1+D_2+D_3)=0$ by \eqref{eq:period-dissipation-work}.
	Hence
	\begin{equation*}
		\Phi(t)\defeq\int_0^tg(s)\cdot\dot x(s)\,\dd s
	\end{equation*}
	is $P$-periodic and satisfies $\lvert\Phi\rvert+\lvert\Phi'\rvert\leq C\eta^2$.
	Up to an additive constant, $\Phi(t)$ is the value that the potential must take at $x(t)$.

	To find a potential with the prescribed gradient along the curve, we begin by constructing a bounded function $U_0$ near the curve.
	A point $y=(y_1,\hat y)\in\RR\times\mathbb C$ with $\hat y\neq0$ is assigned the time $\tau(\hat y)$ at which the phase of $\hat x$ equals the phase of $\hat y$, and $U_0(y)$ is a truncated first-order expansion of the potential around $x(\tau(\hat y))$.
	By \eqref{eq:sharpness-winding}, $\vartheta$ induces a smooth bijection from $\RR/P\mathbb Z$ onto $\RR/2\pi\mathbb Z$.
	Hence, for $\hat y\neq0$, there is a unique $\tau(\hat y)\in\RR/P\mathbb Z$ such that
	\begin{equation*}
		\vartheta(\tau(\hat y))\equiv\arg\hat y
		\pmod{2\pi}.
	\end{equation*}
	The map $\tau$ is smooth on $\mathbb C\setminus\{0\}$, satisfies $\tau(\hat x(t))=t$ modulo $P$, and, since $\lvert\nabla\arg\hat y\rvert=1/\lvert\hat y\rvert$,
	\begin{equation*}
		\lvert\nabla\tau(\hat y)\rvert
		=\frac{1}{\vartheta'(\tau(\hat y))\lvert\hat y\rvert}
		\leq\frac{2}{\omega\lvert\hat y\rvert}.
	\end{equation*}

	Choose smooth functions $\psi\from[0,\infty)\to[0,1]$ and $h\from\RR\to\RR$ such that $\psi$ vanishes outside $(1/2,2)$ and equals one on $[3/4,3/2]$, while $h(u)=u$ near zero, $0\leq h'\leq1$, and $h'=0$ outside $[-1,1]$.
	Set $h_\eta(u)\defeq\eta^4h(u/\eta^4)$, so that $\lvert h_\eta\rvert\leq\eta^4$.
	For $y\in\RR\times\mathbb C$ and $t\in\RR$, define
	\begin{align*}
		Y(y,t)&\defeq\Phi(t)+g_1(t)h_\eta\bigl(y_1-x_1(t)\bigr)+\hat g(t)\cdot\bigl(\hat y-\hat x(t)\bigr),\\
		U_0(y)&\defeq\psi\bigl(\lvert\hat y\rvert/R\bigr)\,Y\bigl(y,\tau(\hat y)\bigr),
	\end{align*}
	with $U_0(y)\defeq0$ when $\hat y=0$.
	Since $Y$ is $P$-periodic in $t$, $U_0$ is well defined, and it is smooth because $\psi(\lvert\hat y\rvert/R)$ vanishes for $\lvert\hat y\rvert\leq R/2$.
	The function $h_\eta$ replaces $y_1-x_1(t)$ by a bounded function.
	This keeps $U_0$ bounded in $y_1$ and makes the term containing the large derivative $\dot g_1$ small.

	Using $\Phi'=g\cdot\dot x=g_1\dot x_1+\hat g\cdot\dot{\hat x}$, we obtain
	\begin{equation*}
		\partial_tY(y,t)
		=\dot g_1(t)h_\eta(u)
		+g_1(t)\bigl(1-h_\eta'(u)\bigr)\dot x_1(t)
		+\dot{\hat g}(t)\cdot\bigl(\hat y-\hat x(t)\bigr),
		\qquad u\defeq y_1-x_1(t).
	\end{equation*}
	We write $\psi$ for $\psi(\lvert\hat y\rvert/R)$ and evaluate $Y$, $\partial_tY$, $g_1$, $\hat g$, and $u$ at $t=\tau(\hat y)$.
	With $\nabla_{\hat y}$ denoting the gradient with respect to $(y_2,y_3)$, the chain rule gives
	\begin{equation*}
		\partial_1U_0=\psi g_1h_\eta'(u),
		\qquad
		\nabla_{\hat y}U_0=Y\nabla\psi+\psi\bigl(\hat g+\partial_tY\,\nabla\tau\bigr).
	\end{equation*}
	At a point $y=x(t)$ of the curve, we have $\tau(\hat y)=t$ and $u=0$, so $h_\eta'(u)=1$ and $\partial_tY=0$.
	By \eqref{eq:sharpness-winding}, $\psi=1$ in a neighborhood of the curve, so $\nabla\psi=0$ there.
	Consequently,
	\begin{equation*}
		\nabla U_0(x(t))=g(t),
		\qquad t\in\RR.
	\end{equation*}

	We next bound the gradient of $U_0$.
	On the support of $\psi(\lvert\hat y\rvert/R)$, we have $R/2\leq\lvert\hat y\rvert\leq2R$, and therefore
	\begin{equation*}
		\lvert\hat y\rvert\geq c\eta,
		\qquad
		\bigl\lvert\hat y-\hat x(\tau(\hat y))\bigr\rvert\leq4R\leq C\eta,
		\qquad
		\lvert\nabla\tau\rvert+\lvert\nabla\psi\rvert\leq\frac C\eta.
	\end{equation*}
	Together with $\lvert\Phi\rvert\leq C\eta^2$, $\lvert h_\eta\rvert\leq\eta^4$, $\lvert\hat g\rvert+\lvert\dot{\hat g}\rvert\leq C\eta$, and \eqref{eq:pulse-bounds}, this gives
	\begin{equation*}
		\lvert Y\rvert\leq C\eta^2,
		\qquad
		\lvert\partial_tY\rvert
		\leq C\eta^{-2}\eta^4+\frac{\lvert g_1m_1\rvert}{\sqrt\varepsilon}+C\eta^2
		\leq C\eta^2.
	\end{equation*}
	Consequently,
	\begin{equation*}
		\lvert\partial_1U_0\rvert\leq\lvert g_1\rvert\leq k,
		\qquad
		\lvert\nabla_{\hat y}U_0\rvert\leq C\eta.
	\end{equation*}
	We now fix $\eta$ so small that all of the preceding estimates hold and $\lvert\nabla_{\hat y}U_0\rvert\leq1/2$.

	The function $U_0$ vanishes when $\lvert\hat y\rvert\geq2R$.
	Set $L_0\defeq\max\lvert x_1\rvert+1$.
	When $\lvert y_1\rvert\geq L_0$, we have $\lvert u\rvert\geq1\geq\eta^4$, so $h_\eta(u)$ does not depend on $y_1$ and $\partial_1U_0=0$.
	Thus, on each of the regions $\{y_1\geq L_0\}$ and $\{y_1\leq-L_0\}$, $U_0$ is a smooth function of $\hat y$ alone, with compact support.
	It follows that $U_0$ and all of its derivatives are bounded.

	We now add confining terms that vanish near the curve.
	Choose $L\geq\max\{L_0,2R\}$ and a smooth even function $\zeta\from\RR\to\RR$ with bounded second derivative such that $\zeta'=0$ on $[-L,L]$, $\lvert\zeta'\rvert\leq1$, and $\zeta'(u)=\operatorname{sign}u$ for $\lvert u\rvert\geq L+1$.
	Set
	\begin{equation*}
		U(y)\defeq U_0(y)+\zeta(y_1)+\zeta(y_2)+\zeta(y_3).
	\end{equation*}
	If $\zeta'(y_i)\neq0$, then $\lvert y_i\rvert>L$, and hence $\partial_iU_0(y)=0$.
	For $i=1$ this follows from $L\geq L_0$, and for $i=2,3$ from $\lvert\hat y\rvert\geq\lvert y_i\rvert>2R$.
	Consequently, $\lvert\partial_iU\rvert\leq1$ everywhere, with equality when $\lvert y_i\rvert\geq L+1$, so the supremum in the proposition equals $1=K$.
	The Hessian of $U$ is bounded.
	Since $U_0$ is bounded and $\zeta(u)\geq\lvert u\rvert-C$,
	\begin{equation*}
		U(y)\geq\lvert y_1\rvert+\lvert y_2\rvert+\lvert y_3\rvert-C\geq\lvert y\rvert-C.
	\end{equation*}
	Along the curve, $\lvert x_1\rvert<L_0\leq L$ and $\lvert x_2\rvert,\lvert x_3\rvert\leq\lvert\hat x\rvert<2R\leq L$, so the confining terms have zero gradient there.
	Hence \eqref{eq:prescribed-gradient} holds, and $Z_t=(x(t),m(t),v(t))$ is a deterministic global characteristic solution of \eqref{eq:adam-system}.
	It is $P$-periodic, and $\lvert m(t)\rvert\geq\lvert\hat m(t)\rvert=B>0$.
	For $d>3$, we add $\sum_{i=4}^d\zeta(y_i)$ to $U$ and set the remaining components of $x$, $m$, and $v$ equal to zero.
	Since $\zeta'(0)=0$, this gives a solution with the same properties.

	\textit{Step 6: verification of the conclusions.}\par\nopagebreak[4]

	We now return to the original parameters.
	The potential $U$ is smooth, its gradient is bounded and globally Lipschitz, and $U(y)\geq K\lvert y\rvert-C$.
	Thus $Q$ is an instance of \cref{ex:interaction} with $W=0$.
	It satisfies \cref{asmp:objective,asmp:force-lipschitz}, and its sublevel sets are tight, so in particular $\mathcal L_{a_0}$ is tight.
	The initial law $\delta_{Z_0}$ satisfies \cref{asmp:bounded-initial-momentum} and therefore \cref{asmp:initial-integrability}.
	The position orbit $\{\delta_{x(t)}:t\geq0\}$ is tight because $x$ is bounded, and $Q(\delta_{x(t)})=U(x(t))$ is bounded from below, which gives \cref{asmp:trajectory-bounds}.
	This proves that all hypotheses of \cref{thm:main,cor:sublevel-convergence} other than \cref{asmp:rate-condition} hold.

	The function
	\begin{equation*}
		t\longmapsto\mathcal A_p(\mu_t)=\sum_{i=1}^d\frac{m_i(t)^2}{(v_i(t)+\varepsilon)^p}
	\end{equation*}
	is continuous, periodic, and positive because $m(t)\neq0$, so it is bounded away from zero.
	Finally, $d_{\mathrm{BL}}(\delta_z,\delta_{z'})\leq\lvert z-z'\rvert$, so $t\mapsto\mu_t=\delta_{Z_t}$ is continuous and the orbit $\{\mu_t:0\leq t\leq P\}$ is compact.
	It is disjoint from $\mathcal S$, since every law in $\mathcal S$ is concentrated on $\{m=0\}$, and $\mathcal S$ is closed by \cref{lem:closed-critical}.
	A compact set and a disjoint closed set have positive distance, so $d_{\mathrm{BL}}(\mu_t,\mathcal S)$ is bounded away from zero.
\end{proof}

\section{Noisy dynamics}\label{sec:noise}

The preceding results describe convergence to equilibrium in the zero-noise system.
This appendix studies how persistent momentum noise changes that behavior.
Under the moment assumptions below, Brownian forcing injects energy at a rate bounded away from zero and prevents convergence to the zero-momentum equilibrium set.
The variance term $\sigma^2$ in the second-moment equation has a different effect and can be absorbed into $\varepsilon$ after translating $v$, as detailed below.
A symmetry argument then rules out a class of free energies for one-dimensional potentials.
The final subsection develops an exponential Foster--Lyapunov estimate in a separate, strongly convex setting.

\subsection{The noisy system}

Fix a momentum noise amplitude $s\geq0$ and a variance $\sigma^2\geq0$, and let $W$ be a $d$-dimensional Brownian motion independent of $Z_0$.
The noisy mean-field Adam system is
\begin{equation}\label{eq:noisy-system}
	\begin{aligned}
		\dd X_t
		&=-\frac{M_t}{\sqrt{V_t+\varepsilon}}\,\dd t,\\
		\dd M_t
		&=\rho_1\bigl(G(\mu_t^x,X_t)-M_t\bigr)\dd t+\rho_1s\,\dd W_t,\\
		\dd V_t
		&=\rho_2\bigl(G(\mu_t^x,X_t)^2+\sigma^2-V_t\bigr)\dd t,
	\end{aligned}
\end{equation}
where $\mu_t=\Law(Z_t)$ and $\sigma^2$ is added to each coordinate.
Formally, \eqref{eq:noisy-system} describes Adam with step size $h$ and minibatch gradient noise of covariance $\sigma^2\mathrm{Id}$, in which case $s^2=h\sigma^2$.
For stochastic differential equation approximations of discrete Adam and their scaling rules, see \cite{MalladiEtAl2022}.
We treat $s$ and $\sigma$ as independent parameters and assume $s>0$ except when discussing the deterministic reduction below.
A \textit{global solution} of \eqref{eq:noisy-system} is an $E$-valued process with continuous sample paths, adapted to a filtration for which $W$ is a Brownian motion, that satisfies \eqref{eq:noisy-system} in integral form.
We do not develop a well-posedness theory for \eqref{eq:noisy-system} in general.
When $G(\nu,x)$ is independent of $\nu$, bounded, and globally Lipschitz in $x$, as in \cref{ex:single-well}, \eqref{eq:noisy-system} is an ordinary stochastic differential equation whose coefficients are locally Lipschitz on a neighborhood of $E$ and grow at most linearly.
It then has a unique global strong solution for every initial law by the standard existence and localization argument for stochastic differential equations \cite[Section~5.2]{Oksendal2003}.

The variance term alone admits a direct reduction to the zero-noise system.
Indeed, if $s=0$ and $V_{0,i}\geq\sigma^2$ for every $i$, then $\widetilde V=V-\sigma^2$ takes values in $[0,\infty)^d$, and \eqref{eq:noisy-system} becomes \eqref{eq:adam-system} with $\varepsilon$ replaced by $\varepsilon+\sigma^2$.
The preceding results therefore apply whenever their hypotheses hold for the translated system; its equilibrium laws are translated by $\sigma^2$ in the $v$-variable.

The second-moment equation is still an ordinary differential equation along each sample path.
The variation-of-constants formula of \cref{lem:variation-of-constants}, with $g_t^2$ replaced by $g_t^2+\sigma^2$, gives
\begin{equation}\label{eq:noisy-v-bound}
	V_{t,i}\leq e^{-\rho_2t}V_{0,i}+\bigl(1-e^{-\rho_2t}\bigr)(K^2+\sigma^2),
	\qquad 1\leq i\leq d.
\end{equation}
For the momentum, the same formula applied to the drift gives
\begin{equation}\label{eq:noisy-momentum}
	\begin{aligned}
		M_t&=e^{-\rho_1t}M_0+\rho_1\int_0^te^{-\rho_1(t-r)}G(\mu_r^x,X_r)\,\dd r+N_t,\\
		N_t&\defeq\rho_1s\int_0^te^{-\rho_1(t-r)}\,\dd W_r.
	\end{aligned}
\end{equation}
Here $N_t$ is a centered Gaussian vector with independent coordinates of variance $\frac{\rho_1s^2}{2}(1-e^{-2\rho_1t})$.
Since the middle term is bounded by $\sqrt dK$, we obtain, for every $k\in\NN$,
\begin{equation}\label{eq:noisy-moment-bound}
	\sup_{t\geq0}\EE\lvert M_t\rvert^{2k}<\infty
	\qquad\text{whenever }\EE\lvert M_0\rvert^{2k}<\infty.
\end{equation}

\subsection{Energy injection}

The noisy counterpart of \cref{prop:dissipation} contains an additional dissipative term from $\sigma^2$ and a positive source from the quadratic variation of $M$.
For the latter, define
\begin{equation*}
	\mathcal J(\mu)\defeq\int_E\sum_{i=1}^d\frac{1}{\sqrt{v_i+\varepsilon}}\,\mu(\dd x,\dd m,\dd v),
	\qquad\mu\in\SCP(E).
\end{equation*}

The Brownian source in the energy balance remains uniformly positive in time for each fixed initial law.
The next proposition gives an explicit lower bound without imposing a moment condition on the initial adaptive state.

\begin{prop}\label[prop]{prop:noisy-energy}
	Suppose that \cref{asmp:objective} holds and that $Z$ is a global solution of \eqref{eq:noisy-system} with $\EE\lvert X_0\rvert<\infty$ and $\EE\lvert M_0\rvert^2<\infty$.
	Then $t\mapsto\FF(\mu_t)$ is locally absolutely continuous, and
	\begin{equation}\label{eq:noisy-energy}
		\frac{\dd}{\dd t}\FF(\mu_t)
		=-\mathcal D(\mu_t)-\alpha\sigma^2\mathcal A_{3/2}(\mu_t)+\frac{\rho_1s^2}{2}\mathcal J(\mu_t)
	\end{equation}
	for almost every $t\geq0$.
	Moreover, for every $t\geq0$,
	\begin{equation}\label{eq:noisy-source-lower-bound}
		\frac{\rho_1s^2}{2}\mathcal J(\mu_t)
		\geq\iota_*(\mu_0)
		\defeq\frac{\rho_1s^2}{2}\EE\sum_{i=1}^d
		\frac{1}{\sqrt{\max\{V_{0,i},K^2+\sigma^2\}+\varepsilon}}>0.
	\end{equation}
\end{prop}

\begin{proof}
	Since $a_{t,i}\geq\sqrt\varepsilon$, \eqref{eq:noisy-moment-bound} bounds $\EE\lvert\dot X_t\rvert^2$ uniformly in time.
	Hence $X$ satisfies the hypothesis of \cref{asmp:chain-rule}, and, as in the proof of \cref{prop:dissipation},
	\begin{equation*}
		\frac{\dd}{\dd t}Q(\mu_t^x)
		=-\EE\sum_{i=1}^d\frac{g_{t,i}M_{t,i}}{a_{t,i}}
	\end{equation*}
	for almost every $t\geq0$.
	The process $a_{t,i}^{-1}$ has absolutely continuous sample paths, while $M_{t,i}$ is a semimartingale with $\dd\langle M_i\rangle_t=\rho_1^2s^2\,\dd t$.
	Itô's product rule therefore gives
	\begin{equation*}
		\dd\left(\frac{M_{t,i}^2}{2\rho_1a_{t,i}}\right)
		=\left(
		\frac{M_{t,i}(g_{t,i}-M_{t,i})}{a_{t,i}}
		-\alpha\frac{M_{t,i}^2(g_{t,i}^2+\sigma^2-V_{t,i})}{a_{t,i}^3}
		+\frac{\rho_1s^2}{2a_{t,i}}
		\right)\dd t
		+\frac{sM_{t,i}}{a_{t,i}}\,\dd W_{t,i}.
	\end{equation*}
	The integrand of the stochastic integral is bounded by $s\lvert M_{t,i}\rvert/\sqrt\varepsilon$, so the stochastic integral is a martingale.
	Using $V_{t,i}/a_{t,i}^2\leq1$, the absolute value of the drift is bounded by $C(1+M_{t,i}^2)$, where $C$ depends only on $K,\sigma,\alpha,\rho_1,s$, and $\varepsilon$.
	Taking expectations and applying Fubini's theorem shows that $\mathcal A_{1/2}(\mu_t)/(2\rho_1)$ is locally absolutely continuous.
	Adding the derivative of $Q(\mu_t^x)$ cancels the force--momentum terms and gives \eqref{eq:noisy-energy}.
	By \eqref{eq:noisy-v-bound}, $V_{t,i}\leq\max\{V_{0,i},K^2+\sigma^2\}$ along every sample path on which the equations hold.
	Taking inverse square roots and expectations gives  the lower bound.
	Each summand defining $\iota_*(\mu_0)$ is positive almost surely and at most $\varepsilon^{-1/2}$.
	Since $s>0$, the constant is finite and strictly positive.
\end{proof}

Integrating \eqref{eq:noisy-energy} and using \eqref{eq:noisy-source-lower-bound} gives
\begin{equation*}
	\FF(\mu_t)+\int_0^t\bigl(\mathcal D(\mu_r)+\alpha\sigma^2\mathcal A_{3/2}(\mu_r)\bigr)\dd r
	\geq\FF(\mu_0)+\iota_*(\mu_0)t.
\end{equation*}

Thus bounded energy and finite accumulated dissipation cannot both hold.
In particular, the estimates used in \cref{sec:omega} cannot be carried over unchanged.

The obstruction also appears directly in the momentum: its time-averaged second moment stays positive for every bounded force.
Let
\begin{equation*}
	\mathcal M_0\defeq\bigl\{\mu\in\SCP(E):\mu(\{m=0\})=1\bigr\},
\end{equation*}
which contains the deterministic equilibrium set $\mathcal S$ and its translate by $\sigma^2$ in the $v$-variable.
The translated laws are equilibria of the variance-only system with $s=0$.
For $s>0$, no stationary law can belong to $\mathcal M_0$.
Indeed, stationarity and continuity would give $M_t=0$ for every $t$, contradicting the quadratic variation $\langle M_i\rangle_t=\rho_1^2s^2t$.
Under a fourth-moment assumption, the next lemma also shows that the distance of $\mu_t$ from $\mathcal M_0$ cannot tend to zero.
A stationary noisy law would instead have persistent momentum fluctuations.

\begin{lem}\label[lem]{lem:persistent-momentum}
	Suppose that \cref{asmp:force-regularity} holds and that $Z$ is a global solution of \eqref{eq:noisy-system} with $\EE\lvert M_0\rvert^2<\infty$.
	Set
	\begin{equation*}
		\kappa\defeq\frac{\sqrt{K^2+2\rho_1s^2}-K}{2}.
	\end{equation*}
	Then:
	\begin{enumerate}[label=\textup{(\roman*)}]
		\item For every $1\leq i\leq d$,
		\begin{equation}\label{eq:persistent-momentum}
			\liminf_{T\to\infty}\frac1T\int_0^T\EE M_{t,i}^2\,\dd t\geq\kappa^2>0.
		\end{equation}
		\item If $\EE\lvert M_0\rvert^4<\infty$, then $d_{\mathrm{BL}}(\mu_t,\mathcal M_0)$ does not converge to zero.
		In particular, $\mu_t$ does not converge to $\mathcal S$.
	\end{enumerate}
\end{lem}

\begin{proof}
	By Itô's formula,
	\begin{equation*}
		\dd M_{t,i}^2
		=\bigl(2\rho_1M_{t,i}g_{t,i}-2\rho_1M_{t,i}^2+\rho_1^2s^2\bigr)\dd t
		+2\rho_1sM_{t,i}\,\dd W_{t,i},
	\end{equation*}
	and the stochastic integral is a martingale by \eqref{eq:noisy-moment-bound}.
	Taking expectations, integrating over $[0,T]$, and dividing by $2\rho_1T$ yields
	\begin{equation*}
		\frac1T\int_0^T\bigl(\EE M_{t,i}^2-\EE[M_{t,i}g_{t,i}]\bigr)\dd t
		=\frac{\rho_1s^2}{2}-\frac{\EE M_{T,i}^2-\EE M_{0,i}^2}{2\rho_1T}.
	\end{equation*}
	The last term tends to zero by \eqref{eq:noisy-moment-bound}.
	Let $y_T^2$ denote the time average in \eqref{eq:persistent-momentum}.
	The Cauchy--Schwarz inequality, applied first in $\Omega$ and then in time, gives
	\begin{equation*}
		\frac1T\int_0^T\bigl\lvert\EE[M_{t,i}g_{t,i}]\bigr\rvert\,\dd t
		\leq\frac KT\int_0^T\bigl(\EE M_{t,i}^2\bigr)^{1/2}\dd t
		\leq Ky_T.
	\end{equation*}
	Hence $\liminf_{T\to\infty}(y_T^2+Ky_T)\geq\rho_1s^2/2$.
	Since $\kappa$ is the positive root of $y^2+Ky=\rho_1s^2/2$, this proves (i).

	For (ii), suppose toward a contradiction that $d_{\mathrm{BL}}(\mu_t,\mathcal M_0)\to0$.
	For $R>0$, the function $z\mapsto\min\{m_i^2,R\}$ is bounded and Lipschitz, and it vanishes almost surely under every law in $\mathcal M_0$.
	Hence $\EE\min\{M_{t,i}^2,R\}\to0$ for every $R>0$.
	By \eqref{eq:noisy-moment-bound} with $k=2$, the family $(M_{t,i}^2)_{t\geq0}$ is uniformly integrable, so $\EE M_{t,i}^2\to0$.
	This contradicts (i).
\end{proof}

The next example shows that $\FF$ need not be nonincreasing, even when $\rho_2\leq4\rho_1$.

\begin{ex}\label[ex]{ex:noisy-single-well}
	Consider the objective of \cref{ex:single-well}, for which $K=1$ and $\mathcal C=\{\delta_0\}$, and let $\rho_2\leq4\rho_1$.
	Start from the deterministic state $Z_0=(0,0,\sigma^2)$.
	When $s=0$, this state is stationary.
	Moreover, by the reduction above together with \cref{prop:well-posedness,cor:sublevel-convergence}, every solution with $\EE\lvert X_0\rvert<\infty$, bounded $M_0$, and $V_{0,i}\geq\sigma^2$ converges weakly to $\delta_{(0,0,\sigma^2)}$.

	Suppose now that $s>0$.
	Since $Q\geq0$ and $\mathcal A_{1/2}\geq0$, the initial energy $\FF(\mu_0)=0$ is the global minimum of $\FF$.
	By \eqref{eq:noisy-v-bound} and \eqref{eq:noisy-moment-bound} with $k=2$, the right-hand side of \eqref{eq:noisy-energy} is continuous in $t$.
	At $t=0$ it equals $d\rho_1s^2/(2\sqrt{\sigma^2+\varepsilon})$, because $M_0=0$.
	Consequently, $\FF(\mu_t)>\FF(\mu_0)$ for all sufficiently small $t>0$.
	By \cref{lem:persistent-momentum}, the time-averaged mean-square momentum remains at least $\kappa^2$ in each coordinate, and $\mu_t$ does not converge to $\delta_{(0,0,\sigma^2)}$.
\end{ex}

\subsection{An obstruction to symmetric free energies}

The increase of $\FF$ in \cref{ex:noisy-single-well} motivates adding an entropy term, which for kinetic Langevin dynamics produces a nonincreasing free energy.
To see why this works when the preconditioner is constant, let $G(\nu,x)=\nabla U(x)$, fix $a>0$, and replace $\sqrt{V_t+\varepsilon}$ by $a$ in \eqref{eq:noisy-system}.
Set
\begin{equation*}
	H(x,p)\defeq U(x)+\frac{\lvert p\rvert^2}{2\rho_1a}.
\end{equation*}
The pair $(X,P)=(X,-M)$ then solves
\begin{align*}
	\dd X_t&=\rho_1\nabla_pH(X_t,P_t)\,\dd t,\\
	\dd P_t&=-\rho_1\nabla_xH(X_t,P_t)\,\dd t-\rho_1^2a\nabla_pH(X_t,P_t)\,\dd t-\rho_1s\,\dd W_t.
\end{align*}
The friction coefficient multiplying $\nabla_pH$ is $\rho_1^2a$, so the corresponding temperature is $\theta_a\defeq s^2/(2a)$.
When $e^{-H/\theta_a}$ is integrable, it defines a stationary Gibbs law, which is invariant under the reflection $x\mapsto-x$ of the position alone when $U$ is even.
For a sufficiently regular density $f$ on $\RR^{2d}$, the free energy
\begin{equation*}
	\int_{\RR^{2d}}Hf\,\dd x\,\dd p
	+\theta_a\int_{\RR^{2d}}f\log f\,\dd x\,\dd p
\end{equation*}
is nonincreasing along the flow \cite{Villani2009}.
For the following class of potentials, noisy Adam has no stationary law with this symmetry.

Throughout the rest of this subsection, let $d=1$ and $G(\nu,x)=U'(x)$, where $U\in C^1(\RR)$ is even, $U'$ is bounded and globally Lipschitz, and $\lvert U'\rvert$ is strictly increasing on $[0,\infty)$.
The objective of \cref{ex:single-well} satisfies these conditions.
A \textit{stationary law} is a law $\pi\in\SCP(E)$ such that the solution of \eqref{eq:noisy-system} with $\Law(Z_0)=\pi$ satisfies $\mu_t=\pi$ for every $t\geq0$.
Define the reflection
\begin{equation*}
	\mathcal R(x,m,v)\defeq(-x,m,v).
\end{equation*}
Since $U'$ is odd, \eqref{eq:noisy-system} is invariant under the joint reflection $(x,m,v)\mapsto(-x,-m,v)$ when $W$ is replaced by $-W$, but it is not invariant under $\mathcal R$.

\begin{prop}\label[prop]{prop:no-symmetric-stationary}
	No stationary law $\pi$ of \eqref{eq:noisy-system} satisfies $\mathcal R_\#\pi=\pi$.
\end{prop}

\begin{proof}
	Let $\pi$ be a stationary law, let $Z$ be the solution with $\Law(Z_0)=\pi$, and write $a(v)\defeq\sqrt{v+\varepsilon}$.

	\textit{A priori bounds.}
	Since $\Law(V_t)=\Law(V_0)$ for every $t$, letting $t\to\infty$ in \eqref{eq:noisy-v-bound} gives $\pi(v\leq K^2+\sigma^2)=1$.
	Similarly, \eqref{eq:noisy-momentum} gives $\lvert M_t\rvert\leq e^{-\rho_1t}\lvert M_0\rvert+K+\lvert N_t\rvert$.
	Letting $t\to\infty$ and using stationarity, we obtain $\pi(\lvert m\rvert>r)\leq\mathbb P(\lvert N_\infty\rvert\geq r-K)$ for $r>K$, where $N_\infty$ is centered Gaussian with variance $\rho_1s^2/2$.
	In particular, $\int\lvert m\rvert\,\dd\pi<\infty$.

	\textit{A stationarity identity.}
	Let $\psi,\varphi\in C_c^\infty(\RR)$.
	Since $X$ and $V$ have absolutely continuous sample paths,
	\begin{equation*}
		\frac{\dd}{\dd t}\bigl(\psi(X_t)\varphi(V_t)\bigr)
		=-\frac{M_t}{a(V_t)}\psi'(X_t)\varphi(V_t)
		+\rho_2\bigl(U'(X_t)^2+\sigma^2-V_t\bigr)\psi(X_t)\varphi'(V_t)
	\end{equation*}
	for almost every $t$.
	The right-hand side is integrable because $\int\lvert m\rvert\,\dd\pi<\infty$.
	Integrating over $[0,1]$, taking expectations, and using $\EE\psi(X_1)\varphi(V_1)=\EE\psi(X_0)\varphi(V_0)$ gives
	\begin{equation}\label{eq:stationarity-identity}
		\int_E\left(
		-\frac{m}{a(v)}\psi'(x)\varphi(v)
		+\rho_2\bigl(U'(x)^2+\sigma^2-v\bigr)\psi(x)\varphi'(v)
		\right)\pi(\dd x,\dd m,\dd v)=0.
	\end{equation}

	\textit{Consequence of the symmetry.}
	Suppose, toward a contradiction, that $\mathcal R_\#\pi=\pi$, and let $\psi$ be even.
	Then $\psi'$ is odd, so the first integrand in \eqref{eq:stationarity-identity} changes sign under $\mathcal R$ and integrates to zero.
	Define the finite signed measure
	\begin{equation*}
		\nu_\psi(B)\defeq\int_E\psi(x)\bigl(U'(x)^2+\sigma^2-v\bigr)\mathbf 1_B(v)\,\pi(\dd x,\dd m,\dd v),
		\qquad B\subset\RR\text{ Borel}.
	\end{equation*}
	By \eqref{eq:stationarity-identity}, $\int\varphi'\,\dd\nu_\psi=0$ for every $\varphi\in C_c^\infty(\RR)$.
	A distribution with vanishing derivative is a constant multiple of Lebesgue measure, and $\nu_\psi$ is finite, so $\nu_\psi=0$.
	Since $U'$ is bounded and $\pi(v\leq K^2+\sigma^2)=1$, dominated convergence extends this to every bounded continuous even $\psi$.
	Let $\nu_1$ correspond to $\psi\equiv1$.
	Using $\nu_{(U')^2}(\RR)=0$ and $\int(\sigma^2-v)\,\nu_1(\dd v)=0$, we obtain
	\begin{equation*}
		\int_E\bigl(U'(x)^2+\sigma^2-v\bigr)^2\,\pi(\dd x,\dd m,\dd v)=0.
	\end{equation*}
	Hence $V_0=U'(X_0)^2+\sigma^2$ almost surely.

	\textit{Contradiction.}
	By stationarity, $V_t=U'(X_t)^2+\sigma^2$ almost surely for each $t$.
	Both sides are continuous in $t$, so, almost surely, this holds for every $t\geq0$.
	The second-moment equation then gives $\dot V_t=0$, so $\lvert U'(X_t)\rvert$ is constant.
	Since $\lvert U'\rvert$ is even and strictly increasing on $[0,\infty)$, $\lvert X_t\rvert$ is constant, and continuity of $X$ gives $X_t=X_0$ for every $t$.
	The position equation then forces $M_t=0$ for almost every $t$, and hence for every $t$ by continuity.
	The momentum equation now reads
	\begin{equation*}
		0=\rho_1U'(X_0)t+\rho_1sW_t,
		\qquad t\geq0,
	\end{equation*}
	which is impossible: the right-hand side has quadratic variation $\rho_1^2s^2t>0$ for $t>0$.
\end{proof}

Stationarity alone allows moving sample paths.
In the proof above, it is the reflection symmetry together with the adaptive second-moment equation that forces them to be constant, placing the stationary law in $\mathcal M_0$.

To state the consequence for Lyapunov functionals, we call $\Phi\from\SCP(E)\to(-\infty,\infty]$ a \textit{strict Lyapunov functional} for \eqref{eq:noisy-system} if, for every solution with $\Phi(\mu_0)<\infty$, the map $t\mapsto\Phi(\mu_t)$ is nonincreasing and is constant only when $\mu_t=\mu_0$ for every $t\geq0$.
This strictness is the property needed to identify stationary laws from constant-energy trajectories.

\begin{cor}\label[cor]{cor:no-free-energy}
	Let $\Phi\from\SCP(E)\to(-\infty,\infty]$ be convex, satisfy $\Phi(\mathcal R_\#\mu)=\Phi(\mu)$ for every $\mu\in\SCP(E)$, and attain a finite minimum.
	Then $\Phi$ is not a strict Lyapunov functional for \eqref{eq:noisy-system}.
\end{cor}

\begin{proof}
	Let $\pi$ be a minimizer of $\Phi$.
	By invariance, $\mathcal R_\#\pi$ is also a minimizer, and by convexity so is $\bar\pi\defeq\frac12(\pi+\mathcal R_\#\pi)$, which satisfies $\mathcal R_\#\bar\pi=\bar\pi$.
	Suppose, toward a contradiction, that $\Phi$ is a strict Lyapunov functional.
	Along the solution with $\Law(Z_0)=\bar\pi$, the function $t\mapsto\Phi(\mu_t)$ is nonincreasing and bounded below by $\Phi(\bar\pi)$, so it is constant.
	Strictness then makes $\bar\pi$ a stationary law, contradicting \cref{prop:no-symmetric-stationary}.
\end{proof}

\Cref{cor:no-free-energy} applies to the usual functionals of free-energy type.
For example, given a finite reference measure $\lambda_v$ on $[0,\infty)$, such as Lebesgue measure restricted to $[0,K^2+\sigma^2]$, a kinetic energy $\mathcal K\from\RR\times[0,\infty)\to\RR$, and $\theta>0$, consider the functional
\begin{equation*}
	\Phi(\mu)
	\defeq\int_EU\,\dd\mu+\int_E\mathcal K(m,v)\,\mu(\dd x,\dd m,\dd v)+\theta\int_E\varrho\log\varrho\,\dd x\,\dd m\,\lambda_v(\dd v)
\end{equation*}
if $\mu$ has density $\varrho$ with respect to $\dd x\,\dd m\,\lambda_v(\dd v)$, and $\Phi(\mu)=\infty$ otherwise.
It is convex and $\mathcal R$-invariant.
If the Gibbs density proportional to $e^{-(U+\mathcal K)/\theta}$ is normalizable and has finite free energy, it is a minimizer.
With $\mathcal K(m,v)=m^2/(2\rho_1\sqrt{v+\varepsilon})$, this is $\FF$ plus a multiple of the entropy.
Convexity and $\mathcal R$-invariance are preserved if the temperature $\theta$ is allowed to depend on $(m,v)$, or if further convex $\mathcal R$-invariant terms, such as moments of $\lvert x\rvert$ or a Fisher information, are added, so the corollary applies to each such functional that attains a finite minimum.
A convex strict Lyapunov functional with a finite minimum must therefore break position-reflection symmetry, for instance through a coupling between position and momentum such as the tracking-error term in \eqref{eq:noisy-lyapunov} below.

\subsection{A Foster--Lyapunov substitute}

Drift inequalities provide another way to control the noisy dynamics, following the Foster--Lyapunov approach to continuous-time Markov processes \cite{MeynTweedie1993}.
Suppose that $\rho_2<4\rho_1$, so that $\alpha<1$.
Each summand of $\mathcal D$ is then at least $(1-\alpha)m_i^2/\sqrt{v_i+\varepsilon}$, so $\mathcal D\geq(1-\alpha)\mathcal A_{1/2}$.
Since also $\mathcal J\leq d\varepsilon^{-1/2}$, \eqref{eq:noisy-energy} gives
\begin{equation}\label{eq:noisy-drift}
	\frac{\dd}{\dd t}\FF(\mu_t)
	\leq-(1-\alpha)\mathcal A_{1/2}(\mu_t)+\frac{d\rho_1s^2}{2\sqrt\varepsilon}
\end{equation}
under the hypotheses of \cref{prop:noisy-energy}.
Thus the kinetic part of the energy dissipates at a rate proportional to itself, while the injection is bounded.

What \eqref{eq:noisy-drift} lacks is dissipation in the position.
In a separate setting with a strongly convex potential, a correction involving the tracking error $m-\nabla U(x)$ provides this dissipation, and exponentiating the corrected energy yields a Foster--Lyapunov estimate.
Let $G(\nu,x)=\nabla U(x)+\chi R(\nu,x)$, where $R$ is bounded and Lipschitz, $U\in C^2(\RR^d)$ satisfies $U(0)=0$ and $\nabla U(0)=0$, and
\begin{equation*}
	m_U\,\mathrm{Id}\preceq\nabla^2U(x)\preceq L_U\,\mathrm{Id},
	\qquad x\in\RR^d,
\end{equation*}
for some $0<m_U\leq L_U<\infty$.
This force is unbounded and therefore lies outside \cref{asmp:force-regularity} and the preceding symmetry result.
For $c_e,c_S>0$ define
\begin{equation}\label{eq:noisy-lyapunov}
	F(z)\defeq1+U(x)+c_e\lvert m-\nabla U(x)\rvert^2
	+\frac{1}{2\rho_1}\sum_{i=1}^d\frac{m_i^2}{\sqrt{v_i+\varepsilon}}
	+c_S\sum_{i=1}^d\sqrt{v_i+\varepsilon}.
\end{equation}
The second and fourth terms form the integrand of $\FF$ when $Q(\nu)=\int U\dd\nu$.
For $\rho_2<4\rho_1$, choose $c_e>0$ sufficiently small, then $c_S>0$ sufficiently small relative to $c_e$, and finally $\gamma>0$ sufficiently small.
These choices give
\begin{equation*}
	\mathcal Le^{\gamma F}\leq-\lambda e^{\gamma F}+B,
\end{equation*}
where $\lambda>0$, $B<\infty$, and $\mathcal L$ is the generator of \eqref{eq:noisy-system} with the law in $G$ frozen.

We sketch the verification.
Write $e\defeq m-\nabla U(x)$, $a_i\defeq\sqrt{v_i+\varepsilon}$, $K_F\defeq\sum_im_i^2/a_i$, and $S_F\defeq\sum_ia_i$.
The chain rule in $X$ and $V$, It\^o's formula in $M$, and Young's inequality give
\begin{align*}
	\mathcal L\left(U+\frac{K_F}{2\rho_1}\right)
	&\leq-\frac{1-\alpha}{2}K_F+C,\\
	\mathcal L\lvert e\rvert^2
	&\leq-\rho_1\lvert e\rvert^2+CK_F+C,\\
	\mathcal LS_F
	&\leq-\frac{\rho_2}{2}S_F+C\bigl(K_F+\lvert e\rvert^2+1\bigr),
\end{align*}
where $C$ depends only on the model parameters and the bounds on $R$ and $\nabla^2U$.
With the above choice of $c_e$ and $c_S$, these bounds give $\mathcal LF\leq-cD_F+C$ for some $c>0$, where $D_F\defeq K_F+\lvert e\rvert^2+S_F$.
Since also $\lvert\nabla_mF\rvert^2\leq C(\lvert e\rvert^2+K_F)$, it follows that
\begin{equation*}
	\mathcal Le^{\gamma F}
	=\gamma e^{\gamma F}\left(\mathcal LF+\frac{\gamma\rho_1^2s^2}{2}\lvert\nabla_mF\rvert^2\right)
	\leq\gamma e^{\gamma F}\left(-\frac c2D_F+C\right)
\end{equation*}
for sufficiently small $\gamma>0$, uniformly over the frozen law.
The sublevel sets of $D_F$ are compact by the strong convexity of $U$, and the expression in parentheses is bounded above by a negative constant outside a sufficiently large one.
Enlarging $B$ on that compact set proves the drift inequality.

Existence, uniqueness, and exponential convergence to an invariant law for sufficiently weak interaction require further analysis, including a coupling argument.
We leave that analysis to future work.

For this choice of $F$, neither $F$ nor any power $F^q$, $q\geq1$, satisfies the same geometric drift inequality.
Indeed, take $\chi=0$ and $z_R\defeq(x,p,p^2+\sigma^2)$, where $\lvert x\rvert=R$ and $p=\nabla U(x)$.
At $z_R$, the tracking error and the drifts of $m$ and $v$ vanish, the It\^o terms are nonnegative, and each coordinate of the position moves with speed $\lvert p_i\rvert/\sqrt{p_i^2+\sigma^2+\varepsilon}<1$.
Since $\lvert p\rvert\leq L_UR$, it follows that
\begin{equation*}
	\mathcal LF(z_R)\geq-\sqrt dL_UR,
	\qquad
	F(z_R)\geq\frac{m_U}{2}R^2,
\end{equation*}
which rules out an inequality $\mathcal LF\leq-\lambda F+B$ with $\lambda>0$.
Since $r\mapsto r^q$ is convex for $q\geq1$, we have $\mathcal LF^q\geq qF^{q-1}\mathcal LF$, so the same holds for $F^q$.
The obstruction is the normalization of the step, which lets $U$ decrease only at a rate proportional to $R$, whereas $U$ itself is of order $R^2$.
By contrast, the identity for $\mathcal Le^{\gamma F}$ above shows that the exponential form needs only that the expression in parentheses be bounded above by a negative constant outside a compact set.
With a constant preconditioner $a$, the corresponding position contribution is $\mathcal LU=-\lvert p\rvert^2/a$, which is of order $-R^2$.

\section{Numerical experiments}\label{sec:numerics}

Our experiments illustrate the energy estimates, the role of the averaging-rate condition, and the convergence rates of \cref{thm:main,thm:kl-rates}.
We begin with an interacting-particle model, which allows us to compare convergence at finite particle count with the approach to a continuous critical law.
Two scalar potentials then isolate failure of convergence outside the parameter range and sharp algebraic decay within it.
Finally, a teacher--student network illustrates the energy identity and exponential convergence in a training problem.
The numerical parameters and fitting procedures are collected in \cref{sec:numerics-methods}.

\subsection{Particle approximation}\label{sec:numerics-setup}

Take $\Omega=\{1,\ldots,N\}$ with the uniform probability measure.
The empirical law $\mu_t=N^{-1}\sum_j\delta_{(X_t(j),M_t(j),V_t(j))}$ is a characteristic solution when each particle solves \eqref{eq:adam-system} with force $G(\mu_t^x,X_t(j))$.
Thus our theorems apply directly to this continuous-time particle system whenever their hypotheses hold.
This does not require taking a limit as $N\to\infty$, although the constants in the convergence-rate theorem may depend on $N$.

To approximate the dynamics, set $\beta_i\defeq1-\rho_ih$, where $h\max\{\rho_1,\rho_2\}\leq1$.
Writing $\nu^k\defeq N^{-1}\sum_j\delta_{X^k(j)}$ and $G_k(j)\defeq G(\nu^k,X^k(j))$, we use
\begin{equation}\label{eq:discrete-adam}
	\begin{aligned}
		M^{k+1}&=\beta_1M^k+(1-\beta_1)G_k,\\
		V^{k+1}&=\beta_2V^k+(1-\beta_2)G_k^2,\\
		X^{k+1}&=X^k-h\frac{M^{k+1}}{\sqrt{V^{k+1}+\varepsilon}}.
	\end{aligned}
\end{equation}
This is a first-order consistent scheme for \eqref{eq:adam-system}.
For $\nu=N^{-1}\sum_j\delta_{x_j}$ and $Q_N(x_1,\ldots,x_N)\defeq Q(\nu)$, we have $G(\nu,x_j)=N\nabla_{x_j}Q_N$.
Consequently, \eqref{eq:discrete-adam} is also Adam without bias correction, applied to $Q_N$ with learning rate $h$, memory variables $m=M/N$, $v=V/N^2$, and regularizer $\hat\varepsilon=\varepsilon/N^2$ inside the square root.
This scaling keeps the effective regularization fixed as $N$ varies.
We write $\mu_t$ for the empirical law of the numerical state at $t=kh$ as well.

\subsection{An interacting-particle model}\label{sec:numerics-meanfield}

We use a one-dimensional model in which confinement competes with repulsion.
For $a\geq1$ and $\varsigma>0$, let
\begin{equation}\label{eq:repulsive-objective}
	\begin{aligned}
		Q(\nu)&\defeq a\int_\RR\log\cosh x\,\nu(\dd x)
		-\frac12\int_{\RR\times\RR}\ell_\varsigma(x-y)\,\nu(\dd x)\nu(\dd y),\\
		\ell_\varsigma(z)&\defeq\lvert z\rvert-\varsigma\bigl(1-e^{-\lvert z\rvert/\varsigma}\bigr).
	\end{aligned}
\end{equation}
The force is
\begin{equation*}
	G(\nu,x)=a\tanh x-\int_\RR\operatorname{sign}(x-y)\bigl(1-e^{-\lvert x-y\rvert/\varsigma}\bigr)\,\nu(\dd y).
\end{equation*}
Its coordinatewise bound is $K=a+1$.
The bounded, Lipschitz derivatives give \cref{asmp:objective,asmp:force-lipschitz}, by the chain-rule argument of \cref{ex:interaction}.
Moreover,
\begin{equation*}
	Q(\nu)\geq(a-1)\int_\RR\lvert x\rvert\,\nu(\dd x)-a\log2,
\end{equation*}
so $Q$ has tight sublevel sets when $a>1$.
We fix $\varsigma=10^{-7}$ in the computations.

\textit{The corrected energy.}
The first experiment illustrates the part of \cref{prop:corrected-energy} that permits $\rho_2>4\rho_1$ and requires no moment of $V_0$.
We take
\begin{equation*}
	a=2,\qquad \rho_1=1,\qquad \rho_2=8,\qquad \varepsilon=10.
\end{equation*}
Then $\alpha=2$, $K=3$, and $\delta=\varepsilon-(\alpha-1)K^2=1>0$, so the averaging-rate condition holds.
The initial variables are independent, with $X_0\sim N(0,1)$, $M_0$ uniform on $[-3,3]$, and $V_0=100\,\Upsilon^{-2}$ for $\Upsilon$ uniform on $(0,1)$.
This distribution has $\EE V_0=\infty$, while the bounded momentum gives the required weighted integrability.
Each finite sample has a finite empirical mean, but samples from this law allow us to examine the effect of increasingly large initial second-moment variables.

\begin{figure}[t]
	\centering
	\includegraphics[width=0.58\linewidth]{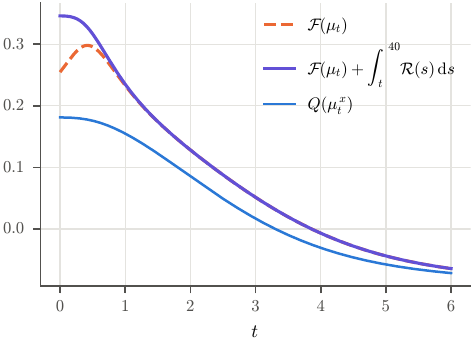}
	\caption{Corrected energy for the interacting-particle model with $\alpha=2$ and initial second-moment variables sampled from an infinite-mean distribution.
	The energy $\FF$ can increase even though the averaging-rate condition holds.
	Adding the remaining error gives the decreasing quantity used in \cref{prop:corrected-energy}.
	Here $N=10^4$; the tail integral is truncated as described in \cref{sec:numerics-methods}.}
	\label{fig:numerics-energy}
\end{figure}

\Cref{fig:numerics-energy} shows that $\FF$ initially increases, whereas the corrected energy
\begin{equation*}
	\FF(\mu_t)+\int_t^\infty\mathcal R(s)\,\dd s
\end{equation*}
is nonincreasing, as in the proof of \cref{prop:corrected-energy}.
Across particle counts from $10^2$ to $10^5$, the empirical mean of $V_0$ varies by several orders of magnitude, while the accumulated correction remains close to $0.09$ and below the bound $\mathcal R_0$.

\textit{Critical laws and finite-particle limits.}
For the remaining experiments with this model we use zero initial memory and parameters satisfying $\rho_2<4\rho_1$.
The critical law can be described explicitly in the limit $\varsigma\to0$.
At continuity points of the distribution function $F_\nu$, the limiting force is $a\tanh x-(2F_\nu(x)-1)$.
Its atomless critical law is
\begin{equation}\label{eq:repulsive-critical-law}
	\nu_\infty(\dd x)=\frac a2\operatorname{sech}^2x\,
	\mathbf 1_{\{\lvert\tanh x\rvert<1/a\}}\,\dd x.
\end{equation}
This law has compact support for $a>1$ and unbounded support for $a=1$.
Its midpoint quantiles are
\begin{equation*}
	x_j^*=\operatorname{arctanh}\frac{2j-1-N}{aN},\qquad
	\nu_{\infty,N}\defeq\frac1N\sum_{j=1}^N\delta_{x_j^*}.
\end{equation*}
They form an exact critical configuration for the limiting interaction and an approximate one for $\varsigma>0$, with a force error bounded by $e^{-\Delta/\varsigma}$, where $\Delta$ is the smallest spacing.
For the particle counts studied here, this error is below numerical precision at the limiting configurations.

Starting from independent Gaussian positions with $a=3/2$, the computed positions approach these quantiles, while momentum and second-moment variables approach zero.
This is the behavior described by \cref{thm:main}, with an empirical equilibrium law at each fixed $N$.
The remaining distance from $\nu_{\infty,N}$ to $\nu_\infty$ is a particle approximation error.

\textit{Dependence of the rate on the number of particles.}
To interpret the rates, let $\mathsf H\defeq N\nabla^2Q_N$ be the Jacobian of the particle force at a critical configuration.
At a strict local minimum with positive definite $\mathsf H$, the inequality \eqref{eq:unweighted-mean-field-kl} holds with $\theta=1/2$ on nearby empirical laws with the same number of particles.
For a trajectory approaching such a minimum with $\rho_2<4\rho_1$ and bounded initial memory, \cref{thm:kl-rates} therefore gives exponential convergence, with constants that may depend on $N$.
The linearization identifies the local decay rate: an eigenmode with eigenvalue $\lambda>0$ satisfies
\begin{equation*}
	\ddot y+\rho_1\dot y+\frac{\rho_1\lambda}{\sqrt\varepsilon}\,y=0,
\end{equation*}
since $V$ vanishes at equilibrium.
This is the heavy-ball equation \cite{Polyak1964}, with exponential envelope rate
\begin{equation}\label{eq:numerics-linear-rate}
	r(\lambda)\defeq\frac12\left(\rho_1-\operatorname{Re}\sqrt{\rho_1^2-4\rho_1\lambda/\sqrt\varepsilon}\right).
\end{equation}
This is the slow exponential rate of the linearized mode, and at critical damping a polynomial prefactor can appear.

\Cref{sec:width-dependence} gives the resulting smallest eigenvalue and the approximation for positive $\varsigma$.
The slowest particles are those in the tails, where confinement is weakest.
\Cref{fig:numerics-meanfield}(a) shows the resulting distinction: the measured rates remain bounded away from zero for $a>1$, but decrease like $1/N$ for $a=1$ over the tested range.
To compare the evolution of the laws, we also initialize the particles at midpoint quantiles of $N(0,1)$.
In \cref{fig:numerics-meanfield}(b), the curves for $a=3/2$ agree until they reach their respective particle approximation errors.
For $a=1$, relaxation is much slower, and the larger systems have not reached their particle approximation errors by the final time.

\subsection{Failure of convergence outside the parameter range}\label{sec:numerics-condition}

The single well of \cref{ex:single-well} isolates the role of the averaging-rate condition.
In dimension one its force is $G(x)=x/\sqrt{1+x^2}$, its only critical law is $\delta_0$, and its sublevel sets are tight.
With zero initial memory, \cref{cor:sublevel-convergence} therefore guarantees convergence whenever the parameter condition $\alpha\leq1+\varepsilon$ holds.

\begin{figure}[t]
	\centering
	\includegraphics[width=0.85\linewidth]{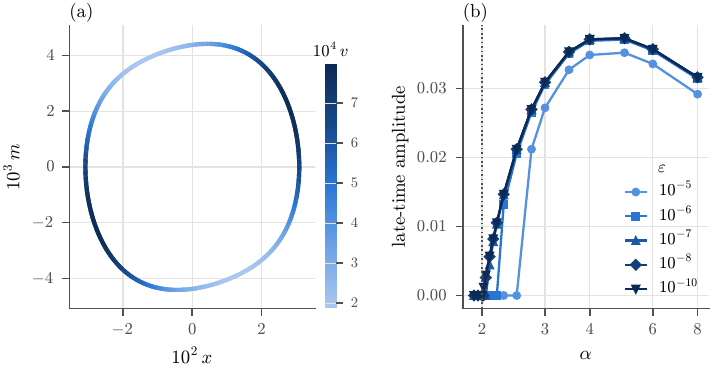}
	\caption{Single-well dynamics outside the sufficient parameter range.
	(a) The computed orbit for $\alpha=3$ and $\varepsilon=10^{-8}$, colored by the second-moment variable.
	(b) Late-time amplitude from $(2,0,0)$ as $\alpha$ varies.
	Amplitudes below the numerical threshold are shown on the axis; the dotted line marks $\alpha=2$.}
	\label{fig:numerics-orbit}
\end{figure}

For $\rho_1=1$, $\alpha=3$, and $\varepsilon=10^{-8}$, the trajectory from $(X_0,M_0,V_0)=(2,0,0)$ approaches the closed orbit in \cref{fig:numerics-orbit}(a).
Its amplitude and period stabilize; halving the integration step leaves the amplitude unchanged at the precision of the plot.
The dissipation changes sign along the computed orbit, allowing the energy lost during part of an oscillation to be recovered during the rest.
Small initial displacements instead converge to the equilibrium: among starts $X_0=c\sqrt\varepsilon$, those with $c=0.1,1$ converge to zero, while those with $c=10,100$ approach the same oscillation.
This gives numerical evidence for coexistence of an attracting equilibrium and an attracting periodic orbit.

\Cref{fig:numerics-orbit}(b) shows that the observed onset of persistent oscillations approaches $\alpha=2$ as $\varepsilon$ decreases.
The convergent runs extend beyond the sufficient condition $\alpha\leq1+\varepsilon$.
Thus the experiment supports the need for a restriction on the averaging rates, although for this potential convergence persists beyond the sufficient range.
This does not contradict \cref{prop:sharpness}, which concerns all objectives with force bound $K$.

\subsection{Sharp algebraic rates at a degenerate minimum}\label{sec:numerics-rates}

The next scalar example illustrates the algebraic bounds in \cref{thm:kl-rates}.
Let
\begin{equation*}
	Q(\nu)=\int_\RR\bigl((1+x^4)^{1/4}-1\bigr)\,\nu(\dd x),
	\qquad
	G(\nu,x)=\frac{x^3}{(1+x^4)^{3/4}}.
\end{equation*}
As in \cref{ex:single-well}, the force is bounded and Lipschitz and the potential is confining, but its curvature vanishes at the unique minimum.
For laws supported in $[-1,1]$, the bounds on the potential and force, followed by Jensen's inequality, give
\begin{equation*}
	Q(\nu)^{3/4}
	\leq\frac1{2^{3/2}}\left(\int_\RR x^6\,\nu(\dd x)\right)^{1/2}
	\leq\frac1{2^{3/4}}\left(\int_\RR G(\nu,x)^2\,\nu(\dd x)\right)^{1/2}.
\end{equation*}
This is \eqref{eq:unweighted-mean-field-kl} with $\theta=3/4$ and a constant independent of the number of particles.
With zero initial memory, $V_t\leq1$, so the weighted inequality holds as well.
For $\rho_2<4\rho_1$, every finite particle system eventually enters this neighborhood, and \cref{thm:kl-rates} gives
\begin{equation*}
	Q(\mu_t^x)\leq Ct^{-2},\qquad
	\lVert X_t\rVert_{L^2(\Omega)}\leq Ct^{-1/2}.
\end{equation*}

\begin{figure}[t]
	\centering
	\includegraphics[width=0.64\linewidth]{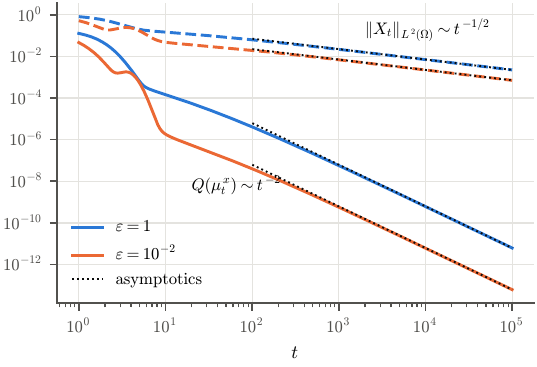}
	\caption{Algebraic convergence for the quartic well with $N=1000$.
	The objective (solid) and position norm (dashed) follow the powers $t^{-2}$ and $t^{-1/2}$ predicted by \cref{thm:kl-rates} for $\theta=3/4$.
	Dotted lines give the asymptotic expressions \eqref{eq:quartic-asymptotics}.}
	\label{fig:numerics-lojasiewicz}
\end{figure}

\Cref{fig:numerics-lojasiewicz} shows that these powers are attained in the computations.
Once $V_t\ll\varepsilon$ and the motion is slow compared with momentum relaxation, the reduced equation is $\dot X_t\approx-X_t^3/\sqrt\varepsilon$.
For particles following its nonzero solutions, this predicts
\begin{equation}\label{eq:quartic-asymptotics}
	Q(\mu_t^x)\sim\frac{\varepsilon}{16t^2},\qquad
	\lVert X_t\rVert_{L^2(\Omega)}\sim\left(\frac{\sqrt\varepsilon}{2t}\right)^{1/2}.
\end{equation}
For Gaussian initial positions, the computed quantities agree with these expressions to within one percent at the final time.
We note that the predicted powers agree with our theorem, while the leading constants come from the reduced equation.

\subsection{Energy and exponential rates in a teacher--student network}\label{sec:numerics-teacher}

We now examine the network objective \eqref{eq:neural-objective} from \cref{sec:neural-convergence}.
Recall that a neuron has parameters $x=(u,w)\in\RR\times\RR^8$, and for Gaussian inputs $\xi\sim N(0,I_8)$ we set
\begin{equation*}
	\begin{aligned}
		f_\nu(\xi)&\defeq\int_{\RR^9}2\tanh(u)\operatorname{erf}(w\cdot\xi)\,\nu(\dd u,\dd w),\\
		Q(\nu)&\defeq\frac12\EE\lvert f_\nu(\xi)-y(\xi)\rvert^2
		+\gamma\int_{\RR^9}\bigl(\sqrt{1+\lvert x\rvert^2}-1\bigr)\,\nu(\dd x).
	\end{aligned}
\end{equation*}
The teacher is $y(\xi)=\sum_{i=1}^3a_i^*\operatorname{erf}(w_i^*\cdot\xi)$, with $w_i^*=3e_i/2$ and $a^*=(0.6,-0.5,0.4)$.
We take $\gamma=10^{-2}$.
The regularity and sublevel hypotheses are verified in the proof of \cref{cor:neural-convergence}.
For zero initial memory and the stated averaging-rate condition, that corollary gives approach to the critical laws and decay of the mean-square force.
The rate comparisons below concern the local behavior at the computed finite-width minima.
The Gaussian expectations are evaluated in closed form \cite{Williams1997}, so no sampling noise enters the force.

\begin{figure}[t]
	\centering
	\includegraphics[width=\linewidth]{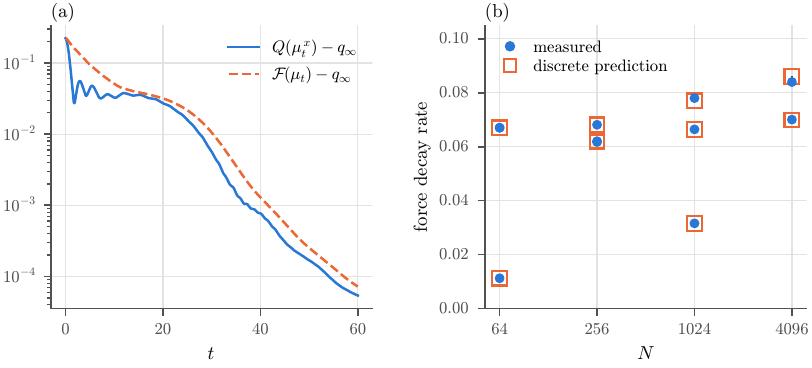}
	\caption{Teacher--student network.
	(a) The objective oscillates while the energy decreases, as in \cref{prop:dissipation}.
	(b) Force decay rates across widths and initializations, with predictions from the discrete linearization at each computed minimum.
	The spread reflects the different minima reached.}
	\label{fig:numerics-teacher}
\end{figure}

For the energy experiment, $\alpha<1$ and \cref{prop:dissipation} gives monotonicity of $\FF$.
\Cref{fig:numerics-teacher}(a) shows the predicted behavior: the objective increases during some oscillations, while the energy decreases.
Under step refinement, the residual in the integrated energy identity decreases at first order, consistent with \eqref{eq:discrete-adam}.

At the computed minima, the force Jacobian $\mathsf H$ is positive definite.
The corresponding local inequality has $\theta=1/2$ for each fixed particle count, so \cref{thm:kl-rates} predicts exponential convergence.
The local linearization \eqref{eq:numerics-linear-rate} gives the more precise rate $r=r(\lambda_{\min}(\mathsf H))$.
When a slow mode governs the motion, position and force norms have envelope rate $r$, while the objective gap, $\mathcal A_{1/2}$, and $\mathcal B$ have rate $2r$.
\Cref{tab:late-rate} compares the discrete linearized prediction with the observed decay.
The slowest eigenmode agrees with the prediction in every row, while the full force can have a different fitted rate when several modes contribute over the observation window.

\begin{table}[t]
	\centering
	\begin{tabular}{cc|ccc}
		$\rho_1$ & $\varepsilon$ & predicted & slowest mode & force\\ \hline
		$1$ & $10^{-1}$ & $0.00583$ & $0.00583$ & $0.00583$\\
		$1$ & $10^{-3}$ & $0.0612$ & $0.0609$ & $0.0612$\\
		$1$ & $10^{-5}$ & $0.503$ & $0.502$ & $0.486$\\
		$0.25$ & $10^{-3}$ & $0.0791$ & $0.0791$ & $0.0928$\\
		$4$ & $10^{-3}$ & $0.0562$ & $0.0556$ & $0.0562$
	\end{tabular}
	\caption{Exponential decay rates for the teacher--student model with $N=256$.
	The prediction uses the discrete linearization along the slowest eigenvector.
	The measured quantities are the displacement along that eigenvector and the $L^2(\Omega)$ norm of the force.}
	\label{tab:late-rate}
\end{table}

We also vary the width while keeping the effective regularization fixed.
\Cref{fig:numerics-teacher}(b) compares the measured force decay rates for several independent initializations at each width.
The rates show no systematic decrease over this range, although different initializations can reach minima with different curvatures.

\subsection{Numerical parameters and diagnostics}\label{sec:numerics-methods}

\Cref{tab:numerics-configurations} records the configurations for \eqref{eq:discrete-adam}.
Unless specified otherwise, memory starts at zero and scalar positions are independent standard Gaussian variables.
For the network, $u\sim N(0,1)$ and $w\sim N(0,I_8/8)$ independently.
The width comparison uses three initializations for $N\leq1024$ and two for $N=4096$.

\begin{table}[t]
	\centering
	\small
	\begin{tabular}{lccccc}
		Experiment & $N$ & $\rho_1$ & $\rho_2$ & $\varepsilon$ & $h$\\ \hline
		Interaction: corrected energy & $10^2$--$10^5$ & $1$ & $8$ & $10$ & $0.00125$\\
		Interaction: rates & $10$--$3000$ & $1$ & $0.5$ & $4$ & $0.05$\\
		Interaction: law comparison & $10^3$--$10^5$ & $1$ & $0.5$ & $4$ & $0.05$\\
		Quartic well & $1000$ & $1$ & $0.5$ & $10^{-2},\ 1$ & $0.05$\\
		Network: energy & $256$ & $0.2$ & $0.5$ & $10^{-3}$ & $0.005$--$0.04$\\
		Network: rates & $256$ & $0.25$--$4$ & $0.5$ & $10^{-5}$--$10^{-1}$ & $0.01$\\
		Network: width & $64$--$4096$ & $1$ & $0.5$ & $10^{-3}$ & $0.02$
	\end{tabular}
	\caption{Parameters for the particle discretization.
	The interacting model uses $a=2$ for the corrected energy, $a\in\{1,1.1,1.5\}$ for the rate comparison, and $a\in\{1,1.5\}$ for the law comparison.}
	\label{tab:numerics-configurations}
\end{table}

For the corrected-energy plot, the integral of $\mathcal R$ is evaluated by the trapezoidal rule up to $T=40$.
Truncating its upper limit subtracts the same constant from the corrected energy at every plotted time, and therefore does not affect its monotonicity.
The particle-law comparisons run to $T=2000$ and the quartic experiment to $T=10^5$.
For the single-well orbit, we instead use classical fourth-order Runge--Kutta integration to $T=6000$, with step
\begin{equation*}
	h_{\mathrm{RK}}=\min\{2\cdot10^{-3},\,0.02\varepsilon^{1/4}\}.
\end{equation*}
The reported amplitude is the maximum over the final fifth of the run; amplitudes below $10^{-6}$ are plotted as zero.
The step-refinement and initial-displacement checks use the same integrator.

Exponential force rates are least-squares slopes of the logarithm of the force norm over the final $60\%$ of the recorded window in which that norm lies between $10^{-13}$ and $10^{-7}$.
The slow-mode fits in \cref{tab:late-rate} use displacements from the final state over the wider force window from $10^{-4}$ to $10^{-13}$, fitting local maxima of the absolute displacement when the mode oscillates.
The network rate-table runs are continued until the mean-square force is at round-off level.
The width comparisons run to $T=400$, with early stopping at that level.
The network Jacobian is computed by central differences of the analytic force with step $10^{-6}$ and symmetrized, and at the largest width, its smallest eigenvalues are computed using directional differences.
All spectral predictions use the linearization of the discrete update, whose rate tends to \eqref{eq:numerics-linear-rate} as $h\to0$.

\bibliographystyle{amsplain}
\bibliography{references}

\providecommand{\bysame}{\leavevmode\hbox to3em{\hrulefill}\thinspace}
\providecommand{\MR}{\relax\ifhmode\unskip\space\fi MR }
\providecommand{\MRhref}[2]{%
  \href{http://www.ams.org/mathscinet-getitem?mr=#1}{#2}
}
\providecommand{\href}[2]{#2}
\begin{thebibliography}{10}

\bibitem{AmbrosioGigliSavare2008}
Luigi Ambrosio, Nicola Gigli, and Giuseppe Savar{\'e}, \emph{Gradient flows: In metric spaces and in the space of probability measures}, second ed., Lectures in Mathematics. ETH Z{\"u}rich, Birkh{\"a}user, Basel, 2008.

\bibitem{BarakatBianchi2021}
Anas Barakat and Pascal Bianchi, \emph{Convergence and dynamical behavior of the {ADAM} algorithm for nonconvex stochastic optimization}, SIAM Journal on Optimization \textbf{31} (2021), no.~1, 244--274.

\bibitem{BegoutBolteJendoubi2015}
Pascal B{\'e}gout, J{\'e}r{\^o}me Bolte, and Mohamed~Ali Jendoubi, \emph{On damped second-order gradient systems}, Journal of Differential Equations \textbf{259} (2015), no.~7, 3115--3143.

\bibitem{BelottoGazeau2020}
Andr{\'e} Belotto~da Silva and Maxime Gazeau, \emph{A general system of differential equations to model first-order adaptive algorithms}, Journal of Machine Learning Research \textbf{21} (2020), no.~129, 1--42.

\bibitem{CarrilloGvalaniWu2020}
Jos{\'e}~A. Carrillo, Rishabh~S. Gvalani, and Jeremy Wu, \emph{An invariance principle for gradient flows in the space of probability measures}, Journal of Differential Equations \textbf{345} (2023), 233--284.

\bibitem{ChenEtAl2019}
Xiangyi Chen, Sijia Liu, Ruoyu Sun, and Mingyi Hong, \emph{On the convergence of a class of {A}dam-type algorithms for non-convex optimization}, International Conference on Learning Representations, 2019.

\bibitem{ChizatBach2018}
L{\'e}na{\"i}c Chizat and Francis Bach, \emph{On the global convergence of gradient descent for over-parameterized models using optimal transport}, Advances in Neural Information Processing Systems, vol.~31, 2018.

\bibitem{deepmind2020jax}
DeepMind, Igor Babuschkin, Kate Baumli, Alison Bell, Surya Bhupatiraju, Jake Bruce, Peter Buchlovsky, David Budden, Trevor Cai, Aidan Clark, Ivo Danihelka, Antoine Dedieu, Claudio Fantacci, Jonathan Godwin, Chris Jones, Ross Hemsley, Tom Hennigan, Matteo Hessel, Shaobo Hou, Steven Kapturowski, Thomas Keck, Iurii Kemaev, Michael King, Markus Kunesch, Lena Martens, Hamza Merzic, Vladimir Mikulik, Tamara Norman, George Papamakarios, John Quan, Roman Ring, Francisco Ruiz, Alvaro Sanchez, Laurent Sartran, Rosalia Schneider, Eren Sezener, Stephen Spencer, Srivatsan Srinivasan, Milo\v{s} Stanojevi\'{c}, Wojciech Stokowiec, Luyu Wang, Guangyao Zhou, and Fabio Viola, \emph{The {D}eep{M}ind {JAX} {E}cosystem}, \url{https://github.com/google-deepmind}, 2020.

\bibitem{DefossezEtAl2022}
Alexandre D{\'e}fossez, L{\'e}on Bottou, Francis Bach, and Nicolas Usunier, \emph{A simple convergence proof of {A}dam and {A}dagrad}, Transactions on Machine Learning Research (2022), arXiv:2003.02395.

\bibitem{DuchiEtAl2011}
John Duchi, Elad Hazan, and Yoram Singer, \emph{Adaptive subgradient methods for online learning and stochastic optimization}, Journal of Machine Learning Research \textbf{12} (2011), no.~61, 2121--2159.

\bibitem{Dudley2002}
Richard~M. Dudley, \emph{Real analysis and probability}, second ed., Cambridge Studies in Advanced Mathematics, vol.~74, Cambridge University Press, Cambridge, 2002.

\bibitem{FeinbergKasyanovZadoianchuk2014}
Eugene~A. Feinberg, Pavlo~O. Kasyanov, and Nina~V. Zadoianchuk, \emph{{F}atou's lemma for weakly converging probabilities}, Theory of Probability and Its Applications \textbf{58} (2014), no.~4, 683--689.

\bibitem{HarauxJendoubi1998}
Alain Haraux and Mohamed~Ali Jendoubi, \emph{Convergence of solutions of second-order gradient-like systems with analytic nonlinearities}, Journal of Differential Equations \textbf{144} (1998), no.~2, 313--320.

\bibitem{HauerMazon2019}
Daniel Hauer and Jos{\'e}~M. Maz{\'o}n, \emph{{K}urdyka--{{\L}ojasiewicz}--{S}imon inequality for gradient flows in metric spaces}, Transactions of the American Mathematical Society \textbf{372} (2019), no.~7, 4917--4976.

\bibitem{KingmaBa2015}
Diederik~P. Kingma and Jimmy Ba, \emph{Adam: A method for stochastic optimization}, International Conference on Learning Representations, 2015, arXiv:1412.6980.

\bibitem{KricheneCaluyaHalder2020}
Walid Krichene, Kenneth~F. Caluya, and Abhishek Halder, \emph{Global convergence of second-order dynamics in two-layer neural networks}, 2020, arXiv:2007.06852.

\bibitem{Kurdyka1998}
Krzysztof Kurdyka, \emph{On gradients of functions definable in o-minimal structures}, Annales de l'Institut Fourier \textbf{48} (1998), no.~3, 769--783.

\bibitem{LaSalle1960}
Joseph~P. LaSalle, \emph{Some extensions of {L}iapunov's second method}, IRE Transactions on Circuit Theory \textbf{7} (1960), no.~4, 520--527.

\bibitem{MaWuE2022}
Chao Ma, Lei Wu, and Weinan E, \emph{A qualitative study of the dynamic behavior for adaptive gradient algorithms}, Proceedings of the 2nd Mathematical and Scientific Machine Learning Conference, Proceedings of Machine Learning Research, vol. 145, 2022, pp.~671--692.

\bibitem{MalladiEtAl2022}
Sadhika Malladi, Kaifeng Lyu, Abhishek Panigrahi, and Sanjeev Arora, \emph{On the {SDEs} and scaling rules for adaptive gradient algorithms}, Advances in Neural Information Processing Systems, vol.~35, 2022, pp.~7697--7711.

\bibitem{MeiMontanariNguyen2018}
Song Mei, Andrea Montanari, and Phan-Minh Nguyen, \emph{A mean field view of the landscape of two-layer neural networks}, Proceedings of the National Academy of Sciences \textbf{115} (2018), no.~33, E7665--E7671.

\bibitem{MeynTweedie1993}
Sean~P. Meyn and Richard~L. Tweedie, \emph{Stability of {M}arkovian processes {III}: {F}oster--{L}yapunov criteria for continuous-time processes}, Advances in Applied Probability \textbf{25} (1993), no.~3, 518--548.

\bibitem{Oksendal2003}
Bernt {\O}ksendal, \emph{Stochastic differential equations: An introduction with applications}, sixth ed., Universitext, Springer, Berlin, 2003.

\bibitem{Polyak1964}
Boris~T. Polyak, \emph{Some methods of speeding up the convergence of iteration methods}, USSR Computational Mathematics and Mathematical Physics \textbf{4} (1964), no.~5, 1--17.

\bibitem{ReddiEtAl2018}
Sashank~J. Reddi, Satyen Kale, and Sanjiv Kumar, \emph{On the convergence of {A}dam and beyond}, International Conference on Learning Representations, 2018.

\bibitem{RotskoffVandenEijnden2022}
Grant~M. Rotskoff and Eric Vanden-Eijnden, \emph{Trainability and accuracy of artificial neural networks: An interacting particle system approach}, Communications on Pure and Applied Mathematics \textbf{75} (2022), no.~9, 1889--1935.

\bibitem{SirignanoSpiliopoulos2020}
Justin Sirignano and Konstantinos Spiliopoulos, \emph{Mean field analysis of neural networks: A law of large numbers}, SIAM Journal on Applied Mathematics \textbf{80} (2020), no.~2, 725--752.

\bibitem{Sznitman1991}
Alain-Sol Sznitman, \emph{Topics in propagation of chaos}, {\'E}cole d'{\'E}t{\'e} de {P}robabilit{\'e}s de {S}aint-{F}lour {XIX}---1989, Lecture Notes in Mathematics, vol. 1464, Springer, Berlin, 1991, pp.~165--251.

\bibitem{Villani2009}
C{\'e}dric Villani, \emph{Hypocoercivity}, Memoirs of the American Mathematical Society \textbf{202} (2009), no.~950, iv+141.

\bibitem{Williams1997}
Christopher K.~I. Williams, \emph{Computing with infinite networks}, Advances in Neural Information Processing Systems (Michael~C. Mozer, Michael~I. Jordan, and Thomas Petsche, eds.), vol.~9, MIT Press, 1997, pp.~295--301.

\bibitem{WuKungurtsevMondelli2023}
Diyuan Wu, Vyacheslav Kungurtsev, and Marco Mondelli, \emph{Mean-field analysis for heavy ball methods: Dropout-stability, connectivity, and global convergence}, Transactions on Machine Learning Research (2023).

\bibitem{YangLittwin2023}
Greg Yang and Etai Littwin, \emph{Adaptive optimization in the $\infty$-width limit}, International Conference on Learning Representations, 2023.

\end{thebibliography}

\end{document}